\documentclass[11pt]{amsart}

\usepackage{amssymb, amsmath, amscd, eucal, rotating,mathtools}
\usepackage{tikz-cd}
\usetikzlibrary{decorations.markings}

\tikzset{double line with arrow/.style args={#1,#2}{decorate,decoration={markings,mark=at position 0 with {\coordinate (ta-base-1) at (0,1pt);
\coordinate (ta-base-2) at (0,-1pt);},
mark=at position 1 with {\draw[#1] (ta-base-1) -- (0,1pt);
\draw[#2] (ta-base-2) -- (0,-1pt);
}}}}

\input xy
 \xyoption{all}

\numberwithin{equation}{section}

\theoremstyle{plain}

\newtheorem{proposition}{Proposition}[section]
\newtheorem{theorem}[proposition]{Theorem}		
\newtheorem*{theorem*}{Theorem}		
\newtheorem{corollary}[proposition]{Corollary}
\newtheorem{lemma}[proposition]{Lemma}

\theoremstyle{definition}

\newtheorem{definition}[proposition]{Definition}
\newtheorem{remark}[proposition]{Remark}
\newtheorem{example}[proposition]{Example}

\theoremstyle{remark}

\newcommand{\C}{\mathbb{C}}

\DeclareMathOperator{\grad}{grad}
\DeclareMathOperator{\im}{im}
\DeclareMathOperator{\rank}{rank}
\DeclareMathOperator{\id}{id}
\DeclareMathOperator{\Gr}{Gr}
\DeclareMathOperator{\Rep}{Rep}
\DeclareMathOperator{\Hom}{Hom}
\DeclareMathOperator{\Vect}{Vect}
\DeclareMathOperator{\Rel}{Rel}
\DeclareMathOperator{\slope}{slope}

\DeclareMathOperator{\Crit}{Crit}

\DeclareMathOperator{\Ad}{Ad}

\DeclareMathOperator{\tr}{tr}
\DeclareMathOperator{\coker}{coker}

\renewcommand{\Im}{\mathsf{Im}}
\renewcommand{\Re}{\mathsf{Re}}

\DeclareMathOperator{\GL}{\mathsf{GL}}
\DeclareMathOperator{\U}{\mathsf{U}}
\DeclareMathOperator{\BU}{\mathsf{BU}}

\begin{document}


\title[Convolution via cup product on a Morse complex]{Convolution for quiver varieties via cup product on a Morse complex}

\author{Graeme Wilkin}
\address{Department of Mathematics,
University of York, 
York YO10 5DD, United Kingdom}
\email{graeme.wilkin@york.ac.uk}

\date{\today}

\thanks{This material is based upon work supported by the National Science Foundation under Grant No. 1440140, while the author was in residence at the Mathematical Sciences Research Institute in Berkeley, California, during the Fall 2022 program ``Analytic and Geometric Aspects of Gauge Theory''.}

\begin{abstract}
Convolution in Borel-Moore homology plays an important role in Nakajima's construction of representations of the Heisenberg algebra and of modified enveloping algebras of Kac-Moody algebras. In its most basic form, convolution between two quiver varieties is given by pullback and then pushforward via the Hecke correspondence for quivers. 

In previous work we showed that the Hecke correspondence has a Morse-theoretic interpretation in terms of spaces of flow lines. The goal of this paper is to show that the topological information that defines generators for Nakajima's representations can be encoded in the cup product for a Morse complex defined on the smooth space of representations of a quiver without relations, and then pulling back to the subvariety of representations that do satisfy a given set of relations. The results are valid for the main motivating example of Nakajima quivers, as well as other quivers with relations derived from these (for example handsaw quivers).

For the norm square of a moment map on the space of representations of a quiver, the usual Morse-Bott-Smale transversality condition on the space of flow lines fails, however a weaker version of transversality is still satisfied. A major part of the paper is spent developing a general theory in this setting of weak transversality from which one can recover the usual construction of the differentials and cup product on the Morse complex by adding an intermediate step of taking cup product with a certain Euler class, which is explicitly computable for the space of representations of a quiver. 
\end{abstract}




\maketitle


\thispagestyle{empty}

\baselineskip=16pt



\section{Introduction}

A well known tool to compute the cohomology of a space $X$ is to begin with a filtration $\emptyset = X_{-1} \subset X_0 \subset \cdots \subset X_n = X$ and then use a spectral sequence to compute the cohomology. The Morse-theoretic point of view is to approach this problem by considering a function $f : X \rightarrow \mathbb{R}$ together with an associated gradient flow or pseudogradient flow $\varphi : X \times \mathbb{R} \rightarrow X$. If the flow has good compactness and convergence properties, then there is a canonical \emph{Morse filtration} of $X$ associated to $f$.  The \emph{main theorem of Morse theory} says that, if the flow also has good local behaviour around the critical points of $X$ (for example $f$ is Morse or Morse-Bott), then the groups on the first page of the spectral sequence can be expressed in terms of  the cohomology of the critical sets together with their Morse indices. 

If, in addition to the properties above, the flow has well-behaved spaces of flow lines between critical sets (for example, if the function is Morse-Bott-Smale so that the stable and unstable manifolds always intersect transversely) then the differentials and the cup product in the spectral sequence can be computed using pullback/pushforward homomorphisms via  spaces of flow lines. The benefit of using Morse theory is that all of the data used above (the critical sets, their indices and spaces of flow lines between critical points) can be described analytically using the function $f$ and the flow $\varphi$, thus giving us a concrete way to compute the terms in the spectral sequence.


There are many interesting examples that do not immediately satisfy the above conditions. For example, transversality may fail, or (even worse) the space may be singular. The main motivating examples for this paper are the space of representations $\Rep(Q, {\bf v})$ of a quiver $Q$ with a fixed dimension vector ${\bf v}$ as well as the subvarieties $\Rep(Q, {\bf v}, \mathcal{R})$ of representations satisfying a given set $\mathcal{R}$ of relations. The space $\Rep(Q, {\bf v})$ is a vector space with a linear action of a complex reductive group $G$ and a symplectic structure for which the action of the maximal compact subgroup $K$ is Hamiltonian, from which we obtain a moment map $\mu : \Rep(Q, {\bf v}) \rightarrow \mathfrak{k}^*$ defined up to a constant central element $\alpha \in Z(\mathfrak{k}^*)$. This setup fits into a large class of examples for which the Morse theory of $\| \mu - \alpha \|^2$ has been well-studied (cf. \cite{AtiyahBott83}, \cite{Kirwan84}) and it turns out that the associated stratification is equivariantly perfect with respect to the action of $K$, which leads to very explicit inductive formulae for the cohomology of the quotient space \cite{Kirwan84}.

The results of \cite{Wilkin17} completely classify the flow lines of $\| \mu - \alpha \|^2$ on the space $\Rep(Q, {\bf v})$ and show that the space of flow lines is related to the Hecke correspondence. It is then natural to ask whether the Morse-Bott-Smale theory can be extended to this setting and whether the pullback/pushforward homomorphisms via the Hecke correspondence that appear in Nakajima's work \cite{Nakajima94}, \cite{Nakajima97}, \cite{Nakajima98} also appear when constructing the cup product on the spectral sequence. The goal of this paper is to prove that this is indeed the case (cf. Theorem \ref{thm:intro-cup-product-induces-convolution}). 

In Section \ref{subsec:flow-lines} we show that transversality fails in this setting, however a weaker form of transversality is still satisfied, which is defined as follows. Given a Riemannian manifold $M$ and a smooth function $f: M \rightarrow \mathbb{R}$, one can define the downwards gradient flow $\varphi_t(x_0)$ satisfying $\frac{\partial}{\partial t} \varphi_t(x_0) = -\nabla f(\varphi_t(x_0))$ and $\varphi_0(x_0) = x_0$. If the flow has good properties such as those in Conditions (1)--(3) of \cite{Wilkin19} ($M$ is real analytic, $f$ is analytic with isolated critical values and the flow satisfies a compactness condition) then for each pair of critical sets $C_\ell$ and $C_u$ with $f(C_\ell) < f(C_u)$, we can define the \emph{stable and unstable sets}
\begin{align*}
W_{C_\ell}^+ & := \{ x \in M \, \mid \, \lim_{t \rightarrow \infty} \varphi_t(x) \in C_\ell \} \\
W_{C_u}^- & := \{ x \in M \, \mid \, \lim_{t \rightarrow - \infty} \varphi_t(x) \in C_u \}, \quad \text{and} \quad W_{C_u,0}^- := W_{C_u}^- \setminus C_u .
\end{align*}
and then the \emph{space of points lying on a flow line}
\begin{equation*}
\mathcal{F}_{\ell,0}^{u,0} := W_{C_\ell}^+ \cap W_{C_u}^- .
\end{equation*}
Then the weak transversality condition is defined as follows.

\begin{definition}[cf. Definition \ref{def:transversality-conditions}]\label{def:intro-transversality-conditions}

Let $M$ be a Riemannian manifold and $f : M \rightarrow \mathbb{R}$ a minimally degenerate smooth function satisfying Conditions (1)--(3) of \cite{Wilkin19}. The spaces of flow lines satisfies \emph{weak transversality} if the following conditions hold. 

\begin{enumerate}

\item[{\bf (T1)}] The space of flow lines $\mathcal{F}_{\ell,0}^{u,0}$ has a tubular neighbourhood in $W_{C_u,0}^-$, denoted $D_\ell^u \rightarrow \mathcal{F}_{\ell,0}^{u,0}$. 

\item[{\bf (T2)}] The stratum $W_{C_\ell}^+$ has a tubular neighbourhood in $M$ denoted $\tilde{V}_\ell \rightarrow W_{C_\ell}^+$, which restricts to a disk bundle $V_\ell \rightarrow \mathcal{F}_{\ell,0}^{u,0}$ such that $D_\ell^u$ is a subbundle of $V_\ell$.

\end{enumerate}

\end{definition}

\begin{remark}
The usual Morse-Bott-Smale transversality condition is equivalent to the condition that $D_\ell^u = V_\ell$.\end{remark}

Proposition \ref{prop:quivers-T2} shows that $\| \mu - \alpha \|^2 : \Rep(Q, {\bf v}) \rightarrow \mathbb{R}$ satisfies weak transverality, and therefore the definition is satisfied for a large class of interesting examples. The key to proving this is the local analysis of \cite{Wilkin17}, which shows that analytic questions about the structure of the space of flow lines near the upper critical set can be reduced (up to homeomorphism) to algebraic questions on a linearisation of the unstable set, which we call the \emph{negative slice}. Moreover, this homeomorphism restricts to a homeomorphism in any $G$-invariant subset \cite[Cor. 4.24]{Wilkin17}. The space of flow lines then corresponds to a subvariety of the negative slice with an explicit description (cf. Section \ref{subsec:flow-lines-neg-slice}).

\begin{remark}
There are many versions of Morse theory where one perturbs the function to obtain a well-behaved Morse-Smale function. We want to avoid this here, since the unperturbed function $\| \mu - \alpha \|^2 : \Rep(Q, {\bf v}) \rightarrow \mathbb{R}$ contains a lot of useful topological information; for example the critical sets correspond to quiver varieties of smaller dimension and the flow lines are related to Nakajima's Hecke correspondence \cite{Wilkin17}. Moreover, the same is true after restricting to a subvariety of representations that satisfy a given set of relations, in which case perturbing the function is more difficult since lack of smoothness means that arbitrary perturbations may create a function whose gradient flow does not preserve the subvariety. 

The main point of this paper is to show that the topological information that defines generators for Nakajima's representations from \cite{Nakajima98} is encoded in the cup product on the Morse complex for the unperturbed function $\| \mu - \alpha \|^2$, which is why we work with weak transversality here, rather than perturbing the function.
\end{remark}

For Morse-Bott-Smale functions, the differentials and cup product can be constructed using pullback/pushforward homomorphisms via spaces of flow lines connecting critical sets (cf. \cite{AustinBraam95}). In the setting of weak transversality, one can still use spaces of flow lines to define differentials and cup product on the Morse complex, however now we need to include an additional step of taking cup product with the Euler class of the quotient bundle $V_\ell / D_\ell^u$  (cf. Section \ref{subsec:cup-product-T1-T2}).

This theory works on the smooth space $\Rep(Q, {\bf v})$, however from the point of view of representation theory (cf. \cite{Nakajima94}, \cite{Nakajima97}, \cite{Nakajima98}, \cite{Nakajima01}, \cite{Nakajima04}) it is more interesting to study the subset of representations $\Rep(Q, {\bf v}, \mathcal{R}) \subset \Rep(Q, {\bf v})$ that satisfy a given set $\mathcal{R}$ of relations. The most important example is the class of Nakajima quiver varieties, where the relations take the form of a complex moment map, however the results below are not restricted to this class of examples.

The gradient flow of $\| \mu - \alpha \|^2$ on $\Rep(Q, {\bf v})$ is generated by the complex reductive group $G$, and therefore preserves the subset $\Rep(Q, {\bf v}, \mathcal{R})$. From this we can construct an injective map from each critical set in $\Rep(Q, {\bf v}, \mathcal{R})$ into the corresponding critical set in $\Rep(Q, {\bf v})$, and so the cup product homomorphism defined above pulls back to a homomorphism between the cohomology of critical sets on the singular space $\Rep(Q, {\bf v}, \mathcal{R})$. The Kirwan surjectivity theorem of McGerty and Nevins \cite{McGertyNevins18} shows that, for Nakajima quivers, \emph{all} of the classes of $H^*(\mathcal{M}(Q, {\bf v_\ell}, \mathcal{R}))$ appear in the image of the pullback $H^*(\mathcal{M}(Q, {\bf v_\ell})) \rightarrow H^*(\mathcal{M}(Q, {\bf v_\ell}, \mathcal{R}))$ for each ${\bf v_\ell} \leq {\bf v}$.

In order to construct a convolution homomorphism between critical sets on the singular space, we need to take cup product with a well-chosen cohomology class on the ambient smooth space. This turns out to be the cup product of a Thom class $\tilde{\tau}$ of an explicitly defined submanifold of the space of flow lines (cf. Corollary \ref{cor:thom-compatible} and Corollary \ref{cor:global-generation-Thom}) with a Chern class $\tilde{\xi}$ (cf. Lemma \ref{lem:global-generation-pushforward}). The main result of this paper is Theorem \ref{thm:cup-product-induces-convolution}, which shows that convolution is induced by cup product.

\begin{theorem}[Theorem \ref{thm:cup-product-induces-convolution}] \label{thm:intro-cup-product-induces-convolution}
Let $Q$ be a quiver with a set $\mathcal{R}$ of complete quadratic relations that are fully restricted from those for a Nakajima quiver (cf. Definitions \ref{def:complete-relations} and \ref{def:fully-restricted}). Consider a pair of adjacent critical sets $C_u$ (upper) and $C_\ell$ (lower) for the moment map energy function $\| \mu - \alpha \|^2 : \Rep(Q, {\bf v}, \mathcal{R}) \rightarrow \mathbb{R}$ corresponding to quiver varieties $\mathcal{M}(Q, {\bf v_u}, \mathcal{R})$ and $\mathcal{M}(Q, {\bf v_\ell}, \mathcal{R})$ with dimension vectors ${\bf v_u} < {\bf v_\ell} = {\bf v_u} + {\bf e_k} \leq {\bf v}$.

Then the Poincar\'e dual of convolution in Borel-Moore homology via the Hecke correspondence 
\begin{equation*}
H^{BM}_*(\mathcal{M}(Q, {\bf v_\ell}, \mathcal{R})) \stackrel{\pi_\ell^*}{\longrightarrow} H^{BM}_*(\mathcal{B}(Q, {\bf v_u}, {\bf v_\ell}, \mathcal{R})) \stackrel{(\pi_u)_*}{\longrightarrow} H^{BM}_*(\mathcal{M}(Q, {\bf v_u}, \mathcal{R}))
\end{equation*}
is induced by cup product with $\tilde{\tau} \smallsmile \tilde{\xi}$ on the Morse complex for the moment map energy function on the ambient smooth space $\Rep(Q, {\bf v}) \times \Gr({\bf v_u}, {\bf v_\ell})$.
\end{theorem}

{\bf Organisation of the paper.} Section \ref{sec:background} gives a self-contained overview of the necessary definitions and constructions needed to develop the Morse theory on the space of representations of a quiver in later sections. Section \ref{sec:restricted-quivers} contains a general construction that can be used to show that the moduli space and negative slice are smooth for a large class of quivers with relations, which we can then use to set up the deformation theory of these spaces in Section \ref{sec:convolution-cup-product}. Section \ref{sec:Morse-spectral} then recalls how to construct the differentials and cup product on the first page of the spectral sequence associated to a Morse filtration. The construction here only uses the existence of a Morse filtration and is written in a way that leads to Section \ref{sec:main-theorem-localisation}, which uses the main theorem of Morse theory to further refine the construction of the cup product and differentials so that it is expressed in terms of relative cohomology groups localised around the critical points and spaces of flow lines. We focus on the case of adjacent critical sets (which is sufficient to construct generators for the representations of \cite{Nakajima98}), as the general case requires compactifying spaces of flow lines which we defer to a later paper.

Up until this point all the results are valid for functions $f : Z \rightarrow \mathbb{R}$ satisfying Conditions (1)--(5) of \cite{Wilkin19} (in particular, we do not require the space $Z$ to be smooth). In Section \ref{sec:transversality} we now restrict to the case where the space is a manifold and the function satisfies Kirwan's minimal degeneracy condition (cf. \cite{Kirwan84}) as well as the weak transversality conditions of Definition \ref{def:transversality-conditions}. Now the relative cohomology groups of Sections \ref{sec:Morse-spectral} and \ref{sec:main-theorem-localisation} can be rewritten as cohomology groups of the critical sets and spaces of flow lines, and the differentials and cup product can be expressed in terms of pullback and cup product maps between these spaces (cf. the diagrams \eqref{eqn:differential-T1-T2} and \eqref{eqn:cup-product-T1-T2}). 

In Section \ref{sec:BM-convolution} we recall some facts about Borel-Moore homology before showing in Section \ref{sec:morse-convolution} that the homomorphisms used to construct the differentials and cup product are Poincar\'e dual to pullback and pushforward in Borel-Moore homology. Section \ref{sec:convolution-cup-product} then contains the main result (Theorem \ref{thm:cup-product-induces-convolution}), which shows that cup product with the Thom class of a certain submanifold then pulls back to a homomorphism which is Poincar\'e dual to convolution in Borel-Moore homology on the space of representations satisfying a given set of relations.

\section{Background}\label{sec:background}

In this section we recall the important definitions and theorems from \cite{Wilkin17} and \cite{Wilkin19} which will be used in the rest of the paper. The original sources used for the material on quivers are \cite{King94}, \cite{Nakajima94}, \cite{Nakajima98} and \cite{Crawley01}. The notation and setup follows that used in \cite{Wilkin17}.

\subsection{The space of representations of a quiver}


A \emph{quiver} is a directed graph, consisting of a set of vertices $\mathcal{I}$, edges $\mathcal{E}$ and head/tail maps $h,t : \mathcal{E} \rightarrow \mathcal{I}$. The quiver is \emph{finite} if the sets $\mathcal{I}$ and $\mathcal{E}$ are finite. A \emph{complex representation of a quiver $Q$} consists of a collection of complex vector spaces $\{ V_k \}_{k \in \mathcal{I}}$ and $\C$-linear homomorphisms $\{ x_a : V_{t(a)} \rightarrow V_{h(a)}\}_{a \in \mathcal{E}}$. The \emph{dimension vector} of a representation is the vector ${\bf v} := \left( \dim_\C V_k \right)_{k \in \mathcal{I}} \in \mathbb{Z}_{\geq 0}^\mathcal{I}$.

The vector space of all complex representations of $Q$ with fixed dimension vector ${\bf v}$ is denoted
\begin{equation}\label{eqn:rep-space-def}
\Rep(Q, {\bf v}) := \bigoplus_{a \in \mathcal{E}} \Hom(V_{t(a)}, V_{h(a)}) .
\end{equation}
The group
\begin{equation}\label{eqn:reductive-group-def}
G_{\bf v} := \prod_{k \in \mathcal{I}} \GL(V_k, \C)
\end{equation}
acts on the space $\Vect(Q, {\bf v}) := \bigoplus_{k\in \mathcal{I}} V_k$ and therefore on $\Rep(Q, {\bf v})$ via the induced action on each vector $\Hom(V_{t(a)}, V_{h(a)})$
\begin{equation}\label{eqn:group-action-def}
(g_i)_{k \in \mathcal{I}} \cdot (x_a)_{a \in \mathcal{E}} := \left( g_{h(a)} x_a g_{t(a)}^{-1} \right)_{a \in \mathcal{E}} .
\end{equation}
Given a Hermitian structure on each vector space $V_k$, we can define the unitary group $\U(V_k)$ and therefore define
\begin{equation*}
K_{\bf v} := \prod_{k \in \mathcal{I}} \U(V_k) .
\end{equation*}
This acts on $\Rep(Q, {\bf v})$ via the inclusion $K_{\bf v} \hookrightarrow G_{\bf v}$. The Lie algebras of $G_{\bf v}$ and $K_{\bf v}$ are denoted $\mathfrak{g}_{\bf v}$ and $\mathfrak{k}_{\bf v}$ respectively. Given a representation $x \in \Rep(Q, {\bf v})$, the infinitesimal action of $\mathfrak{g}_{\bf v}$ on the tangent space $T_x \Rep(Q, {\bf v})$ is denoted $\rho_x^\C : \mathfrak{g}_{\bf v} \rightarrow T_x \Rep(Q, {\bf v})$ and given by the following formula
\begin{equation}\label{eqn:complex-inf-action}
\rho_x^\C(u) := \left. \frac{d}{dt} \right|_{t=0} e^{tu} \cdot x = \bigoplus_{a \in \mathcal{E}} (u_{h(a)} x_a - x_a u_{t(a)} ) \in \bigoplus_{a \in \mathcal{E}} \Hom(V_{t(a)}, V_{h(a)}) \cong \Rep(Q, {\bf v})
\end{equation}
The action of $\mathfrak{k}_{\bf v}$ is denoted $\rho_x : \mathfrak{k}_{\bf v} \rightarrow \Rep(Q, {\bf v})$ and is defined via the inclusion $\mathfrak{k}_{\bf v} \hookrightarrow \mathfrak{g}_{\bf v}$. The distinction between $\rho_x^\C$ and $\rho_x$ will be important when using the adjoint of these homomorphisms; for example when defining the local slice in Definition \ref{def:slice}.

Given a quiver $Q$, let ${\bf v_1} = (v_k^{(1)})_{k \in \mathcal{I}} \in \mathbb{Z}^\mathcal{I}$ and ${\bf v_2} = (v_k^{(2)})_{k \in \mathcal{I}} \in \mathbb{Z}^\mathcal{I}$ be dimension vectors for $Q$. For each $k \in \mathcal{I}$, let $V_k^{(1)}$ and $V_k^{(2)}$ denote complex vector spaces of dimension $v_k^{(1)}$ and $v_k^{(2)}$ respectively. We define the spaces
\begin{equation*}
\Hom^0(Q, {\bf v_1}, {\bf v_2}) := \bigoplus_{k \in \mathcal{I}} \Hom(V_k^{(1)}, V_k^{(2)}) \quad \text{and} \quad \Hom^1(Q, {\bf v_1}, {\bf v_2}) := \bigoplus_{a \in \mathcal{E}} \Hom(V_{t(a)}, V_{h(a)}) .
\end{equation*}
These are spaces of homomorphisms along length $0$ and length $1$ paths in $Q$. Note that $\mathfrak{g}_{\bf v} \cong \Hom^0(Q, {\bf v}, {\bf v})$. The dimension of $\Hom^0(Q, {\bf v_1}, {\bf v_2})$ is denoted by
\begin{equation*}
{\bf v_1} \cdot {\bf v_2} := \dim_\C \Hom^0(Q, {\bf v_1}, {\bf v_2}) = \sum_{k \in \mathcal{I}} v_k^{(1)} v_k^{(2)} .
\end{equation*}
The \emph{Ringel form} (cf. \cite[Sec. 2]{Crawley01}) on $\mathbb{Z}^\mathcal{I}$ is
\begin{equation}\label{eqn:ringel-def}
\left< {\bf v_1}, {\bf v_2} \right> := \dim_\C \Hom^0(Q, {\bf v_1}, {\bf v_2}) - \dim_\C \Hom^1(Q, {\bf v_1}, {\bf v_2}) = \sum_{k \in \mathcal{I}} v_k^{(1)} v_k^{(2)} - \sum_{a \in \mathcal{E}} v_{t(a)}^{(1)} v_{h(a)}^{(2)} .
\end{equation}
The Ringel form leads to a simple expression for the Euler characteristic of the deformation complex \eqref{eqn:neg-def-complex}, which we use in \eqref{eqn:neg-dimension}.

\subsection{Representations of a quiver with relations}

A \emph{path} in $Q$ is a concatenation of a finite number of edges, denoted $p = a_n \cdots a_1$ with $h(a_j) = t(a_{j+1})$ for each $j = 1, \ldots, n-1$. The head and tail of the path $p = a_n \cdots a_1$ are $h(p) := h(a_n)$ and $t(p) := t(a_1)$. A path $p = a_n \cdots a_1$ determines a homomorphism $\Rep(Q, {\bf v}) \rightarrow \Hom(V_{t(p)}, V_{h(p)})$ given by 
\begin{equation*}
x \mapsto x_p := x_{a_n} \cdots x_{a_1} .
\end{equation*}
Given a representation $x = x_1 \oplus x_2 \in \Rep(Q, {\bf v_1}) \oplus \Rep(Q, {\bf v_2})$ with vector spaces denoted by $\Vect(Q, {\bf v_j}) = \bigoplus_{k \in \mathcal{I}} V_k^{(j)}$ for $j = 1,2$, a path $p = a_n \cdots a_1$ determines a homomorphism $dp_x : \Hom^1(Q, {\bf v_2}, {\bf v_1}) \rightarrow \Hom(V_{t(p)}^{(2)}, V_{h(p)}^{(1)})$ defined as follows. For each $a \in \mathcal{E}$ and $j=1,2$, let $(x_j)_a : V_{t(a)}^{(j)} \rightarrow V_{h(a)}^{(j)}$ denote the homomorphisms in the representations $x_1$ and $x_2$, and for any $\delta x \in \Hom^1(Q, {\bf v}_2, {\bf v}_1)$ define
\begin{equation}\label{eqn:path-derivative-def}
dp_x(\delta x) := \sum_{\ell=1}^n (x_1)_{a_n} \cdots (x_1)_{a_{\ell+1}} (\delta x)_{a_\ell} (x_2)_{a_{\ell-1}} \cdots (x_2)_{a_1} .
\end{equation}

Let $\mathcal{P}_{t,h}$ denote the set of all paths with a given tail $t \in \mathcal{I}$ and head $h \in \mathcal{I}$. Given a quiver $Q$, a \emph{relation $r$ with head $h(r) \in \mathcal{I}$ and tail $t(r) \in \mathcal{I}$} is a finite $\C$-linear combination of paths, denoted by 
\begin{equation*}
r = \sum_{p \in \mathcal{P}_{t(r), h(r)}} \lambda_p p .
\end{equation*}
A \emph{quadratic relation} is a relation for which the only paths $p$ with $\lambda_p \neq 0$ have length two. Any relation determines an algebraic map $\nu_r : \Rep(Q, {\bf v}) \rightarrow \Hom(V_{t(r)}, V_{h(r)})$ given by 
\begin{equation*}
x \mapsto \sum_{p} \lambda_p x_p .
\end{equation*}
Given a finite set of relations $\mathcal{R}$, define $\nu : \Rep(Q, {\bf v}) \rightarrow \bigoplus_{r \in \mathcal{R}} \Hom(V_{t(r)}, V_{h(r)})$ by
\begin{equation}\label{eqn:relation-map-def}
\nu(x) = \bigoplus_{r \in \mathcal{R}} \nu_r(x) .
\end{equation}
A \emph{representation of $(Q, \mathcal{R})$} is a representation in the kernel of $\nu$. The \emph{space of all representations of $(Q, \mathcal{R})$ with dimension vector ${\bf v}$} is the subset $\Rep(Q, {\bf v}, \mathcal{R}) := \nu^{-1}(0) \subset \Rep(Q, {\bf v})$. For example, quivers with relations are used to construct hyperk\"ahler ALE 4 manifolds \cite{Kronheimer89}, moduli spaces of instantons on hyperk\"ahler ALE 4 manifolds \cite{KronheimerNakajima90}, and more generally the hyperk\"ahler quotients of \cite{Nakajima94}, \cite{Nakajima98}, \cite{Nakajima01}, \cite{Nakajima04}. The handsaw quivers of \cite{Nakajima12} are another important example of a quiver with relations. It will be useful to note that all of these examples have quadratic relations.

With a view towards studying the deformation complex \eqref{eqn:neg-def-complex}, we define
\begin{equation}\label{eqn:relation-hom-def}
\Rel(Q, {\bf v_2}, {\bf v_1}, \mathcal{R}) := \bigoplus_{r \in \mathcal{R}} \Hom(V_{t(r)}^{(2)}, V_{h(r)}^{(1)})
\end{equation}
for a finite set of relations $\mathcal{R}$. For each $r \in \mathcal{R}$, write $r = \sum_{p\in \mathcal{P}_{t(r), h(r)}} \lambda_p p$. Given a representation $x = x_1 \oplus x_2 \in \Rep(Q, {\bf v_1}, \mathcal{R}) \oplus \Rep(Q, {\bf v_2}, \mathcal{R})$, we can extend the homomorphism \eqref{eqn:path-derivative-def} to a homomorphism $d \nu_x : \Hom^1(Q, {\bf v_2}, {\bf v_1}) \rightarrow \Rel(Q, {\bf v_2}, {\bf v_1}, \mathcal{R})$ given by 
\begin{equation}\label{eqn:relation-derivative-def}
d\nu_x (\delta x) := \sum_{r \in \mathcal{R}} \sum_{p \in \mathcal{P}_{t(r),h(r)}} \lambda_p(r) dp_x(\delta x) .
\end{equation}
The usual deformation theory shows that a point in the moduli space \eqref{eqn:symplectic-quotient-def} is smooth if $\rho_x^\mathbb{C}$ is injective and $d \nu_x$ is surjective in the following deformation complex
\begin{equation}\label{eqn:deformation-tangent-space}
\begin{tikzcd}
\Hom^0(Q, {\bf v}, {\bf v}) \arrow{r}{\rho_x^\mathbb{C}} & \Hom^1(Q, {\bf v}, {\bf v}) \arrow{r}{d \nu_x} & \Rel(Q, {\bf v}, {\bf v}, \mathcal{R}) ,
\end{tikzcd}
\end{equation}
in which case the tangent space to the moduli space is given by the middle cohomology of the above complex.

The following lemma is used in Section \ref{subsec:canonical-central-element} to show that the negative slice can be defined using the deformation complex \eqref{eqn:neg-def-complex}.

\begin{lemma}\label{lem:linearise-slice}
Let $\mathcal{R}$ be a finite set of relations, let $x = x_1 \oplus x_2 \in  \Rep(Q, {\bf v_1}) \oplus \Rep(Q, {\bf v_2})$ and let $\delta x \in \Hom^1(Q, {\bf v_2}, {\bf v_1})$. Then $\nu(x) = \nu(x + \delta x)$ if and only if $\delta x \in \ker d\nu_x$.
\end{lemma}

\begin{proof}
With respect to the decomposition $\Rep(Q, {\bf v_1}, \mathcal{R}) \oplus \Rep(Q, {\bf v_2}, \mathcal{R})$, for each $k \in \mathcal{I}$ let $V_k^{(1)}$ and $V_k^{(2)}$ denote the vector spaces at vertex $k$.

First note that for any pair of edges $a_1, a_2 \in \mathcal{E}$ with $t(a_2) = h(a_1)$ we have
\begin{equation*}
(x_2)_{a_2} (\delta x)_{a_1} = 0, \quad (\delta x)_{a_2} (x_1)_{a_1} = 0, \quad (\delta x)_{a_2} (\delta x)_{a_1} = 0 .
\end{equation*}
The first equation above follows from the fact that the domain of $(x_1)_{a_2}$ is $V_{t(a_2)}^{(2)}$ and the image of $(\delta x)_{a_1}$ is contained in $V_{h(a_1)}^{(1)} = V_{t(a_2)}^{(1)}$, so the composition is zero. The other equations follow from similar reasoning. Therefore for each path $p = a_n \cdots a_1$ we have
\begin{align*}
(x + \delta x)_p & = (x_1)_{a_n} \cdots (x_1)_{a_1} + (x_2)_{a_n} \cdots (x_2)_{a_1} + \sum_{\ell=1}^n (x_1)_{a_n} \cdots (x_1)_{a_{\ell+1}} (\delta x)_{a_\ell} (x_2)_{a_{\ell-1}} \cdots (x_2)_{a_1} \\
 & = (x_1)_p + dp_x(\delta x) + (x_2)_p .
\end{align*}
Summing the above formula for each relation $r = \sum_p \lambda_p p \in \mathcal{R}$ gives us
\begin{equation*}
\nu_r(x+\delta x) = \sum_p \lambda_p (x + \delta x)_p = \sum_p \lambda_p \left( (x_2)_p + (x_1)_p \right) + \sum_p \lambda_p dp_x(\delta x) 
\end{equation*}
and so summing over all relations $r \in \mathcal{R}$ leads to
\begin{equation*}
\nu(x + \delta x) = \nu(x) + d\nu_x(\delta x) .
\end{equation*}
Therefore $\nu(x) = \nu(x + \delta x)$ if and only if $d \nu_x(\delta x) = 0$. 
\end{proof}


Now consider the case where the relations are all quadratic. Each relation then has the form
\begin{equation*}
r = \sum_{p \in \mathcal{P}_{t(r), h(r)}} \lambda_{p}(r) p, \quad p = a_2(p) a_1(p) .
\end{equation*}
\begin{definition}\label{def:complete-relations}
A set $\mathcal{R}$ of quadratic relations is \emph{complete} if and only if

\begin{enumerate}

\item for each edge $a \in \mathcal{E}$ and every relation $r \in \mathcal{R}$ such that $h(a) = h(r)$, there exists a path $p \in \mathcal{P}_{t(r), h(r)}$ such that $\lambda_p(r) \neq 0$ and $a = a_2(p)$, and

\item for each edge $a \in \mathcal{E}$ such that $t(a) = t(r)$ for some relation $r \in \mathcal{R}$, then $r$ is the unique relation with this property, and there exists a unique path $p \in \mathcal{P}_{t(r), h(r)}$ such that $\lambda_p(r) \neq 0$ and $a = a_1(p)$.

\end{enumerate}

\end{definition}

Completeness is used in Lemma \ref{lem:cokernel-dimension}, which gives a formula to compute the cokernel of the homomorphism $d\nu_x$ used in the deformation complex \eqref{eqn:neg-def-complex}. This formula is much easier to use in examples, since it only depends on the dimension of the image of a representation. The following examples show that completeness occurs for many examples of interest; in particular all of the quivers from \cite{Kronheimer89}, \cite{KronheimerNakajima90}, \cite{Nakajima94}, \cite{Nakajima98}, \cite{Nakajima01}, \cite{Nakajima04} and \cite{Nakajima12} have complete sets of quadratic relations.

\begin{example}[Nakajima quivers have complete relations]\label{ex:nakajima-complete}
A \emph{Nakajima quiver} is a quiver $Q$ with vertices $\mathcal{I}$, edges $\mathcal{E}$ and head/tail maps $h,t: \mathcal{E} \rightarrow \mathcal{I}$ such that each edge $a \in \mathcal{E}$ has a conjugate $\bar{a} \neq a$ such that $h(\bar{a}) = t(a)$, $t(\bar{a}) = h(a)$ and $\bar{\bar{a}} = a$. Choose a subset of edges $\mathcal{E}^{0,1} \subset \mathcal{E}$ such that $\mathcal{E} = \mathcal{E}^{0,1} \cup \overline{\mathcal{E}^{0,1}}$ and $\mathcal{E}^{0,1} \cap \overline{\mathcal{E}^{0,1}}$ is empty. For each vertex $k \in \mathcal{I}$ there is a single relation 
\begin{equation}\label{eqn:hyperkahler-relation}
r_k = \sum_{a \in \mathcal{E}^{0,1} \, : \, h(a)=k} a \bar{a} - \sum_{a \in \mathcal{E}^{0,1} \, : \, t(a) = k} \bar{a} a .
\end{equation}
Therefore $t(r_k) = h(r_k)$ for each relation $r_k$. Let $\mathcal{R}$ denote the finite set of relations $\{ r_k \}_{k \in \mathcal{I}}$. From \eqref{eqn:hyperkahler-relation} it is clear that for each edge $a \in \mathcal{E}^{0,1}$ there is a unique relation $r = r_{h(a)}$ such that $a = a_2(p)$ for exactly one path $p$ with $\lambda_p(r) \neq 0$. The same reasoning applies to each edge $\bar{a} \in \overline{\mathcal{E}^{0,1}}$ since $h(\bar{a}) = t(a)$. Therefore the first condition of Definition \ref{def:complete-relations} is satisfied. By explicitly examining \eqref{eqn:hyperkahler-relation} again, it is clear that the second condition of Definition \ref{def:complete-relations} is also satisfied, and so the set of relations is complete.
\end{example}

\begin{example}[Handsaw quivers have complete relations]\label{ex:handsaw-complete}
Let $Q$ be a ``handsaw'' quiver as in \cite{Nakajima12} with edges labeled as below.
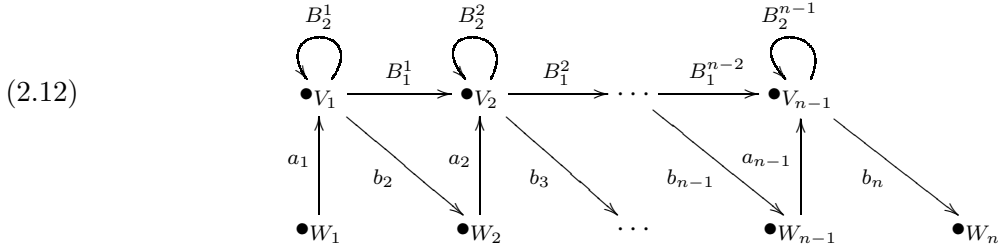
\begin{equation}\label{eqn:handsaw-quiver}
\xymatrixrowsep{0.5in}
\xymatrixcolsep{0.5in}
\xymatrix{
\bullet_{V_1} \ar[r]^{B_1^1} \ar[dr]_{b_2} \ar@`{(10,10),(-10,10)}_{B_2^1} & \bullet_{V_2} \ar[r]^{B_1^2} \ar[dr]_{b_3} \ar@`{(30,10),(10,10)}_{B_2^2} & \cdots \ar[r]^{B_1^{n-2}} \ar[dr]_{b_{n-1}} & \bullet_{V_{n-1}} \ar[dr]_{b_n} \ar@`{(72,10),(52,10)}_{B_2^{n-1}} & \\
\bullet_{W_1} \ar[u]^{a_1} & \bullet_{W_2} \ar[u]^{a_2} & \cdots & \bullet_{W_{n-1}} \ar[u]^{a_{n-1}} & \bullet_{W_n} 
}
\end{equation} 
 For each $k = 1, \ldots, n-2$ there is a relation
\begin{equation}\label{eqn:handsaw-relation}
r_k = B_1^k B_2^k - B_2^{k+1} B_1^k + a_{k+1} b_{k+1}
\end{equation}
Each relation induces a map $\nu_{r_k} : \Rep(Q, {\bf v}) \rightarrow \Hom(V_k, V_{k+1})$. Therefore we see that for each $k=1, \ldots, n-2$, the vertex $V_{k+1}$ satisfies $V_{k+1} = h(r_k)$. For each of these vertices, we have $\{ a \in \mathcal{E} \, : \, h(a) = V_{k+1} \} = \{ a_{k+1}, B_1^k, B_2^{k+1} \}$. One can then see from \eqref{eqn:handsaw-relation} that the first condition of Definition \ref{def:complete-relations} is satisfied. Similarly, for each $k=1, \ldots, n-2$, the vertex $V_k$ satisfies $V_k = t(r_k)$ and $\{ a \in \mathcal{E} \, : \, t(a) = V_k \} = \{ b_{k+1}, B_1^k, B_2^k \}$. Again, one can see from \eqref{eqn:handsaw-relation} that the second condition of Definition \ref{def:complete-relations} is satisfied, and so the relations are complete.
\end{example}

The next example is an extended version of the ADHM quiver.

\begin{example}
Let $Q$ be a quiver with two vertices, labelled $V$ and $W$ in the diagram below, an arbitrary number of loops at the vertex $V$ (labelled $a_1, \ldots, a_n$) and edges $b_1$ and $b_2$ between $V$ and $W$.
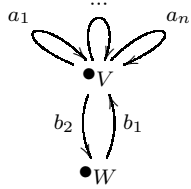
\begin{equation*}
\xymatrix{
\bullet_{V} \ar@`{(-15,5),(-8,10)}^{a_1} \ar@`{(-5,10),(5,10)}^{\cdots} \ar@`{(8,10),(15,5)} ^{a_n} \ar@/_0.5pc/[d]_{b_2} \\
\bullet_{W} \ar@/_0.5pc/[u]_{b_1} 
}
\end{equation*} 
Given any permutation $\sigma \in S_n$, define the relations
\begin{equation}\label{eqn:extended-ADHM-relations}
r = a_1 a_{\sigma(1)} + a_2 a_{\sigma(2)} + \cdots + a_n a_{\sigma(n)} + b_1 b_2 \quad \text{and} \quad r' = b_2 b_1 .
\end{equation}
Then it is easy to verify directly that these relations are complete, since each edge $a_1, \ldots, a_n, b_1, b_2$ is the leading edge in a path appearing nontrivially in one of the relations $r$ or $r'$, and also the tail of a unique path appearing in either $r$ or $r'$.
\end{example}

As in Lemma \ref{lem:linearise-slice} above, consider a representation $x = x_1 \oplus x_2 \in \Rep(Q, {\bf v_1}, \mathcal{R}) \oplus \Rep(Q, {\bf v_2}, \mathcal{R})$ and let $\delta x \in \Hom^1(Q, {\bf v_2}, {\bf v_1})$. For each path $p = a_2 a_1$ of length $2$, from \eqref{eqn:path-derivative-def} we have
\begin{equation*}
dp_x(\delta x) = (x_1)_{a_2} (\delta x)_{a_1} + (\delta x)_{a_2} (x_2)_{a_1} .
\end{equation*}
and therefore the adjoint of 
\begin{equation*}
\begin{tikzcd}
\Hom^1(Q, {\bf v_2}, {\bf v_1}) \arrow{r}{dp_x} & \Hom(V_{t(a_1)}, V_{h(a_2)})
\end{tikzcd}
\end{equation*}
is 
\begin{equation*}
(dp_x)^*(u) = (x_1)_{a_2}^* u + u (x_2)_{a_1}^* .
\end{equation*}
Therefore \eqref{eqn:relation-derivative-def} becomes
\begin{align*}
d \nu_{x}(\delta x) & = \sum_{r \in \mathcal{R}} \sum_{p \in \mathcal{P}_{t(r), h(r)}} \lambda_p(r) d p_x(\delta x) \\
 & = \sum_{r \in \mathcal{R}} \sum_{p \in \mathcal{P}_{t(r), h(r)}} \lambda_p(r) \left( (x_1)_{a_2(p)} (\delta x)_{a_1(p)} + (\delta x)_{a_2(p)} (x_2)_{a_1(p)} \right)
\end{align*}
and so for $u = (u_r)_{r \in \mathcal{R}}$ the adjoint 
\begin{equation*}
\begin{tikzcd}
\Rel(Q, {\bf v_2}, {\bf v_1}, \mathcal{R}) \arrow{r}{d \nu_x^*} & \Hom^1(Q, {\bf v_2}, {\bf v_1})
\end{tikzcd}
\end{equation*}
is given by
\begin{equation}\label{eqn:adjoint-derivative-relations}
d \nu_x^*(u) = \sum_{r \in \mathcal{R}} \sum_{p \in \mathcal{P}_{t(r), h(r)}} \overline{\lambda_p(r)} \left( (x_1)_{a_2}^* u_r + u_r (x_2)_{a_1}^* \right) .
\end{equation}

The next lemma gives a formula for the cokernel of $d\nu_x$, which only depends on the dimension of the image of $x_2$.
\begin{lemma}\label{lem:cokernel-dimension}
Let $\mathcal{R}$ be a set of complete quadratic relations, let $x = x_1 \oplus x_2 \in \Rep(Q, {\bf v_1}) \oplus \Rep(Q, {\bf v_2})$ and consider
\begin{equation*}
\begin{tikzcd}
\Hom^1(Q, {\bf v_2}, {\bf v_1}) \arrow{r}{d \nu_x} & \Rel(Q, {\bf v_2}, {\bf v_1}, \mathcal{R})  
\end{tikzcd}
\end{equation*}
from \eqref{eqn:relation-derivative-def}. Let ${\bf r}$ be the dimension vector of $(\im x_1)^\perp$. Then if $x_2 = 0$ we have 
\begin{equation}\label{eqn:cokernel-dimension}
\coker d \nu_x \cong \Rel(Q, {\bf v_2}, {\bf r}, \mathcal{R}) .
\end{equation}
In particular
\begin{equation}\label{eqn:explicit-cokernel-dimension}
\dim_\C \coker d \nu_x = \sum_{r \in \mathcal{R}} (v_2)_{t(r)} \dim_\C {\bf r}_{h(r)}
\end{equation}
\end{lemma}

\begin{proof}
First we show that $u \in \Rel(Q, {\bf v_2}, {\bf r}, \mathcal{R})$ implies that $d \nu_x^*(u) = 0$. From the definition of ${\bf r}$, $u \in \Rel(Q, {\bf v_2}, {\bf r}, \mathcal{R})$ implies that $\im u \perp \im x_1$, therefore $(x_1)_a^* u_r = 0$ for all edges $a \in \mathcal{E}$ and all relations $r \in \mathcal{R}$. Then \eqref{eqn:adjoint-derivative-relations} with $x_2 = 0$ shows that $d \nu_x^*(u) = 0$. 

Now we use the completeness of the relations to show the converse. If $u \in \Rel(Q, {\bf v_2}, {\bf v_1}, \mathcal{R}) \setminus \Rel(Q, {\bf v_2}, {\bf r}, \mathcal{R})$ then the image of $u$ is not perpendicular to the image of $x$, and so there exists an edge $a \in \mathcal{E}$ and a relation $r \in \mathcal{R}$ such that $(x_1)_a^* u_r \neq 0$. One consequence of this is that $h(a) = h(r)$, and so the assumption that the relations are complete and quadratic implies that there is a path $p = a a'$ with $\lambda_p(r) \neq 0$. Note that this implies that $t(a') = t(r)$, and therefore there exists a homomorphism $(\delta x)_{a'} \in \Hom(V_{t(a')}, V_{h(a')}) \subset \Hom^1(Q, {\bf v_2}, {\bf v_1})$ such that
\begin{equation*}
0 \neq \left< (x_1)_a (\delta x)_{a'}, u_r \right> = \left< (\delta x)_{a'}, (x_1)_a^* u_r \right> .
\end{equation*}
The completeness of the relations then implies that $r$ is the unique relation with $t(a') = t(r)$ and $p$ is the unique path with $\lambda_p(r) \neq 0$ and $a_1(p) = a'$. This uniqueness means that in the following expression for $d \nu_x^*$, there is only one nonzero term in the sum 
\begin{equation*}
\left< (\delta x)_{a'}, d \nu_x^*(u) \right> = \left< d \nu_x (\delta x)_{a'}, u \right> = \sum_r \sum_p \left< \lambda_p(r) (x_1)_a (\delta x)_{a'}, u_r \right> =  \left< (\delta x)_{a'}, \overline{\lambda_p(r)} (x_1)_a^* u_r \right> .
\end{equation*}
Therefore $u \in \Rel(Q, {\bf v_2}, {\bf v_1}, \mathcal{R}) \setminus \Rel(Q, {\bf v_2}, {\bf r}, \mathcal{R})$ implies that $d \nu_x^*(u) \neq 0$, or equivalently $d \nu_x^*(u) = 0$ implies that $u \in \Rel(Q, {\bf v_2}, {\bf r}, \mathcal{R})$. 
\end{proof}


Finally, there are a number of examples where $d \nu_x$ is surjective, which will be useful in Section \ref{sec:convolution-cup-product} to show that certain subvarieties are smooth. An important example is the case of the relations for a Nakajima quiver from Example \ref{ex:nakajima-complete}, where $\nu = \mu_\mathbb{C}$ is the complex moment map. In this case $\Rel(Q, {\bf v}, {\bf v}, \mathcal{R}) = \Hom^0(Q, {\bf v}, {\bf v})$ and so the deformation complex for the tangent space has the form
\begin{equation*}
\begin{tikzcd}
\Hom^0(Q, {\bf v}, {\bf v}) \arrow{r}{\rho_x^\mathbb{C}} & \Hom^1(Q, {\bf v}, {\bf v}) \arrow{r}{d \nu_x} & \Hom^0(Q, {\bf v}, {\bf v}) 
\end{tikzcd}
\end{equation*}
and therefore homomorphisms in the image of $d\nu_x$ correspond to elements of the Lie algebra $\mathfrak{g}_{\bf v} \cong \Hom^0(Q, {\bf v}, {\bf v})$. The adjoint of $d \nu_x$ is then
\begin{align*}
d \nu_x^*(u) & = \sum_{r \in \mathcal{R}} \sum_{p \in \mathcal{P}_{t(r), h(r)}} \lambda_p(r) \left( (x_1)_{a_2}^* u_r + u_r (x_2)_{a_1}^* \right) \\
 & = \left( u_{t(a)} x_a^* - x_a^* u_{h(a)} \right)_{a \in \mathcal{E}} .
\end{align*}
Therefore the transpose with respect to the Hermitian inner product on $\Hom^1(Q, {\bf v}, {\bf v})$ is
\begin{equation*}
\left( d \nu_x^*(u) \right)^* = \left( x_a u_{t(a)}^* - u_{h(a)}^* x_a \right)_{a \in \mathcal{E}} = \rho_x^\mathbb{C}(u^*) .
\end{equation*}
In particular, $\ker d \nu_x^* \cong \ker \rho_x^\mathbb{C}$. If the representation is stable then $\ker \rho_x^\mathbb{C}$ consists of the diagonal elements of $\mathfrak{g}_{\bf v}$ and therefore $d \nu_x$ is surjective onto $\mathfrak{g}_{\bf v}/\mathbb{C}$.

\subsection{Properties of the norm-square of the moment map}\label{subsec:moment-map-def}


Let $Q$ be a quiver with dimension vector ${\bf v} = ( v_k )_{k \in \mathcal{I}}$, and fix a Hermitian structure on the vector spaces $V_k \cong \C^{v_k}$. There is an associated symplectic structure on $\Rep(Q, {\bf v})$, defined as follows. Given $x \in \Rep(Q, {\bf v})$ and tangent vectors $\delta x_1, \delta x_2 \in T_x \Rep(Q, {\bf v}) \cong \Rep(Q, {\bf v})$, define the metric
\begin{equation}\label{eqn:metric-def}
g(\delta x_1, \delta x_2) := \sum_{a \in \mathcal{E}} \Re \tr \left( (\delta x_1)_a (\delta x_2)_a^* \right) ,
\end{equation}
and symplectic structure
\begin{equation}\label{eqn:symplectic-structure}
\omega(\delta x_1, \delta x_2) := \sum_{a \in \mathcal{E}} \Im \tr \left( (\delta x_1)_a (\delta x_2)_a^* \right) .
\end{equation}
Note that $\omega(\delta x_1, \delta x_2) = g(-i \delta x_1, \delta x_2)$, so that the complex structure $I = -i \cdot \id$ is compatible with the metric. With this complex structure and metric, the vector space $\Rep(Q, {\bf v})$ has the structure of a K\"ahler manifold.

On the Lie algebra $\mathfrak{g}_{\bf v}$, define the inner product
\begin{equation}\label{eqn:inner-product-def}
\left< u_1, u_2 \right> = \sum_{k \in \mathcal{I}} \tr \left( (u_1)_k (u_2)_k^* \right)
\end{equation}
and note that this is invariant under the adjoint action of $K_{\bf v} \subset G_{\bf v}$. Using the inner product we identify $\mathfrak{k} \cong \mathfrak{k}^*$. We will use $(\rho_x^\C)^* : T_x \Rep(Q, {\bf v}) \rightarrow \mathfrak{g}_{\bf v}$ to denote the adjoint of the infinitesimal action $\rho_x^\C : \mathfrak{g}_{\bf v} \rightarrow T_x \Rep(Q, {\bf v})$ with respect to the metric $g$ and the inner product on $\mathfrak{g}_{\bf v}$.

The action of $K_{\bf v}$ on $\Rep(Q, {\bf v})$ preserves the symplectic form $\omega$ and is Hamiltonian with moment map given by
\begin{equation}\label{eqn:moment-map-def}
\mu(x) = \frac{1}{2i} \sum_{a \in \mathcal{E}} [x_a, x_a^*] \in \mathfrak{k}_{\bf v}^* \cong \mathfrak{k}_{\bf v} .
\end{equation}
The moment map is $K_{\bf v}$-equivariant $\mu(k \cdot x) = \Ad_k \mu(x)$ for all $k \in K_{\bf v}$ and $x \in \Rep(Q, {\bf v})$. The centre of $\mathfrak{k}_{\bf v}$ is
\begin{equation*}
Z(\mathfrak{k}_{\bf v}) = \left\{ \left( \alpha_k \cdot \id_{V_k} \right)_{k \in \mathcal{I}} \mid \alpha_k \in \mathbb{C} \, \text{for all $k \in \mathcal{I}$} \right\} .
\end{equation*}
We say that $\alpha = (i \alpha_k \cdot \id_{V_k})_{k \in \mathcal{I}} \in Z(\mathfrak{k}_{\bf v})$ is \emph{admissible} if $\sum_{k \in \mathcal{I}} \alpha_k \dim_\C V_k = 0$. Note that $\tr(\mu(x)) = 0$ and so $\mu^{-1}(\alpha)$ is empty if $\alpha$ is not admissible. Given an admissible $\alpha \in Z(\mathfrak{k}_{\bf v})$, the \emph{symplectic quotient} is
\begin{equation}\label{eqn:symplectic-quotient-def}
\mathcal{M}_\alpha(Q, {\bf v}) := \mu^{-1}(\alpha) / K_{\bf v} .
\end{equation}
From now on $\alpha$ will always refer to an admissible central element of $\mathfrak{k}_{\bf v}$. Define the \emph{rank} and \emph{$\alpha$-degree} of a dimension vector ${\bf v}$ by 
\begin{align*}
\deg_\alpha(Q, {\bf v}) & := \sum_{k \in \mathcal{I}} \alpha_k v_k \\
\rank(Q, {\bf v}) & := \sum_{k \in \mathcal{I}} v_k
\end{align*}
The \emph{$\alpha$-slope} is $\slope_\alpha(Q, {\bf v}) := \deg_\alpha(Q, {\bf v}) / \rank(Q, {\bf v})$. A representation $x \in \Rep(Q, {\bf v})$ is \emph{$\alpha$-stable} (resp. $\alpha$-semistable) if and only if every proper non-zero subrepresentation satisfies
\begin{equation*}
\slope_\alpha(Q, {\bf v'}) < 0 \quad \text{(resp. $\slope_\alpha(Q, {\bf v'}) \leq 0$)} .
\end{equation*}
King \cite[Prop. 3.1]{King94} shows that slope stability in the above sense coincides with stability from affine GIT and that $\alpha$-polystable representations are isomorphic to minimisers of $\| \mu - \alpha \|^2$.

The following definition is used in Lemma \ref{lem:crit-set-cohomology} to define critical points as direct sums of minimisers of the norm-square of a shifted moment map.
\begin{definition}\label{def:induced-central}
Let $Q$ be a quiver, ${\bf v}$ a dimension vector, and $\alpha = (\alpha_k)_{k \in \mathcal{I}}$ an admissible central element of $\mathfrak{k}_{\bf v}$. Given any dimension vector ${\bf 0} \leq {\bf v'} \leq {\bf v}$, the \emph{induced admissible central element} on $(Q, {\bf v'})$ is 
\begin{equation}\label{eqn:induced-central}
\alpha' = \left( (\alpha_k - \slope_\alpha(Q, {\bf v'})) \cdot \id_{V_k'} \right)_{k \in \mathcal{I}} .
\end{equation}

\end{definition}

Define the energy function $f : \Rep(Q, {\bf v}) \rightarrow \mathbb{R}$ given by
\begin{equation*}
f(x) = \| \mu(x) - \alpha \|^2 .
\end{equation*}

With respect to the metric \eqref{eqn:metric-def}, the gradient of $f$ is
\begin{equation*}
\grad f(x) = I \rho_x(\mu(x) - \alpha) .
\end{equation*}

In analogy with critical points of the Yang-Mills functional studied in \cite[Sec. 5]{AtiyahBott83}, the critical point equation $\grad f(x) = 0$ implies that the representation $x$ splits into subrepresentations $x = \bigoplus_{\ell=1}^n x_\ell$ such that each $x_\ell \in \Rep(Q, {\bf v_\ell})$ satisfies $\mu(x_\ell) = \alpha_\ell$, where the admissible central element $\alpha_\ell \in Z(\mathfrak{k}_{\bf v_\ell})$ is induced from $\alpha \in Z(\mathfrak{k}_{\bf v})$ by the construction of \eqref{eqn:induced-central}. Equivalently, a critical point is a direct sum of minimisers of $\| \mu(x_\ell) - \alpha_\ell \|^2$ on $\Rep(Q, {\bf v_\ell})$. Therefore each critical set is labelled by a vector of dimension vectors
\begin{equation}\label{eqn:crit-set-vector}
({\bf v}_1, \ldots, {\bf v}_k) \quad \text{such that ${\bf v}_1 + \cdots + {\bf v}_k$} .
\end{equation}
This is explained in more detail in \cite[Sec. 2.4]{Wilkin17}. In the following we abuse the notation and use $\mu^{-1}(\alpha_\ell)$ to denote the set of minimisers of $\| \mu - \alpha_\ell \|^2$ on $\Rep(Q, {\bf v_\ell})$.

Let $\Crit(f) \subset \Rep(Q, {\bf v})$ be the set of all critical points of $f$. Given any critical point $x$, let $\beta = \mu(x) - \alpha$ and define
\begin{equation}\label{eqn:crit-set-def}
C_\beta := \Crit(f) \cap \mu^{-1}(\beta), \quad C_{K_{\bf v} \cdot \beta} := K_{\bf v} \cdot C_\beta .
\end{equation}
For representations of quivers, an argument directly analogous to that of Atiyah and Bott \cite{AtiyahBott83} classifies the critical sets in terms of Harder-Narasimhan type. In analogy with Atiyah and Bott's calculations for the Yang-Mills functional \cite{AtiyahBott83}, we can inductively compute the equivariant cohomology of the critical sets in terms of energy minimisers on the space of representations with smaller dimension vector.

\begin{lemma}\label{lem:crit-set-cohomology}
The $K_{\bf v}$-equivariant cohomology of $C_{K \cdot \beta}$ is
\begin{equation*}
H_{K_{\bf v}}^*(C_{K \cdot \beta}) \cong \bigotimes_{\ell=1}^n H_{K_{\bf v_\ell}}^*(\mu^{-1}(\alpha_\ell)) .
\end{equation*}
\end{lemma}


Given an initial condition $x_0$, define $\phi(x_0,t)$ to be the solution to the downwards gradient flow equation for the energy function $\| \mu - \alpha \|^2$
\begin{equation}\label{eqn:neg-grad-flow-def}
\frac{d}{dt} \phi(x_0, t) = -I \rho_{\phi(x_0,t)} (\mu(\phi(x_0,t)) - \alpha), \quad \phi(x_0,0) = x_0 .
\end{equation}
For each initial condition $x_0$, there exists a unique minimal $T \in [-\infty, 0)$ such that the solution $\phi(x_0, t)$ exists for all $t \in (T, \infty)$ and converges to a unique critical point as $t \rightarrow \infty$. If $T$ is finite then $\lim_{t \rightarrow T^+} f(\phi(x_0, t)) = \infty$. The results of \cite{Wilkin17} classify the isomorphism classes of solutions which converge to a critical point as $t \rightarrow - \infty$.

For each $t$ such that a solution of \eqref{eqn:neg-grad-flow-def} exists, there is a solution $g_t \in G_{\bf v}$ of the equation
\begin{equation}\label{eqn:group-neg-flow}
\frac{d g_t}{dt} g_t^{-1} = -i (\mu(g_t \cdot x_0) - \alpha)
\end{equation}
which satisfies $\phi(x_0, t) = g_t \cdot x_0$. In particular, the gradient flow preserves any subset $Z \subset \Rep(Q, {\bf v})$ preserved by the action of $G_{\bf v}$. If the subset $Z$ is closed then the limit of the flow is also contained in $Z$. Therefore we can define the gradient flow and its limit on any closed subset preserved by the action of $G_{\bf v}$, even if this subset is singular and the usual definition of the gradient vector field in terms of derivatives does not make sense.


\begin{definition}\label{def:unstable-def}
Let $Z \subset \Rep(Q, {\bf v})$ be a closed subset preserved by $G_{\bf v}$ and define $f : Z \rightarrow \mathbb{R}$ by $f(x) = \| \mu(x) - \alpha \|^2$. A \emph{critical point} of $f$ is a stationary point for the gradient flow \eqref{eqn:neg-grad-flow-def}. The \emph{unstable set} of $x$ is
\begin{equation}\label{eqn:unstable-set-def}
W_x^- := \left\{ y \in Z \mid \lim_{t \rightarrow -\infty} \phi(y, t) = x \right\} .
\end{equation}
Given a critical set $C_{K \cdot \beta}$, the \emph{unstable bundle} is
\begin{equation}\label{eqn:unstable-bundle-def}
W_{K \cdot \beta}^- := \left\{ y \in Z \mid \lim_{t \rightarrow - \infty} \phi(y, t) \in C_{K \cdot \beta} \right\} 
\end{equation}
and the \emph{stable set} or \emph{Morse stratum} is denoted
\begin{equation}\label{eqn:stable-set-def}
W_{K \cdot \beta}^+ := \left\{ y \in Z \mid \lim_{t \rightarrow \infty} \phi(y, t) \in C_{K \cdot \beta} \right\}  .
\end{equation}
\end{definition}

\subsubsection{The negative slice associated to the critical set}


Let $x$ be a critical point of $\| \mu - \alpha \|^2$, and let $\beta = \mu(x) - \alpha$. Since $I \rho_x(\beta) = 0$ then $e^{i \beta t} \cdot x = x$ for all $t \in \mathbb{R}$ and so there is an induced action of the one-parameter subgroup $\{ e^{i \beta t} \mid t \in \mathbb{R} \} \subset G_{\bf v}$ on the tangent space $T_x \Rep(Q, {\bf v})$. 

\begin{definition}\label{def:slice}
Let $Z \subset \Rep(Q, {\bf v})$ be a closed subset preserved by $G_{\bf v}$ and define $f : Z \rightarrow \mathbb{R}$ by $f(x) = \| \mu(x) - \alpha \|^2$. Let $x$ be a critical point of $f$ with $\beta = \mu(x) - \alpha$. The \emph{local slice} at $x$ is
\begin{equation*}
S_x := \left\{ \delta x \in T_x \Rep(Q, {\bf v}) \mid \delta x \in \ker (\rho_x^\C)^* \, \text{and} \, x + \delta x \in Z \right\} .
\end{equation*}
The \emph{negative slice} at $x$ is
\begin{equation}\label{eqn:neg-slice-def}
S_x^- := \left\{ \delta x \in S_x \mid \lim_{t \rightarrow \infty} e^{i \beta t} \cdot \delta x = 0 \right\} .
\end{equation}
Given a critical set $C_{K \cdot \beta}$ as in \eqref{eqn:crit-set-def}, consider the trivial bundle $C_{K \cdot \beta} \times \Rep(Q, {\bf v}) \rightarrow C_{K \cdot \beta}$. By identifying the fibre over $x$ with $T_x \Rep(Q, {\bf v}) \cong \Rep(Q, {\bf v})$, define the \emph{negative slice bundle}
\begin{equation}\label{eqn:neg-slice-bundle}
S_{K \cdot \beta}^- := \left\{ (x, \delta x) \in C_{K \cdot \beta} \times \Rep(Q, {\bf v}) \mid \delta x \in S_x^- \right\}
\end{equation}
together with the projection $p : S_{K \cdot \beta}^- \rightarrow C_{K \cdot \beta}$.
\end{definition}

Since the infinitesimal action $\rho_x^\C$ is $K_{\bf v}$-equivariant and the metric \eqref{eqn:metric-def} and inner product \eqref{eqn:inner-product-def} are both $K_{\bf v}$-invariant, then the adjoint $(\rho_x^\C)^*$ is also $K_{\bf v}$-equivariant. More explicitly, the action is given by
\begin{equation*}
\Ad_k ((\rho_x^\C)^*(\delta x)) = (\rho_{k \cdot x}^\C)^* (k \cdot \delta x) .
\end{equation*}
In particular, each $k \in K_{\bf v}$ defines an isomorphism $S_x \cong S_{k \cdot x}$ given by $\delta x \mapsto k \cdot \delta x$. Since the moment map is $K_{\bf v}$-equivariant then $\lim_{t \rightarrow \infty} e^{i \mu(x) t} \cdot \delta x = 0$ if and only if $\lim_{t \rightarrow \infty} e^{i \mu(k \cdot x)} \cdot (k \cdot \delta x) = 0$, and therefore this isomorphism restricts from the slice $S_x$ to the negative slice $S_x^- \stackrel{\cong}{\longrightarrow} S_{k \cdot x}^-$, which defines an action of $K_{\bf v}$ on $S_{K \cdot \beta}^-$. 


We can describe the points in the negative slice more explicitly as follows. Let $x$ be a critical point with $\mu(x) = \beta$, which implies that the representation is a direct sum of subrepresentations $x = \bigoplus_{\ell} x_\ell$ with dimension vector ${\bf v_\ell}$, where each subspace $\Vect(Q, {\bf v_\ell}) \hookrightarrow \Vect(Q, {\bf v})$ is an eigenspace for the action of $\beta$ and the eigenvalue is determined by $\slope_\alpha(Q, {\bf v_\ell})$. In the following we use the convention that the subspaces are ordered by increasing $\alpha$-slope, i.e. $\ell_1 < \ell_2$ if and only if $\slope_\alpha(Q, {\bf v_{\ell_1}}) < \slope_\alpha(Q, {\bf v_{\ell_2}})$. The condition $\lim_{t \rightarrow \infty} e^{i \beta t} \cdot \delta x = 0$ then implies that $\delta x$ is contained in the subspace
\begin{equation*}
\delta x \in \bigoplus_{\ell_1 < \ell_2} \Hom^1(Q, {\bf v_{\ell_2}}, {\bf v_{\ell_1}}) \subset \Hom^1(Q, {\bf v}, {\bf v}) = \Rep(Q, {\bf v}) .
\end{equation*}


The point of introducing the negative slice is that one can explicitly describe the pair $(S_{K \cdot \beta}^-, S_{K \cdot \beta}^- \setminus C_{K \cdot \beta})$, however for the purposes of Morse theory the natural object of study is the pair $(W_{K \cdot \beta}^-, W_{K \cdot \beta}^- \setminus C_{K \cdot \beta})$ (see \cite{Wilkin19} and Theorem \ref{thm:main-thm-morse} in this paper). The following theorem is one of the main results of \cite{Wilkin17} and shows that the topology of the negative slice is the same as that of the unstable set, and therefore we can reduce the study of flow lines inside the unstable set to the study of certain explicit subspaces of the negative slice (cf. Section \ref{subsec:flow-lines}).  

\begin{theorem}\label{thm:slice-homeo}
Let $Z \subset \Rep(Q, {\bf v})$ be a closed subset preserved by $G_{\bf v}$ and define $f : Z \rightarrow \mathbb{R}$ by $f(x) = \| \mu(x) - \alpha \|^2$. Then for each critical set $C_{K \cdot \beta}$ there exists a neighbourhood $U$ of $C_{K \cdot \beta}$ in $W_{K \cdot \beta}^-$, a neighbourhood $V$ of $C_{K \cdot \beta} \times \{0\}$ in $S_{K \cdot \beta}^-$ and a $K_{\bf v}$-equivariant homeomorphism of pairs 
\begin{equation}\label{eqn:slice-homeo}
H_\beta : (U, U \setminus C_{K \cdot \beta}) \stackrel{\cong}{\longrightarrow} (V, V \setminus (C_{K \cdot \beta} \times \{0\})) .
\end{equation}
\end{theorem}

\subsection{Critical sets and the negative slice for framed quivers}\label{subsec:canonical-central-element}

In this section we consider a quiver $Q$ with vertices $\mathcal{I}$ and edges $\mathcal{E}$, and let ${\bf v} = ( v_k )_{k \in \mathcal{I}}$ a dimension vector such that one vertex (which we label $\infty$) has dimension $1$.  Define $\mathcal{I}' = \mathcal{I} \setminus \{ \infty \}$ to be the set of remaining vertices of $Q$. 

\begin{definition}\label{def:aasp}
For such a quiver $Q$ and dimension vector ${\bf v} = (v_k)_{k \in \mathcal{I}}$, the \emph{canonical central element} $\alpha(Q, {\bf v}) := ( \alpha_k)_{k \in \mathcal{I}}$ is given by
\begin{equation}\label{eqn:Nakajima-stability}
\alpha_k := \left\{ \begin{matrix} 1 & k \in \mathcal{I}' \\ -\sum_{j \in \mathcal{I}'} v_j  & k = \infty \end{matrix} \right. 
\end{equation}
\end{definition}

Given a dimension vector $0 \leq {\bf v'} = ( v_k' )_{k \in \mathcal{I}} \leq {\bf v}$ with $v_\infty' = 1$, the induced central element $\alpha'$ from \eqref{eqn:induced-central} is $\alpha' = (i \alpha_k' \cdot \id_{V_k'})_{k \in \mathcal{I}}$, where 
\begin{equation}\label{eqn:induced-canonical-scalar}
\alpha_k' = \left( \frac{1 + \sum_{k \in \mathcal{I}'} v_k}{1 + \sum_{k \in \mathcal{I}'} v_k'} \right) \alpha_k \quad \text{for each $k \in \mathcal{I}$} .
\end{equation}
If $v_\infty' = 0$ then the induced central element $\alpha'$ is zero.

\begin{remark}\label{rem:canonical-central-element}
Note that any subrepresentation containing the vertex $\infty$ must have negative slope and any subrepresentation that does not contain the vertex $\infty$ must have positive slope. Therefore, when constructing the Harder-Narasimhan filtration (cf. \cite{Reineke03}) with respect to this stability parameter, we see that 
\begin{enumerate}

\item  the maximal $\alpha$-semistable subrepresentation is the largest subrepresentation that does not contain the vertex $\infty$,

\item the quotient by the maximal semistable subrepresentation is an $\alpha$-stable representation, for which every subrepresentation must contain the vertex $\infty$, and

\item the Harder-Narasimhan filtration of an $\alpha$-unstable representation is a two step filtration. We will denote the type of this filtration by the dimension vector of the quotient by the maximal semistable subrepresentation.

\item A construction of Crawley-Boevey \cite{Crawley01} shows that the moduli space of stable representations of any framed quiver is homeomorphic to the unframed moduli space $\mathcal{M}_\alpha(Q, {\bf v})$ as defined in \eqref{eqn:symplectic-quotient-def} for a quiver which has an extra vertex of dimension $1$, and therefore fits into the above construction. Via this correspondence, the central element of Definition \ref{def:aasp} determines the same stability condition as for the hyperk\"ahler quiver varieties \cite{Nakajima94} and the handsaw quiver varieties \cite{Nakajima12} studied by Nakajima.

\end{enumerate}
\end{remark}

With respect to the canonical central element $\alpha$ from Definition \ref{def:aasp}, \cite[Prop. 3.13]{Wilkin17} shows that any critical point $x$ splits as a direct sum of two representations $x = x_1 \oplus x_2$ with corresponding dimension vectors ${\bf v_1}$, ${\bf v_2}$ with moment maps $\mu(x_1)$ given by \eqref{eqn:induced-canonical-scalar} and $\mu(x_2) = 0$. Therefore the vector labelling the critical set (described in \eqref{eqn:crit-set-vector}) is $({\bf v_1}, {\bf v_2}) = ({\bf v_1}, {\bf v} - {\bf v_1})$ and so, with respect to this stability parameter, the critical sets are labelled by the vector ${\bf v_1}$.


Given such a critical point $x$ on a closed $G_{\bf v}$-invariant subset $Z \subset \Rep(Q, {\bf v})$, the \emph{negative slice} from \eqref{eqn:neg-slice-def} becomes
\begin{equation}\label{eqn:neg-slice-framed}
S_x^- = \left\{ \delta x \in \Hom^1(Q, {\bf v_2}, {\bf v_1}) \cap \ker (\rho_x^\C)^* \mid x + \delta x \in Z \right\} .
\end{equation}

Now let $\mathcal{R}$ be a finite set of relations on the quiver and consider the case where $Z$ is the subset $\Rep(Q, {\bf v}, \mathcal{R}) = \nu^{-1}(0) \subset \Rep(Q, {\bf v})$ as defined in \cite[Sec. 3]{Wilkin17}. Since $x = x_1 \oplus x_2 \in \Rep(Q, {\bf v_1}, \mathcal{R}) \oplus \Rep(Q, {\bf v_2}, \mathcal{R})$ and $S_x^- \subset \Hom^1(Q, {\bf v_2}, {\bf v_1})$, then \cite[Lem. 3.29]{Wilkin17} states that
\begin{equation}\label{eqn:slice-as-kernel}
S_x^- = \Hom^1(Q, {\bf v_2}, {\bf v_1}) \cap \ker (\rho_x^\C)^* \cap \ker d \nu_x .
\end{equation}

At each critical point we have the following deformation complex
\begin{equation}\label{eqn:neg-def-complex}
\Hom^0(Q, {\bf v_2}, {\bf v_1}) \stackrel{\rho_x^\C}{\xrightarrow{\hspace{0.5cm}}} \Hom^1(Q, {\bf v_2}, {\bf v_1}) \stackrel{d\nu_x}{\xrightarrow{\hspace{0.5cm}}} \Rel(Q, {\bf v_2}, {\bf v_1}, \mathcal{R}) .
\end{equation}
Define the cohomology groups $\mathcal{H}^0(Q, {\bf v_2}, {\bf v_1}) := \ker \rho_x^\C$, $\mathcal{H}^1(Q, {\bf v_2}, {\bf v_1}, \mathcal{R}) := \ker (\rho_x^\C)^* \cap \ker d \nu_x$ and $\mathcal{H}^2(Q, {\bf v_2}, {\bf v_1}, \mathcal{R}) := \ker (d \nu_x)^*$ at each term of the complex and define 
\begin{equation*}
h^p(Q, {\bf v_2}, {\bf v_1}, \mathcal{R}) := \dim_\C \mathcal{H}^p(Q, {\bf v_2}, {\bf v_1}, \mathcal{R})
\end{equation*}
for each $p=0,1,2$. From \eqref{eqn:slice-as-kernel} we have
\begin{equation}\label{eqn:slice-as-cohomology}
S_x^- = \mathcal{H}^1(Q, {\bf v_2}, {\bf v_1}, \mathcal{R}) .
\end{equation}


\begin{definition}
The index of the complex \eqref{eqn:neg-def-complex} can be written in terms of the Ringel form of \eqref{eqn:ringel-def}, and is denoted as follows
\begin{align}\label{eqn:relation-index-def}
\begin{split}
\left< {\bf v_2}, {\bf v_1} \right>_{\mathcal{R}} & := \dim_\C \Hom^0(Q, {\bf v_2}, {\bf v_1}) - \dim_\C \Hom^1(Q, {\bf v_2}, {\bf v_1}) + \dim_\C \Rel(Q, {\bf v_2}, {\bf v_1}, \mathcal{R}) \\
 & = \left< {\bf v_2}, {\bf v_1} \right> + \dim_\C \Rel(Q, {\bf v_2}, {\bf v_1}, \mathcal{R}) .
\end{split}
\end{align}
\end{definition}


The following lemma allows us to compute the dimension of the negative slice in terms of the index of the complex and the dimension of the second cohomology group of the deformation complex.

\begin{lemma}\label{lem:neg-slice-dimension}
\begin{equation}\label{eqn:neg-dimension}
\dim_\C S_x^- = h^2(Q, {\bf v_2}, {\bf v_1}, \mathcal{R}) - \left< {\bf v_2}, {\bf v_1} \right>_{\mathcal{R}} .
\end{equation}
\end{lemma}

\begin{proof}
A homomorphism in $\ker \rho_x^\C \subset \Hom^0(Q, {\bf v_2}, {\bf v_1})$ is necessarily zero since it maps a semistable representation into a stable representation of smaller slope (cf. \cite[Lem. 2.3]{Reineke03}), and therefore $h^0 = 0$.  Since $S_x^- =  \mathcal{H}^1(Q, {\bf v_2}, {\bf v_1}, \mathcal{R})$, then the result follows from the index formula
\begin{equation*}
\left< {\bf v_2}, {\bf v_1} \right>_{\mathcal{R}} = h^0 - h^1 + h^2 \quad \Leftrightarrow \quad h^1 = h^2 + h^0 - \left< {\bf v_2}, {\bf v_1} \right>_{\mathcal{R}} . \qedhere
\end{equation*}
\end{proof}


From now on we drop the notation for the dimension vector ${\bf v}$ from $K_{\bf v}$ since the meaning will be clear from the context. Given a critical set $C_{K \cdot \beta}$, each critical point $x$ is a direct sum $x_1 \oplus x_2$, which determines a direct sum $\Vect(Q, {\bf v}) \cong \Vect(Q, {\bf v_1}) \oplus \Vect(Q, {\bf v_2})$ of the vector spaces at each vertex according to the eigenvalues of $\beta = \mu(x) - \alpha$. Define the following bundles over $C_{K \cdot \beta}$.
\begin{align*}
\underline{\Hom}^0(Q, {\bf v_2}, {\bf v_1}) & = \left\{ (x, u) \in C_{K \cdot \beta} \times \Hom^0(Q, {\bf v}, {\bf v}) \mid u \in \Hom^0(Q, {\bf v_2}, {\bf v_1}) \right\} \\
\underline{\Hom}^1(Q, {\bf v_2}, {\bf v_1}) & = \left\{ (x, \delta x) \in C_{K \cdot \beta} \times \Rep(Q, {\bf v}) \mid \delta x \in \Hom^1(Q, {\bf v_2}, {\bf v_1}) \right\} .
\end{align*}
Since the moment map is $K$-equivariant, then so is this decomposition, and so there is an induced action of $K$ on these bundles. The complex \eqref{eqn:neg-def-complex} extends to a complex of bundle homomorphisms
and \eqref{eqn:slice-as-cohomology} shows that $S_{K \cdot \beta}^-$ is the middle cohomology of this complex.




We conclude this section with a result on the relative equivariant cohomology of the negative slice which is a singular space analog of the results of Atiyah \& Bott \cite[Sec. 13]{AtiyahBott83} and Kirwan \cite[Sec. 4.23]{Kirwan84} for the critical sets and negative eigenbundle of the Hessian.

Let $K_\beta$ be the isotropy group of $\beta \in \mathfrak{k}^*$ with respect to the coadjoint action. The critical set $C_{K \cdot \beta}$ has the structure of a fibre product $C_{K \cdot \beta} \cong C_\beta \times_{K_\beta} K$ (cf. \cite[Sec. 4.22]{Kirwan84}). In particular, we have $H_K^*(C_{K \cdot \beta}) \cong H_{K_\beta}^*(C_\beta)$. 

Let $S_\beta^-$ be the pullback of $S_{K\cdot \beta}^-$ by the inclusion $C_\beta \hookrightarrow C_{K \cdot \beta}$. Then $S_{K \cdot \beta}^-$ and $S_{K \cdot \beta}^- \setminus C_{K \cdot \beta}$ also have a fibre product  structure $S_{K \cdot \beta}^- \cong S_\beta^- \times_{K_\beta} K$ and $S_{K \cdot \beta}^- \setminus C_{K \cdot \beta} \cong (S_\beta^- \setminus C_\beta) \times_{K_\beta} K$, and so we have the following commutative diagram
\begin{equation}\label{eqn:reduce-structure-group}
\xymatrix{
H_K^*(S_{K \cdot \beta}^-, S_{K \cdot \beta}^- \setminus C_{K \cdot \beta}) \ar[r] \ar[d]^\cong & H_K^*(S_{K \cdot \beta}^-) \cong H_K^*(C_{K \cdot \beta}) \ar[d]^\cong \\
H_{K_\beta}^*(S_{\beta}^-, S_{\beta}^- \setminus C_{\beta}) \ar[r] & H_{K_\beta}^*(S_{\beta}^-) \cong H_{K_\beta}^*(C_{\beta})
}
\end{equation} 

The motivation for this is explained in the next section, where the main theorem of Morse theory shows that the question of studying the terms in the spectral sequence for the Morse stratification by the norm square of the moment map reduces to studying the relative cohomology groups for the pair $(W_{K \cdot \beta}^-, W_{K \cdot \beta}^- \setminus C_{K \cdot \beta})$. The above results show that this is equivalent to studying the relative cohomology groups of the pair $(S_\beta^-, S_\beta^- \setminus C_\beta)$ in $K_\beta$-equivariant cohomology, which is simpler because the negative slice is a linearised version of the unstable set which can be studied explicitly.


\subsection{Reduction to the first component of the critical set}

In this section we restrict to the case where the relations in the quiver are determined by paths of the same length, and therefore the relation map $\nu$ is a homogeneous polynomial. In this case we can prove that the relative equivariant cohomology groups of the previous section simplify, which will be more convenient for constructing the cup product and differentials on the Morse complex. This restriction is not too severe, as the class of quivers with homogeneous relations includes important examples such as Nakajima quivers and the handsaw quivers studied in \cite{Nakajima12}.

\begin{definition}
A relation is \emph{homogeneous} if it is defined by paths of the same length. A \emph{homogeneous set of relations} is a set $\mathcal{R}$ in which every relation is homogeneous.
\end{definition}

\begin{remark}
The above definition of a homogeneous set of relations only requires that the paths defining each relation have the same length. Two different relations in a homogeneous set may have different lengths.

If a relation $r$ is homogeneous then the algebraic map $\nu_r : \Rep(Q, {\bf v}) \rightarrow \Hom(V_{t(r)}, V_{h(r)})$ (cf. \cite[(3.1)]{Wilkin17}) is defined by a homogeneous polynomial in the components $\bigoplus_{a \in \mathcal{E}} \Hom(V_{t(a)}, V_{h(a)})$ of the representation.
\end{remark}

Recall from \cite[Prop. 3.13]{Wilkin17} that $C_\beta = C_{\beta_1} \times C_{\beta_2}$, where $C_{\beta_1} / K_{\beta_1} = \mathcal{M}(Q, {\bf v_1})$ and $C_{\beta_2} = \mu^{-1}(0)$ on $\Rep(Q, {\bf v_2})$. Define $S_{\beta_1}^-$ to be the restriction of the negative slice to the subset $C_{\beta_1} \times \{0\} \subset C_\beta$. In this section we show that when the quiver has homogeneous relations, the relative cohomology groups $H_{K_\beta}^*(S_\beta^-, S_\beta^- \setminus C_\beta)$ and $H_{K_\beta}^*(S_{\beta_1}^-, S_{\beta_1}^- \setminus C_{\beta_1})$ are isomorphic (cf. Corollary \ref{cor:reduce-to-zero} below). Therefore, in defining the differentials and cup product in Section \ref{sec:main-theorem-localisation}, it is sufficient to restrict attention to the subset $C_{\beta_1} \times \{0\} \subset C_{\beta_1} \times C_{\beta_2} = C_\beta$ on which the computations are much simpler.

\begin{lemma}
If the relations are homogeneous then the map $x \mapsto tx$ for $t \in [0,1]$ defines a $K$-equivariant deformation retract of $\nu^{-1}(0)$ and $\mu^{-1}(0) \cap \nu^{-1}(0)$ to the zero representation.
\end{lemma}

\begin{proof}
Since $\nu = \bigoplus_{r \in \mathcal{R}} \nu_r$ is a direct sum of homogeneous polynomials then $\nu^{-1}(0)$ is preserved by the map $x \mapsto tx$. Similarly, the components of $\mu$ are homogeneous polynomials of order $2$ (cf. \eqref{eqn:moment-map-def}), and so $\mu^{-1}(0)$ is also preserved by $x \mapsto tx$. 
\end{proof}

\begin{corollary}\label{cor:critical-level-isomorphism}
\begin{equation}\label{eqn:critical-level-isomorphism}
H_{K_\beta}^*(S_\beta^-) \cong H_{K_\beta}^*(C_\beta) \cong H_{K_\beta}^*(C_{\beta_1}) \cong H_{K_\beta}^*(S_{\beta_1}^-).
\end{equation}
Moreover, since $K_{\beta_1}$ acts on $C_{\beta_1}$ with isotropy group the diagonal $\U(1) \subset K_{\beta_1}$ then we have 
\begin{equation*}
H_{K_\beta}^*(C_{\beta_1}^-) \cong H^*(\mathcal{M}(Q, {\bf v_1})) \otimes H^*(\mathsf{BU}(1)) \otimes H^*(BK_{\beta_2}) .
\end{equation*}
\end{corollary}

For small values of $\varepsilon > 0$, the action of the one-parameter subgroup $\{ e^{i t \beta} \}$ defines a $K_\beta$-equivariant homotopy equivalence $S_{\beta_1}^- \setminus C_{\beta_1} \cong S_\beta^- \cap f^{-1}(c - \varepsilon)$. 


For any $t \in \C^*$, there is a $K_{\beta}$-equivariant isomorphism of the middle cohomology of the complex \eqref{eqn:neg-def-complex} over the critical points $x = (x_1, x_2)$ and $x_t = (x_1, t x_2)$. Therefore, for any neighbourhood $U$ of $S_{\beta_1}^- \cap f^{-1}(c - \varepsilon)$ in $S_\beta^- \cap f^{-1}(c - \varepsilon)$, we can define a $K_\beta$-equivariant deformation retract of $S_\beta^- \cap f^{-1}(c - \varepsilon)$ onto $U$.


Since $S_\beta^-$, $S_{\beta_1}^-$ and $f^{-1}(c - \varepsilon)$ are all real analytic, then \cite[Thm. 1.1]{PflaumWilkin19} shows that there is a neighbourhood $U$ of $S_{\beta_1}^- \cap f^{-1}(c - \varepsilon)$ in $S_\beta^- \cap f^{-1}(c - \varepsilon)$ and a $K_\beta$-equivariant deformation retract of $U$ onto $S_{\beta_1}^- \cap f^{-1}(c - \varepsilon)$.

Combining the above two deformation retracts gives us an isomorphism
\begin{equation}\label{eqn:lower-level-set-isomorphism}
H_{K_\beta}^*(S_\beta^- \setminus C_\beta) \cong H_{K_\beta}^*(S_\beta^- \cap f^{-1}(c - \varepsilon)) \cong H_{K_\beta}^*(S_{\beta_1}^- \cap f^{-1}(c - \varepsilon)) \cong H_{K_\beta}^*(S_{\beta_1}^- \setminus C_{\beta_1})
\end{equation}

The isomorphisms \eqref{eqn:critical-level-isomorphism} and \eqref{eqn:lower-level-set-isomorphism} together with the five lemma gives us the following result.
\begin{corollary}\label{cor:reduce-to-zero}
For any quiver with homogeneous relations, we have
\begin{equation}\label{eqn:reduce-relative-cohomology}
H_{K_\beta}^*(S_\beta^-, S_\beta^- \setminus C_\beta) \cong H_{K_{\beta}}^*(S_{\beta_1}^-, S_{\beta_1}^- \setminus C_{\beta_1})  .
\end{equation}
\end{corollary}

\subsection{Spaces of flow lines between critical points}\label{subsec:flow-lines}

Throughout this section we will consider the space of flow lines between two given critical sets. The flow is the \emph{negative} gradient flow $\phi(x,t)$ defined in \eqref{eqn:neg-grad-flow-def}. The lower critical set will always be denoted $C_\ell$ (with Harder-Narasimhan type ${\bf v}_\ell$; cf. Remark \ref{rem:canonical-central-element}) and the upper critical set by $C_u$ (with type ${\bf v}_u$).


\begin{definition}\label{def:space-of-flow-lines}
Given two critical sets $C_\ell$ and $C_u$ with $f(C_\ell) < f(C_u)$, the space of representations that flow up to $C_u$ and down to $C_\ell$ is
\begin{equation*}
\mathcal{F}_{\ell,0}^{u,0} := \left\{ x \in Z \, \mid \, \lim_{t \rightarrow \infty} \phi(x, t) \in C_\ell, \quad \lim_{t \rightarrow -\infty} \phi(x,t) \in C_u \right\} .
\end{equation*}
The flow defines an $\mathbb{R}$ action on $\tilde{\mathcal{F}}_\ell^u$ and the \emph{space of flow lines} is
\begin{equation*}
\tilde{\mathcal{F}}_{\ell,0}^{u,0} := \mathcal{F}_{\ell,0}^{u,0} / \mathbb{R} .
\end{equation*}
\end{definition}


The main technical result of \cite{Wilkin17} is that any representation $x$ such that $\lim_{t \rightarrow -\infty} \phi(x,t) \in C_u$ is isomorphic to a representation in the negative slice of $C_u$ \cite[Thm. 4.22]{Wilkin17}. On restricting to a neighbourhood of the critical set, this determines a homeomorphism of pairs $H : (W_u, W_{u,0}) \stackrel{\cong}{\rightarrow} (S_u^-, S_{u,0}^-)$ (cf. \cite[Cor. 4.24]{Wilkin17}) and via this homeomorphism we can consider the space of flow lines as a subset of the negative slice bundle $H(\mathcal{F}_{\ell,0}^{u,0}) \subset S_{u,0}$ in a neighbourhood of the critical set. 

The condition $x \in \mathcal{F}_{\ell,0}^{u,0}$ provides two algebraic restrictions on the representation $x$. The first is that the type of the Harder-Narasimhan filtration is determined by the critical set $C_\ell$, since the graded object of this filtration is isomorphic to that of the limit of the downwards flow. The limit of the upwards flow determines a second condition, which is that the representation $x$ must admit a different filtration, this time with increasing slopes and graded object determined by the the type of $C_u$, since $x$ is isomorphic to a representation in the negative slice of $C_u$ (cf. \cite[Lem. 3.24]{Wilkin17}).

\subsubsection{The image of the space of flow lines in the negative slice}\label{subsec:flow-lines-neg-slice}

By definition, the space of flow lines is contained in the unstable set $\mathcal{F}_{\ell,0}^{u,0} \subset W_{u,0}$. Via the homeomorphism $H : (W_u, W_{u,0}) \stackrel{\cong}{\rightarrow} (S_u, S_{u,0})$, we can study the image $H(\mathcal{F}_{\ell,0}^{u,0}) \subset S_{u,0}$, which turns out to be a more tractable object that we can describe explicitly and then transport back to the unstable set via the homeomorphism $H^{-1}$. In this section we do exactly that for the case where the stability parameter is the canonical central element of Definition \ref{def:aasp}.


In this case the lower critical point $x_\ell \in C_\ell$ decomposes as
\begin{equation}\label{eqn:lower-critical-decomp}
x_\ell = x_\ell^{(1)} + x_\ell^{(2)}, \quad x_\ell^{(1)} \in \Rep(Q, {\bf v}_\ell), \quad x_\ell^{(2)} \in \Rep(Q, {\bf v} - {\bf v}_\ell) .
\end{equation}
Note that $\slope_\alpha(Q, {\bf v}_\ell) < \slope_\alpha(Q, {\bf v} - {\bf v}_\ell)$. Moreover, $x_\ell^{(2)}$ is isomorphic to the graded object of the Jordan-H\"older filtration of the maximal semistable subrepresentation of $x$ and $x_\ell^{(1)}$ is isomorphic to the quotient of $x$ by its maximal semistable subrepresentation.

Similarly, the upper critical point $x_u \in C_u$ decomposes as
\begin{equation}\label{eqn:upper-critical-decomp}
x_u = x_u^{(1)} + x_u^{(2)}, \quad x_u^{(1)} \in \Rep(Q, {\bf v}_1 - {\bf d}), \quad x_\ell^{(2)} \in \Rep(Q, {\bf v} - {\bf v}_1 + {\bf d}) .
\end{equation}
Again we have $\slope_\alpha(Q, {\bf v}_\ell - {\bf d}) < \slope_\alpha(Q, {\bf v} - {\bf v}_\ell + {\bf d})$ and that $x_u^{(1)}$ is stable, however this time $x_u^{(1)}$ is isomorphic to a stable subrepresentation of $x$, instead of a quotient representation. Since the two critical points are connected by a flow line, then \cite[Lem. 4.32]{Wilkin17} shows that $x_u^{(1)}$ is isomorphic to a subrepresentation of $x_\ell^{(1)}$. In particular, this implies that ${\bf d} > 0$. This is summarised in the following lemma.

\begin{lemma}\label{lem:reduce-structure-flow-line}
A flow line connecting $x_\ell$ and $x_u$ determines a reduction of structure group from $K_{\bf v}$ to $K_{{\bf v}_\ell - {\bf d}} \times K_{\bf d} \times K_{{\bf v} - {\bf v}_\ell}$ for which $x_\ell^{(1)} \in \Rep(Q, {\bf v}_\ell)$ and $x_u^{(1)} \in \Rep(Q, {\bf v}_\ell - {\bf d})$. The polystable representation $x_u^{(2)} \in \Rep(Q, {\bf v} - {\bf v_\ell} + {\bf d})$ then splits into subrepresentations
\begin{equation}\label{eqn:upper-crit-point-decompose}
x_u^{(2)} = x_u^{(2,1)} + x_u^{(2,2)} \in \Rep(Q, {\bf d}) \oplus \Rep(Q, {\bf v} - {\bf v_\ell}) . 
\end{equation}
The homeomorphism $H$ from the unstable set to the negative slice maps the representation $x$ to $\delta x \in S_{u,0}$.
\end{lemma}

Using Corollary \ref{cor:reduce-to-zero}, from now on we restrict to the case where $x_\ell^{(2)} = 0$. This condition imposes the following restriction on $x \in \mathcal{F}_{\ell,0}^{u,0}$.

\begin{lemma}\label{lem:isomorphic-to-lower-limit}
If $x_\ell^{(2)} = 0$, then $x$ is isomorphic to $x_{\ell}^{(1)}$, which is then stable as a representation in $\Rep(Q, {\bf v_\ell})$.
\end{lemma}

\begin{proof}
Using Theorem \cite[Thm. 4.22]{Wilkin17}, the representation $x$ is isomorphic to a representation in $S_{x_u}^-$ for some $x_u = x_u^{(1)} + x_u^{(2)} \in C_u$.

Since $x$ flows down to $x_\ell$, then the Harder-Narasimhan type is determined, therefore the quotient by the maximal semistable subrepresentation, which we denote $x_\ell^{(1)} = x / x_\ell^{(2)}$ must have dimension vector ${\bf v_\ell}$. Moreover, \cite[Thm. 1.2]{Wilkin17} (the algebraic structure of the limit of the upwards flow) shows that $x_u^{(1)}$ is a stable subrepresentation of $x$ which contains the vertex $\infty$, and so it must be a subrepresentation of $x_\ell^{(1)}$. The representation $x_u^{(2)}$ satisfies $\mu(x_u^{(2)}) = 0$ by \cite[Prop. 3.13]{Wilkin17} and therefore $x_u^{(2)}$ is semisimple. With respect to these subrepresentations, the negative slice decomposes as 
\begin{equation*}
S_u^- \subset \Hom^1(Q, {\bf v} - {\bf v}_\ell + {\bf d}, {\bf v}_\ell - {\bf d}) \cong \Hom^1(Q, {\bf d}, {\bf v}_\ell - {\bf d}) \oplus \Hom^1(Q, {\bf v} - {\bf v}_\ell, {\bf v}_\ell - {\bf d}) .
\end{equation*}
With respect to this decomposition, write $\delta x = \delta x_1 + \delta x_2$ and $x_u^{(2)} = x_u^{(2,d)} + x_u^{(2,v_\ell-d)}$. Therefore, $x$ is isomorphic to $x_u + \delta x_1 + \delta x_2$ and $x_u^{(1)} + \delta x_1 + x_u^{(2,d)}$ is the quotient by the maximal semistable subrepresentation. The choice of stability parameter then implies that $x_u^{(1)} + \delta x_1 + x_u^{(2,d)}$ is stable (see Remark \ref{rem:canonical-central-element}), and therefore it is isomorphic to $x_\ell^{(1)}$. Therefore
\begin{equation*}
x \cong x_u^{(1)} + \delta x_1 + x_u^{(2,d)} + \delta x_2 + x_u^{(2,v_\ell-d)} = x_\ell^{(1)} + \delta x_2 + x_u^{(2,v_\ell-d)} .
\end{equation*}

Since $\delta x_2 \in S_u^- \subset \Hom^1(Q, {\bf v} - {\bf v}_\ell, {\bf v}_\ell - {\bf d})$, then \cite[Lem. 3.26]{Wilkin17} shows that if $\delta x_2 \neq 0$, then $x_\ell^{(1)} + \delta x_2$ is isomorphic to a non-zero element of the negative slice $S_{\ell,0}$ of the lower critical point, and therefore cannot flow down to $x_\ell$ since the energy satisfies $\| \mu(x_\ell^{(1)} + \delta x_2) - \alpha \|^2 < \| \mu(x_\ell) - \alpha \|^2$. Therefore we must have $\delta x_2 = 0$. The condition $x_\ell^{(2)} = 0$ implies that $x_u^{(2,v_\ell-d)}$, and so $x \cong x_u^{(1)} + \delta x_1 + x_u^{(2,d)} \cong x_\ell^{(1)}$.
\end{proof}

This setup can be represented by the following diagram. The notation $Q_{\bf v_\ell}$ means the quiver $Q$ with dimension vector ${\bf v_\ell}$, an arrow from $Q_{\bf v_d}$ to $Q_{\bf v_{\ell-d}}$ represents a homomorphism in the space $\Hom^1(Q, {\bf v_d}, {\bf v_{\ell-d}})$, a special case of which is a loop from $Q_{\bf v_\ell}$ to itself which corresponds to a representation in $\Rep(Q, {\bf v_\ell})$. At the upper and lower critical points, the arrows have been labelled with the corresponding subrepresentations from \eqref{eqn:lower-critical-decomp} and \eqref{eqn:upper-critical-decomp}.
\begin{equation*}
\begin{tikzcd}[row sep = 2cm, column sep = 1cm]
x_u & \bullet_{Q_{{\bf v_\ell}-{\bf d}}} \arrow[loop, swap, "x_u^{(1)}"] & & & \bullet_{Q_{{\bf v} - {\bf v_\ell} + {\bf d}}} \arrow[loop, swap, "x_u^{(2)}"] \\
x & \bullet_{Q_{{\bf v_\ell}-{\bf d}}} \arrow[loop, swap, "x_u^{(1)}"] & & \bullet_{Q_{{\bf d}}} \arrow[loop, swap, "x_u^{(2,d)}"] \arrow[bend right = 15]{ll}{\delta x_1} & & \bullet_{Q_{{\bf v} - {\bf v_\ell}}} \arrow[loop, swap, "x_u^{(2,v_\ell-d)}"] \arrow[bend left=10]{llll}{\delta x_2=0} \\
x_\ell & & \bullet_{Q_{{\bf v_\ell}}} \arrow[loop, swap, "x_\ell^{(1)}"] & & & \bullet_{Q_{{\bf v} - {\bf v_\ell}}} \arrow[loop, swap, "x_\ell^{(2)}=0"]
\end{tikzcd}
\end{equation*}

\begin{corollary}\label{cor:flow-lines-subspace}
After applying the reduction of structure group from $K_{\bf v}$ to $K_{{\bf v}_\ell - {\bf d}} \times K_{\bf d} \times K_{{\bf v} - {\bf v}_\ell}$, the image $H(\mathcal{F}_{\ell,0}^{u,0}) \subset S_{u,0}^-$ is contained in the subspace $S_{u,0}^- \cap \Hom^1(Q, {\bf d}, {\bf v}_\ell - {\bf d})$.
\end{corollary}

\subsubsection{Adjacent critical sets}

Corollary \ref{cor:flow-lines-subspace} shows that the image of the flow lines $H(\mathcal{F}_{\ell,0}^{u,0}) \subset S_{u,0}^-$ is contained in a given subspace of the negative slice that can be described in terms of the algebraic data (Harder-Narasimhan type) of the upper and lower critical sets. In general, the image may not be equal to this subspace, due to the possible existence of flow lines that converge to some intermediate critical set.

In this section we consider the special case where there are no intermediate critical sets, and so the image $H(\mathcal{F}_{\ell,0}^{u,0}) \subset S_{u,0}$ is equal to the subspace identified in Corollary \ref{cor:flow-lines-subspace}. We can then explicitly describe a tubular neighbourhood of $H(\mathcal{F}_{\ell,0}^{u,0}) \subset S_{u,0}$ as a disk bundle over $H(\mathcal{F}_{\ell,0}^{u,0})$. 

\begin{definition}
Two critical sets $C_\ell$ and $C_u$ with corresponding Harder-Narasimhan types ${\bf v}_u$ and ${\bf v}_\ell$ are called \emph{adjacent} iff ${\bf v}_\ell = {\bf v}_u + {\bf e}_k$ for some $k \in \mathcal{I}$.
\end{definition}

In the notation of Corollary \ref{cor:flow-lines-subspace}, this is equivalent to ${\bf d} = {\bf e}_k$. Once again we can decompose $x_\ell = x_\ell^{(1)} + x_\ell^{(2)} \in C_\ell$, $x_u = x_u^{(1)} + x_u^{(2)} \in C_u$  and note that (for the purposes of computing the relative cohomology groups \eqref{eqn:reduce-relative-cohomology}) Corollary \ref{cor:reduce-to-zero} allows us to restrict attention to the case where $x_\ell^{(2)} = 0$ and $x_u^{(2)} = 0$.

Since $x_u^{(2)} = 0$ and $x_u^{(1)}$ is $\alpha$-stable, then the isotropy group of $x_u$ is $\U(1) \times K_{{\bf v} - {\bf v_u}}$, which then acts on the negative slice $S_{x_u}^-$. Recall that the diagonal $\U(1) \subset K_{\bf v}$ acts trivially on all of $\Rep(Q, {\bf v})$, and therefore the action of the isotropy group on the slice is induced from the action of $K_{{\bf v} - {\bf v_u}}$ on $\Hom^1(Q, {\bf v} - {\bf v_u}, {\bf v_u})$.



\begin{lemma}\label{lem:adjacent-all-subspace}
If $C_\ell$ and $C_u$ are adjacent critical sets, then for every $x_u = x_u^{(1)} + x_u^{(2)} \in C_u$ with $x_u^{(2)} = 0$ and every representation $\delta x \in S_{u,0} \subset \Hom^1(Q, {\bf v} - {\bf v_u}, {\bf v_u})$ such that kernel of $\delta x$ has dimension vector ${\bf v} - {\bf v_u} - {\bf e_k} = {\bf v} - {\bf v_\ell}$, the limit of the downwards flow with initial condition $x_u + \delta x$ is contained in $C_\ell$. Conversely, every element of the negative slice which flows down to $x_\ell = x_\ell^{(1)} + x_\ell^{(2)} \in C_\ell$ with $x_\ell^{(2)} = 0$ must have this form.
\end{lemma}

\begin{proof}
By the choice of stability parameter, the stability of $x_u^{(1)} \in \Rep(Q, {\bf v_u})$ and the fact that $x_u^{(2)} = 0$ imply that the maximal semistable subrepresentation of $x_u + \delta x$ has dimension vector equal to the kernel of $\delta x$. Therefore $x_u + \delta x$ flows down to $x_\ell \in C_\ell$ if and only if the dimension vector of the kernel of $\delta x$ is equal to ${\bf v} - {\bf v_u} - {\bf e_k} = {\bf v} - {\bf v_\ell}$.
\end{proof}

If $C_\ell$ and $C_u$ are adjacent critical sets then Lemma \ref{lem:reduce-structure-flow-line} shows that there is a reduction of structure group to $K_{\bf v_u} \times \mathsf{U}(1) \times K_{{\bf v} - {\bf v_\ell}}$. Modulo this reduction, there is an explicit description of the space of flow lines as a subspace of the negative slice of the upper critical set.

\begin{corollary}\label{cor:adjacent-vector-space}
Let $C_\ell$ and $C_u$ be adjacent critical sets and restrict to the subset of upper critical points of the form $x_u = x_u^{(1)} + x_u^{(2)} \in C_u$ with $x_u^{(2)} = 0$ and lower critical points of the form $x_\ell = x_\ell^{(1)} + x_\ell^{(2)} \in C_\ell$ with $x_\ell^{(2)} = 0$. Then modulo the reduction of structure group to $K_{\bf v_u} \times \mathsf{U}(1) \times K_{{\bf v} - {\bf v_\ell}}$, the space of flow lines $\mathcal{F}_{\ell,0}^{u,0}$ fibres over $C_u$ such that the fibre over $x_u \in C_u$ is isomorphic to $V \setminus \{0\}$, where $V$ is the cokernel of
\begin{equation}\label{eqn:flow-lines-cokernel}
\begin{tikzcd}
\Hom^0(Q, {\bf e_k}, {\bf v_u}) \arrow{r}{\rho_{x_u}^{\mathbb{C}}} & \Hom^1(Q, {\bf e_k}, {\bf v_u}) .
\end{tikzcd}
\end{equation}
Moreover, the subgroup $\mathsf{U}(1) \subset K_{\bf v_u} \times \mathsf{U}(1) \times K_{{\bf v} - {\bf v_\ell}}$ acts freely on these fibres with weight one.
\end{corollary}

\begin{proof}
This follows from Lemma \ref{lem:adjacent-all-subspace} and the explicit description of Lemma \ref{lem:reduce-structure-flow-line}.
\end{proof}


Therefore we can describe tubular neighbourhoods of the space of flow lines as follows.

\begin{lemma}\label{lem:flow-lines-tubular-nbhd}
Given two adjacent critical sets $C_\ell$ and $C_u$ with Harder-Narasimhan types satisfying ${\bf v}_\ell = {\bf v}_u + {\bf e}_k$, the fibres of the tubular neighbourhood of $H(\mathcal{F}_{\ell,0}^{u,0}) \subset S_{u,0}$ at a given $\delta x \in H(\mathcal{F}_{\ell,0}^{u,0})$ are isomorphic to the cokernel of
\begin{equation}\label{eqn:normal-to-flow-lines}
\begin{tikzcd}
\Hom^0(Q, {\bf v} - {\bf v}_u - {\bf e}_k,  {\bf v}_u) \arrow{r}{\rho_{x_u}^\mathbb{C}} \arrow[draw=none]{d}[scale=1.2]{\bigoplus} & \Hom^1(Q,{\bf v} - {\bf v}_u - {\bf e}_k,  {\bf v}_u) \\
\Hom^0(Q, {\bf v} - {\bf v_u} - {\bf e}_k, {\bf e}_k) \arrow{ur}[swap]{\rho_{\delta x}^\mathbb{C}}
\end{tikzcd}
\end{equation}
which is a direct summand of the negative slice at $x_u$.
\end{lemma}

\begin{proof}
The negative slice of the upper critical point $x_u$ corresponds to the cokernel of
\begin{equation*}
\Hom^0(Q, {\bf v} - {\bf v}_u,  {\bf v}_u) \stackrel{\rho_{x_u}^\mathbb{C}}{\xrightarrow{\hspace{1 cm}}} \Hom^1(Q,{\bf v} - {\bf v}_u,  {\bf v}_u) .
\end{equation*}
Within this slice, Lemma \ref{lem:adjacent-all-subspace} shows that the space of flow lines to $C_\ell$ is the Grassmann submanifold of all representations $\delta x \in \Hom^1(Q,{\bf v} - {\bf v}_u,  {\bf v}_u)$ for which the kernel has dimension vector ${\bf v} - {\bf v}_u - {\bf e_k}$. This space corresponds to fixing a reduction of structure group of the isotropy group $K_{{\bf v} - {\bf v_u}}$ of $x_u$ to the subgroup $K_{{\bf v}-{\bf v}_u - {\bf e_k}} \times K_{\bf e_k}$ so that $\delta x \in \Hom^1(Q, {\bf e_k}, {\bf v_u})$. 

Given a fixed choice of reduction of structure group, the fibres of the normal bundle correspond to the cokernel of
\begin{equation*}
\begin{tikzcd}
\Hom^0(Q, {\bf v} - {\bf v}_u - {\bf e}_k,  {\bf v}_u) \arrow{r}{\rho_{x_u}^\mathbb{C}} & \Hom^1(Q, {\bf v} - {\bf v}_u - {\bf e}_k,  {\bf v}_u) ,
\end{tikzcd}
\end{equation*}
and varying this choice requires finding the normal bundle of
\begin{equation*}
K_{{\bf v}-{\bf v_u}} \times_{(K_{{\bf v}-{\bf v}_u - {\bf e_k}} \times K_{\bf e_k})} \Hom^1(Q, {\bf e_k}, {\bf v_u}) \subset S_{x_u}^- .
\end{equation*}
which corresponds to the cokernel of \eqref{eqn:normal-to-flow-lines}, consisting of all representations orthogonal to the tangent space of above fibre product.
\end{proof}


Applying the homeomorphism $H : (W_u, W_{u,0}) \stackrel{\cong}{\rightarrow} (S_u, S_{u,0})$ then determines a tubular neighbourhood $U_\ell \rightarrow \mathcal{F}_{\ell,0}^{u,0} \subset W_{u,0}$, which is homeomorphic to a disk bundle.

Now the lower critical point (which is isomorphic to $x = x_u + \delta x$ by Lemma \ref{lem:isomorphic-to-lower-limit}) lies inside the associated Harder-Narasimhan stratum, which also has a tubular neighbourhood inside the total space $\Rep(Q, {\bf v})$ whose fibres are given by the cokernel of the following complex
\begin{equation*}
\begin{tikzcd}
\Hom^0(Q, {\bf v_\ell}, {\bf v} - {\bf v_\ell}) \arrow{r}{\rho_{x_\ell}^\mathbb{C}} & \Hom^1(Q, {\bf v_\ell}, {\bf v} - {\bf v_\ell}) .
\end{tikzcd}
\end{equation*}
With respect to the reduction of structure group from $K_{\bf v}$ to $K_{\bf v_u} \times K_{{\bf e}_k} \times K_{{\bf v} - {\bf v_u} - {\bf e}_k}$, this complex splits
\begin{equation}\label{eqn:normal-to-lower-stratum}
\begin{tikzcd}
\Hom^0(Q, {\bf v} - {\bf v}_\ell,  {\bf v}_u) \arrow{r}{\rho_{x_u}^\mathbb{C}} \arrow[draw=none]{d}[scale=1.2]{\bigoplus} & \Hom^1(Q,{\bf v} - {\bf v}_\ell,  {\bf v}_u) \arrow[draw=none]{d}[scale=1.2]{\bigoplus} \\
\Hom^0(Q, {\bf v} - {\bf v}_\ell, {\bf e}_k) \arrow{ur}[swap]{\rho_{\delta x}^\mathbb{C}} & \Hom^1(Q, {\bf v} - {\bf v}_\ell, {\bf e}_k) 
\end{tikzcd}
\end{equation}
Note that since $x_\ell$ is isomorphic to $x = x_u + \delta x$ by Lemma \ref{lem:isomorphic-to-lower-limit}, then the image of $\rho_{x_\ell}^\mathbb{C}$ is orthogonal to $\Hom^1(Q, {\bf v} - {\bf v}_\ell, {\bf e}_k)$. Therefore, since ${\bf v}_\ell = {\bf v}_u + {\bf e}_k$, then the above complex contains \eqref{eqn:normal-to-flow-lines} as a direct summand. 

Since $\mathcal{F}_{\ell,0}^{u,0}$ is also contained in the Harder-Narasimhan stratum of the critical set $C_\ell$, one can then consider a tubular neighbourhood of the stratum and restrict it to a disk bundle $\pi : N_\ell \rightarrow \mathcal{F}_{\ell,0}^{u,0}$. Given $x \in \mathcal{F}_{\ell,0}^{u,0}$, the fibre $\pi^{-1}(x)$ corresponds to a neighbourhood of the origin in the cokernel of
\begin{equation}\label{eqn:pullback-normal-stratum}
\begin{tikzcd}
\Hom^0(Q, {\bf v} - {\bf v}_\ell, {\bf v}_\ell) \arrow{r}{\rho_x^\mathbb{C}} & \Hom^1(Q, {\bf v} - {\bf v}_\ell, {\bf v}_\ell) .
\end{tikzcd}
\end{equation}

Therefore we have two disk bundles over $\mathcal{F}_{\ell,0}^{u,0}$. 
\begin{enumerate}

\item A tubular neighbourhood $D \rightarrow \mathcal{F}_{\ell,0}^{u,0}$ inside the unstable set $W_{C_u}^-$ with fibres given by the cokernel of \eqref{eqn:normal-to-flow-lines}, and

\item a tubular neighbourhood of the Harder-Narasimhan stratum of $C_\ell$ inside the ambient manifold $\Rep(Q, {\bf v})$, which restricts to a disk bundle $V \rightarrow \mathcal{F}_{\ell,0}^{u,0}$ with fibres given by the cokernel of \eqref{eqn:normal-to-lower-stratum}.

\end{enumerate}

In the following we will abuse the notation by referring to the disk bundle and the associated vector bundle with the same notation. The direct sum decomposition of \eqref{eqn:normal-to-lower-stratum} shows that $D$ is a subbundle of $V$. This is summarised in the following proposition.

\begin{proposition}\label{prop:quivers-T2}
The Harder-Narasimhan stratum of $C_\ell$ has a tubular neighbourhood, which restricts to a disk bundle $V$ over $\mathcal{F}_{\ell,0}^{u,0}$ with fibres given by the cokernel of \eqref{eqn:normal-to-lower-stratum}. The space of flow lines $\mathcal{F}_{\ell,0}^{u,0}$ has a tubular neighbourhood $D$ in $W_{C_u}^-$, for which the fibres are given by the cokernel of \eqref{eqn:normal-to-flow-lines}. Moreover, after shrinking the above neighbourhoods if necessary, the bundle $D$ is a subbundle of $V$.

After reducing the structure group from $K_{\bf v}$ to $K_{{\bf v}_\ell - {\bf e}_k} \times K_{{\bf e}_k} \times K_{{\bf v} - {\bf v}_\ell}$ (cf. Lemma \ref{lem:reduce-structure-flow-line}), the quotient bundle $T := V/D$ is trivial with fibres $\Hom^1(Q, {\bf v} - {\bf v}_\ell, {\bf e}_k)$. 
\end{proposition}


\subsubsection{The equivariant Euler class of the quotient of the tubular neighbourhoods}

The above proposition shows that the tubular neighbourhood of the lower stratum restricts to a disk bundle $V \rightarrow \mathcal{F}_{\ell,0}^{u,0}$, which contains the normal bundle $D \rightarrow \mathcal{F}_{\ell,0}^{u,0}$ inside the unstable set $W_{C_u}^-$ as a subbundle. We would now like to compute the equivariant Euler class of the quotient bundle, which will then be used in the construction of the differential and cup product on the Morse complex in Section \ref{subsec:cup-product-T1-T2}. 

After applying the reduction of structure group from $K_{\bf v}$ to $K_{\bf v_u} \times K_{{\bf e}_k} \times K_{{\bf v} - {\bf v_u} - {\bf e_k}}$ (cf. Lemma \ref{lem:reduce-structure-flow-line}), there is a trivial bundle $T$ over $\mathcal{F}_{\ell,0}^{u,0}$ with fibres given by $\Hom^1(Q, {\bf v} - {\bf v}_\ell, {\bf e}_k) $. Since the fibres of $V$ are given by the cokernel of \eqref{eqn:normal-to-lower-stratum} and the fibres of $D$ are given by the cokernel of \eqref{eqn:normal-to-flow-lines}, which appears as a direct summand in \eqref{eqn:normal-to-lower-stratum}, then $V \cong D \oplus T$.

For the construction of the differential and cup product in Section \ref{sec:transversality}, we need to compute the equivariant Euler class of the bundle $T$. The following lemma applies in this setting (cf. \cite{AtiyahBott84}).

\begin{lemma}\label{lem:equiv-Euler-general}
Let $X$ be a connected topological space and $\pi : E \rightarrow X$ a trivial complex vector bundle of rank $n$ with a Hermitian metric on the fibres. Suppose that a compact Lie group $K$ acts fibrewise on $E$ via a faithful representation $\rho : K \rightarrow \mathsf{U}(n)$, so that $K$ fixes the zero section $X \hookrightarrow E$ and acts on each fibre $\pi^{-1}(x)$ by $k \cdot v = \rho(k) v$. Then the $K$-equivariant Euler class of $E$ is
\begin{equation*}
e_K = 1 \otimes \rho^*(c_n) \in H_K^*(X) \cong H^*(X) \otimes H^*(BK) ,
\end{equation*}
where $1$ represents the identity element in $H^0(X) \cong \mathbb{C}$, the projection $BK \rightarrow \mathsf{BU}(n)$ induces $\rho^* : H^*(\mathsf{BU}(n)) \rightarrow H^*(BK)$ and $c_n$ is the degree $2n$ generator of the polynomial ring $\mathbb{C}[c_1, \ldots, c_n] \cong H^*(\mathsf{BU}(n))$.
\end{lemma}

Now we can apply this result to the tubular neighbourhood of Proposition \ref{prop:quivers-T2}.

\begin{corollary}\label{cor:Euler-class}
In the setting of Proposition \ref{prop:quivers-T2}, the disk $V \rightarrow \mathcal{F}_\ell^u$ contains the tubular neighbourhood $D \rightarrow \mathcal{F}_{\ell,0}^{u,0}$ as a subbundle, with trivial quotient bundle $T$. With respect to the reduction of structure group from $K_{\bf v}$ to $K_{\bf v_u} \times K_{{\bf e}_k} \times K_{{\bf v} - {\bf v_u} - {\bf e_k}}$, the subgroup $K_{{\bf e}_k} \times K_{{\bf v} - {\bf v_u} - {\bf e_k}} = K_{{\bf e}_k} \times K_{{\bf v} - {\bf v_\ell}}$ acts on the fibres $\Hom^1(Q, {\bf v} - {\bf v}_\ell, {\bf e}_k)$ of the quotient bundle $T$ via a faithful representation $\rho$, and the equivariant Euler class is $e_K = \rho^*(c_n)$, where $n = \dim_\mathbb{C} \Hom^1(Q, {\bf v} - {\bf v}_\ell, {\bf e}_k)$ and $c_n$ is the degree $2n$ generator of the polynomial ring $\mathbb{C}[c_1, \ldots, c_n] \cong H^*(\mathsf{BU}(n))$.
\end{corollary}

This will be used in Section \ref{sec:transversality}, where the Euler class of the quotient bundle $T$ is needed to compensate for the failure of transversality when constructing the differentials and cup product on the Morse complex (cf. \eqref{eqn:differential-T1-T2} and \eqref{eqn:cup-product-T1-T2}).

\section{Reduction of structure group and restricted quiver representations}\label{sec:restricted-quivers}

The results of Section \ref{sec:convolution-cup-product} use an embedding of the Hecke correspondence into the negative slice bundle on an ambient smooth manifold (cf. Corollary \ref{cor:thom-compatible}), or equivalently (by \cite{Wilkin17}) an embedding of a space of flow lines between adjacent critical sets. In order to describe the Thom class of this embedding, it is necessary to develop the deformation theory associated to this problem, which is the goal of this section (cf. Lemmas \ref{lem:smoothness-restricts} and \ref{lem:stable-implies-adjoint-injective}). Much of the theory works in more generality than the case of Nakajima quivers, and the majority of this section is taken up with the task of defining some general conditions on the quiver and the relations for which this theory works and then illustrating this with examples. The reader who is only interested in Nakajima quiver varieties can skip to the next section.

First we describe a general construction to reduce the structure group of a quiver representation by decomposing the vector spaces at each vertex and then eliminating some edges from the new quiver representation. Some examples of this are: handsaw quivers (cf. \cite{Nakajima12} and Example \ref{ex:handsaw-reduction}), the fixed points of a circle action on the moduli space (cf. \cite{Nakajima04}, \cite{Hitchin87} and Example \ref{ex:fixed-circle-action}) and the quiver associated to the negative slice (cf. Example \ref{ex:neg-slice-restricted} below). The goal here is to prove Lemma \ref{lem:smoothness-restricts} and Lemma \ref{lem:stable-implies-adjoint-injective}, which are general results showing that certain properties of the original representation and the derivative of the relation map are preserved after reducing the structure group in this way. 

The following constructions are needed to prove these results for a general class of quivers. Throughout the reader is referred to Examples \ref{ex:handsaw-reduction} and \ref{ex:fixed-circle-action} which illustrate this theory.

Given a quiver $Q = (\mathcal{I}, \mathcal{E})$, consider a finite set of vertices $\mathcal{I}'$ with a surjective map $v : \mathcal{I}' \rightarrow \mathcal{I}$. There is an induced set of edges $\tilde{\mathcal{E}}'$ constructed as follows. For each $a \in \mathcal{E}$ and each $t_i \in v^{-1}(t(a))$, $h_j \in v^{-1}(h(a))$, define an edge $a_{ij} \in \tilde{\mathcal{E}}'$ such that $t(a_{ij}) = t(a)$ and $h(a_{ij}) = h(a)$. Note that the construction determines a surjective map $e : \tilde{\mathcal{E}}' \rightarrow \mathcal{E}$. This data defines a new quiver $\tilde{Q}' := (\mathcal{I}', \tilde{\mathcal{E}}')$, which we call an \emph{expansion} of $Q$.

Given a dimension vector ${\bf v}$ for $Q$, together with vector spaces $\{ V_k \}_{k \in \mathcal{I}}$, decompose each $V_k \cong V_{k_1} \oplus \cdots \oplus V_{k_{m_k}}$ where $m_k = \# (v^{-1}(k))$ and each $V_{k_i}$ has positive dimension. Note that this last condition bounds the number of vertices in $v^{-1}(k)$; in particular if $\dim_\mathbb{C} V_k = 1$ then $\#(v^{-1}(k)) = 1$. This data defines a dimension vector ${\bf v'}$ for the quiver $\tilde{Q}'$. The group $K_{\bf v'}$ acting on $\Rep(Q, {\bf v'})$ is a reduction of structure group of the original action of $K_{\bf v}$ on $\Rep(Q, {\bf v})$.

A representation $x' \in \Rep(\tilde{Q}', {\bf v'})$ then induces a representation $x \in \Rep(Q, {\bf v})$ by using the direct sum $V_k \cong V_{k_1} \oplus \cdots \oplus V_{k_{m_k}}$, and the converse is true since the edge set $\tilde{\mathcal{E}}'$ includes every possible edge $a_{ij}$ mapping to each $a \in \mathcal{E}$. Therefore there is an isomorphism
\begin{equation}\label{eqn:isomorphism-decomposed-reps}
\Rep(\tilde{Q}', {\bf v'}) \stackrel{\cong}{\longrightarrow} \Rep(Q, {\bf v}) 
\end{equation}
as well as an inclusion of the Lie algebra of the structure group $\Hom^0(\tilde{Q}', {\bf v'}) \hookrightarrow \Hom^0(Q, {\bf v})$.

Now choose a subset $\mathcal{E}' \subset \tilde{\mathcal{E}'}$ and define a new quiver $Q' = (\mathcal{I}', \mathcal{E}')$, which we call a \emph{restriction} of $\tilde{Q}'$. Given a dimension vector ${\bf v'}$ as above, there is an isomorphism $\Hom^0(Q', {\bf v'}) \stackrel{\cong}{\longrightarrow} \Hom^0(\tilde{Q}', {\bf v'})$ and an injective homomorphism
\begin{equation}\label{eqn:inclusion-decomposed-reps}
S : \Rep(Q', {\bf v'}) \hookrightarrow \Rep(\tilde{Q}', {\bf v'}) \cong \Rep(Q, {\bf v}) . 
\end{equation}


A set of relations $\mathcal{R}$ for the original quiver $Q$ then induces a set of relations $\mathcal{R}'$ for $Q'$. We will need a precise description of this for Lemma \ref{lem:smoothness-restricts}, and so the details of this process are as follows.
  
For each path $p = a_n \cdots a_1$ in the original quiver $Q$, there is a collection $\tilde{\mathcal{P}}'$ of paths in an expansion $\tilde{Q}'$ given by $p_{i_1, \ldots, i_n} = a_n^{(i_n)} \cdots a_1^{(i_1)}$ along edges $a_k^{(i_k)} \in e^{-1}(a_k)$ for each $k = 1, \ldots, n$.  Note that a restriction $Q'$ has a subset of paths $\mathcal{P}' \subset \tilde{\mathcal{P}}'$ given by taking only considering paths with edges in $\mathcal{E}' \subset \tilde{\mathcal{E}}'$ and so $\tilde{\mathcal{P}}' \setminus \mathcal{P}'$ consists of all paths with at least one edge in $\tilde{\mathcal{E}}' \setminus \mathcal{E}'$. We abuse the notation and use $e$ to denote both of the projections $e : \mathcal{P}' \rightarrow \mathcal{P}$ and $e : \tilde{\mathcal{P}}' \rightarrow \mathcal{P}$.

Therefore the relations $\mathcal{R}$ for $Q$ determine a new set of relations $\tilde{\mathcal{R}}'$ for $\tilde{Q}'$ given as follows. For each relation 
\begin{equation*}
r  = \sum_{p \in \mathcal{P}_{t(r), h(r)}} \lambda_p p \in \mathcal{R}
\end{equation*}
with head and tail $h(r), t(r) \in \mathcal{I}$ and for each choice of $h(r') \in v^{-1}(h(r))$ and $t(r') \in v^{-1}(t(r))$ there is a relation $r'$ defined by
\begin{equation*}
r' = \sum_{p' \in \tilde{\mathcal{P}}_{t(r'), h(r')}'} \lambda_{v(p')} p' .
\end{equation*}
Define $\tilde{\mathcal{R}}'$ to be the set of all such relations. For a restriction $Q'$, we consider the subset $\mathcal{R}' \subset \tilde{\mathcal{R}}'$ of relations for which all of the paths are in the subset of edges $\mathcal{E}' \subset \tilde{\mathcal{E}}'$. Equivalently, we construct $\mathcal{R}'$ by removing paths from $\tilde{\mathcal{R}}'$, where a path appearing nontrivially in some relation $r' \in \tilde{\mathcal{R}}'$ is removed if and only if it contains at least one edge in $\tilde{\mathcal{E}}' \setminus \mathcal{E}'$. This construction is illustrated in Examples \ref{ex:handsaw-reduction} and \ref{ex:fixed-circle-action}.

Given a representation $x' \in \Rep(Q', {\bf v'}, \mathcal{R}')$ mapping to $x = S(x') \in \Rep(Q, {\bf v}, \mathcal{R})$ as above, a subrepresentation of $x'$ then induces a subrepresentation of $x$. Note that the converse is not necessarily true, since a subrepresentation for $x$ may not be compatible with the decomposition of each $V_k \cong V_{k_1} \oplus \cdots \oplus V_{k_{m_k}}$.


A stability parameter $\alpha$ for $\Rep(Q, {\bf v})$ induces a parameter $\alpha'$ for $\Rep(Q, {\bf v'})$ by pullback $\alpha' = \alpha \circ v$. Note that if a vertex $\infty \in \mathcal{I}$ has dimension one with respect to ${\bf v}$, then $v^{-1}(\infty)$ consists of a single vertex which has dimension one with respect to ${\bf v'}$, therefore the canonical stability parameter for $(Q, {\bf v})$ from Definition \ref{def:aasp} induces a canonical stability parameter for $(Q', {\bf v'})$. 

In particular, if $x' \in \Rep(Q', {\bf v'})$ induces a stable representation $x = S(x')  \in \Rep(Q, {\bf v})^{\alpha-st}$, then every subrepresentation of $x$ must contain the vertex $\infty$. Since any subrepresentation of $x'$ induces a subrepresentation of $x$ then each subrepresentation of $x'$ must also contain the vertex $v^{-1}(\infty)$, and therefore $x'$ is then $\alpha'$-stable. Therefore we have shown the following
\begin{lemma}\label{lem:stable-induces-stable}
If $\alpha$ is the stability parameter for $(Q, {\bf v})$ from Definition \ref{def:aasp} and $\alpha'$ the induced stability parameter on $(Q', {\bf v'})$, then $S(x') \in \Rep(Q, {\bf v}, \mathcal{R})^{\alpha-st}$ implies that $x' \in \Rep(Q, {\bf v'}, \mathcal{R}')^{\alpha'-st}$.
\end{lemma}

The proof of Lemma \ref{lem:smoothness-restricts} requires an extra condition on a restriction of quivers, namely that a path is removed if and only if it contains at least two edges in $\tilde{\mathcal{E}}' \setminus \mathcal{E}'$.

\begin{definition}\label{def:fully-restricted}
Let $Q' = (\mathcal{I}', \mathcal{E}')$ be a restriction of an expansion of $Q = (\mathcal{I}, \mathcal{E})$ and $\mathcal{R}$ a set of quadratic relations. Then the induced set of relations $\mathcal{R}'$ is called \emph{fully restricted} if and only if for each relation $r =  \sum_{p \in \mathcal{P}} \lambda_p p \in \mathcal{R}$ and each path $p' \in \tilde{\mathcal{P}}' \setminus \mathcal{P}'$ such that $\lambda_{v(p')} \neq 0$, there are at least two edges in $p'$ contained in $\tilde{\mathcal{E}}' \setminus \mathcal{E}'$.
\end{definition}

\begin{remark}\label{rem:fully-restricted}
The point of the above definition is to restrict the possible homomorphisms in the image of the derivative of the relation map \eqref{eqn:relation-derivative-def}. If a relation in $r' \in \mathcal{R}'$ corresponds to a homomorphism $V_{t(r')} \rightarrow V_{h(r')}$ then a fully restricted set of relations implies that $d\nu_x(\delta x)$ must be zero if $\delta x$ is a sum of homomorphisms along edges in $\tilde{\mathcal{E}}' \setminus \mathcal{E}'$. This is used in the proof of Lemma \ref{lem:smoothness-restricts} below.
\end{remark}


Since the relations $\mathcal{R}'$ are induced from $\mathcal{R}$ then there is an induced homomorphism of representations
\begin{equation}\label{eqn:inclusion-decomposed-relations}
S : \Rep(Q', {\bf v'}, \mathcal{R}') \hookrightarrow \Rep(\tilde{Q}', {\bf v'}, \tilde{\mathcal{R}}') \cong \Rep(Q, {\bf v}, \mathcal{R}) 
\end{equation}
given by \eqref{eqn:inclusion-decomposed-reps} and the same is true for the spaces $\Hom^0$, $\Hom^1$ and $\Rel$. Moreover, since the induced map of relations is injective, then there is a surjective splitting homomorphism 
\begin{equation*}
\begin{tikzcd}
\Rel(Q', {\bf v'}, {\bf v'}, \mathcal{R}') \arrow[hookrightarrow, shift left=3pt]{r} & \Rel(\tilde{Q}', {\bf v'}, {\bf v'}, \tilde{\mathcal{R}}') \arrow[twoheadrightarrow, shift left=3pt]{l}{q} 
\end{tikzcd}
\end{equation*}
for which $\ker(q)$ consists of all homomorphisms $\Hom(V_{t(r')}, V_{h(r')}) \subset \Rel(\tilde{Q}', {\bf v'}, {\bf v'}, \tilde{\mathcal{R}}')$ such that $r' \in \tilde{\mathcal{R}}' \setminus \mathcal{R}'$.

Then the deformation complexes \eqref{eqn:deformation-tangent-space} for $(\tilde{Q}', \tilde{\mathcal{R}}')$ and $(Q', \mathcal{R}')$ are related by
\begin{equation}\label{eqn:relate-restricted-deformations}
\begin{tikzcd}
\Hom^0(\tilde{Q}', {\bf v'}, {\bf v'}) \arrow{r}{\rho_{S(x')}^\mathbb{C}}& \Hom^1(\tilde{Q}', {\bf v'}, {\bf v'}) \arrow{r}{d \nu_{S(x')}} & \Rel(\tilde{Q}', {\bf v'}, {\bf v'}, \tilde{\mathcal{R}}') \arrow[twoheadrightarrow, shift left=5pt]{d}{q} \\
\Hom^0(Q', {\bf v'}, {\bf v'}) \arrow{r}{\rho_{x'}^\mathbb{C}} \arrow[hookrightarrow]{u} & \Hom^1(Q', {\bf v'}, {\bf v'}) \arrow{r}{d\nu_{x'}} \arrow[hookrightarrow]{u} & \Rel(Q', {\bf v'}, {\bf v'}, \mathcal{R}') \arrow[hookrightarrow, shift left=5pt]{u}
\end{tikzcd}
\end{equation}

\begin{remark}
Since the group $K_{\bf v}$ contains the one parameter subgroup of scalar multiples of the identity which act trivially on $\Rep(Q, {\bf v})$, then $x \in \Rep(Q, {\bf v})^{\alpha-st}$ only implies that $\rho_x^\mathbb{C} : \Hom^0(Q, {\bf v}, {\bf v}) \rightarrow \Hom^1(Q, {\bf v}, {\bf v})$ is injective on the trace-free part of $\Hom^0(Q, {\bf v}, {\bf v})$. Similarly, given a set $\mathcal{R}$ of relations, the map $d \nu_x : \Hom^1(Q, {\bf v}, {\bf v}, \mathcal{R})$ can only be surjective onto the trace free part (which is only nontrivial if there are relations $r \in \mathcal{R}$ with $t(r) = h(r)$). This is sufficient for the standard deformation theory to prove smoothness of the moduli space. To emphasise this, we will denote the trace-free part by
\begin{equation*}
\Hom_0^0(Q, {\bf v}, {\bf v}) \subset \Hom^0(Q, {\bf v}, {\bf v}) \quad \text{and} \quad \Rel_0(Q, {\bf v}, {\bf v}, \mathcal{R}) \subset \Rel(Q, {\bf v}, {\bf v}, \mathcal{R}) .
\end{equation*}
\end{remark}

\begin{lemma}\label{lem:smoothness-restricts}
Suppose that $Q'$ is defined via expansion and restriction of $Q$ as above, that the relations $\mathcal{R}'$ are fully restricted and that $x = S(x') \in \Rep(Q, {\bf v}, \mathcal{R})$ is induced from $x' \in \Rep(Q', {\bf v'}, \mathcal{R}')$. Suppose that $d \nu_x : \Hom^1(Q, {\bf v}, {\bf v}) \rightarrow \Rel_0(Q, {\bf v}, {\bf v}, \mathcal{R})$ is surjective. Then 
\begin{equation*}
d \nu_{x'} : \Hom^1(Q', {\bf v'}, {\bf v'}) \rightarrow \Rel_0(Q, {\bf v'}, {\bf v'}, \mathcal{R}')
\end{equation*}
is also surjective.
\end{lemma}

\begin{proof}
If $d \nu_x$ is surjective, then so is the same map for the expanded quiver $d \nu_x : \Hom^1(\tilde{Q}', {\bf v'}, {\bf v'}) \rightarrow \Rel_0(\tilde{Q}', {\bf v'}, {\bf v'}, \tilde{\mathcal{R}}')$, and so $q \circ d \nu_x$ is surjective. 

There are canonical splittings so that 
\begin{equation*}
\Hom^1(\tilde{Q}', {\bf v'}, {\bf v'}) \cong \Hom^1(Q', {\bf v'}, {\bf v'}) \oplus \Hom^1(\tilde{Q}', {\bf v'}, {\bf v'}) / \Hom^1(Q', {\bf v'}, {\bf v'}) ,
\end{equation*}
where the quotient consists of all  tangent vectors $\delta x \in \Hom(V_{t(a)}, V_{h(a)}) \subset \Hom^1(\tilde{Q}', {\bf v'}, {\bf v'})$ for edges $a \in \tilde{\mathcal{E}}' \setminus \mathcal{E}'$.

The condition that the relations $\mathcal{R}'$ are fully restricted means that all such tangent vectors are in the kernel of $q \circ d \nu_{S(x')}$. Therefore the image of $d \nu_{x'}$ is equal to the image of the composition
\begin{equation*}
\begin{tikzcd}
\Hom^1(Q', {\bf v'}, {\bf v'}) \arrow[hookrightarrow]{r} & \Hom^1(\tilde{Q}', {\bf v'}, {\bf v'}) \arrow{r}{d \nu_{S(x')}} & \Rel_0(\tilde{Q}', {\bf v'}, {\bf v'}, \tilde{\mathcal{R}}') \arrow{r}{q} & \Rel_0(Q', {\bf v}', {\bf v}', \mathcal{R}') .
\end{tikzcd}
\end{equation*}
and so surjectivity of $q$ and $d \nu_{S(x')}$ implies surjectivity of $d \nu_{x'}$.
\end{proof}

\begin{corollary}\label{cor:smoothness-restricts}
Let $Q$ be a quiver with dimension vector ${\bf v}$, and suppose that $d \nu_x : \Hom^1(Q, {\bf v}, {\bf v}) \rightarrow \Rel_0(Q, {\bf v}, {\bf v}, \mathcal{R})$ is surjective for all $x \in \Rep(Q, {\bf v})^{\alpha-st}$. Then if $Q'$ is defined via expansion and restriction of $Q$ and the relations $\mathcal{R}'$ are fully restricted, then the moduli space $\mathcal{M}_{\alpha'}(Q', {\bf v}, \mathcal{R}')$ is smooth.
\end{corollary}

\begin{proof}
Lemma \ref{lem:stable-induces-stable} shows that if $x' \in \Rep(Q', {\bf v'})$ is stable then $S(x')$ is stable and therefore $d \nu_{x'}$ is surjective by Lemma \ref{lem:smoothness-restricts}. The result follows from the standard deformation theory applied to the complex
\begin{equation*}
\begin{tikzcd}
\Hom^0(Q', {\bf v'}, {\bf v'}) \arrow{r}{\rho_{x'}^\mathbb{C}} & \Hom^1(Q', {\bf v'}, {\bf v'}) \arrow{r}{d\nu_{x'}} & \Rel_0(Q', {\bf v'}, {\bf v'}, \mathcal{R}') . 
\end{tikzcd} \qedhere
\end{equation*}
\end{proof}

The main example of interest is that of fully restricted expansions of Nakajima quiver varieties. Recall that for a Nakajima quiver variety, the hyperk\"ahler structure induces a duality on the deformation complex \eqref{eqn:deformation-tangent-space} such that $\rho_x^\mathbb{C}$ is injective if and only if $d \nu_x$ is surjective (cf. \cite{Hitchin87}, \cite[Lem. 3.10]{Nakajima98}). In particular, if $x$ is a stable representation then $\rho_x^\mathbb{C}$ is injective, and so the moduli space of stable representations is smooth (cf. \cite[Cor. 3.12]{Nakajima98}). The Lemmas \ref{lem:stable-induces-stable} and \ref{lem:smoothness-restricts} show that a representation $x' \in \Rep(Q', {\bf v}, \mathcal{R}')$ of a fully restricted expansion of a Nakajima quiver variety corresponds to a smooth point in the moduli space if it induces a stable representation $x = S(x') \in \Rep(Q, {\bf v}, \mathcal{R})$ for the original Nakajima quiver.

The following examples illustrate all of the above constructions.

\begin{example}[Handsaw quivers from the ADHM quiver]\label{ex:handsaw-reduction}
Consider the ADHM quiver with a given dimension vector
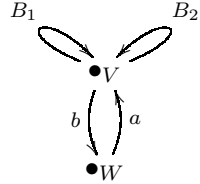
\begin{equation*}
\xymatrix{
\bullet_{V} \ar@`{(-15,5),(-8,10)}^{B_1} \ar@`{(15,5),(8,10)}_{B_2} \ar@/_0.5pc/[d]_{b} \\
\bullet_{W} \ar@/_0.5pc/[u]_{a} 
}
\end{equation*}
For notation, let $k$ denote the vertex with vector space $V$ and $\ell$ denote the vertex with vector space $W$. Define $\mathcal{I}' = \{ k_1, \ldots, k_{n-1}, \ell_1, \ldots, \ell_n \}$ with $v : \mathcal{I}' \rightarrow \mathcal{I}$ defined by $v(k_i) = k$ and $v(\ell_i) = \ell$. Now decompose $V \cong V_1 \oplus \cdots \oplus V_{n-1}$ and $W_1 \oplus \cdots \oplus W_n$. From the full set of all possible edges, choose a subset $\mathcal{E}'$ so that 
\begin{itemize}

\item $B_1$ decomposes into components $B_1^k : V_k \rightarrow V_{k+1}$ for $k = 1, \ldots, n-2$,

\item $B_2$ decomposes into components $B_2^k : V_k \rightarrow V_k$ for $k = 1, \ldots, n$,

\item $a$ decomposes into components $a_k : W_k \rightarrow V_k$ for $k = 1, \ldots, n-1$, and

\item $b$ decomposes into components $b_k : V_{k-1} \rightarrow W_k$ for $k = 2, \ldots, n$.

\end{itemize}
Then the above ADHM quiver expands into the handsaw quiver \eqref{eqn:handsaw-quiver} and the relations \eqref{eqn:hyperkahler-relation} decompose into the relations \eqref{eqn:handsaw-relation} for the handsaw quiver.

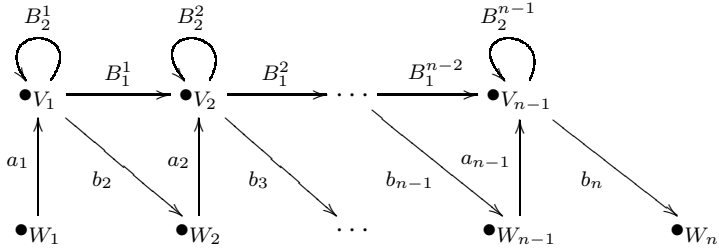
\begin{equation*}
\xymatrixrowsep{0.5in}
\xymatrixcolsep{0.5in}
\xymatrix{
\bullet_{V_1} \ar[r]^{B_1^1} \ar[dr]_{b_2} \ar@`{(10,10),(-10,10)}_{B_2^1} & \bullet_{V_2} \ar[r]^{B_1^2} \ar[dr]_{b_3} \ar@`{(30,10),(10,10)}_{B_2^2} & \cdots \ar[r]^{B_1^{n-2}} \ar[dr]_{b_{n-1}} & \bullet_{V_{n-1}} \ar[dr]_{b_n} \ar@`{(72,10),(52,10)}_{B_2^{n-1}} & \\
\bullet_{W_1} \ar[u]^{a_1} & \bullet_{W_2} \ar[u]^{a_2} & \cdots & \bullet_{W_{n-1}} \ar[u]^{a_{n-1}} & \bullet_{W_n} 
}
\end{equation*}
Note that there is one handsaw relation $V_k \rightarrow V_{k+1}$ for each $k = 1, \ldots, n-2$ (cf. \eqref{eqn:handsaw-relation}). In the original ADHM quiver with $V \cong V_1 \oplus \cdots \oplus V_{n-1}$ and $W_1 \oplus \cdots \oplus W_n$, there are many paths with a nontrival ADHM relation that map $V_k \rightarrow V_{k+1}$ (for example $B_1 : V_k \rightarrow V_\ell$ followed by $B_2 : V_\ell \rightarrow V_{k+1}$ for $\ell \neq k+1$), however we see that the \emph{only} such paths for which \emph{all} edges are in the handsaw quiver are those that appear in the handsaw relations. 

Moreover, the handsaw relations \eqref{eqn:handsaw-relation} are a \emph{full restriction} of the ADHM relations, since all the paths in the relations for the expanded ADHM quiver that do not appear in the handsaw relations must contain two edges that are not in the set of  handsaw quiver edges $\mathcal{E}'$.
\end{example}

\begin{example}[Fixed points of the circle action for Nakajima quivers]\label{ex:fixed-circle-action}
Recall that for a Nakajima quiver $Q$ with decomposition of edges denoted $\mathcal{E} = \mathcal{E}^{0,1} \cup \overline{\mathcal{E}^{0,1}}$, there is a well-studied circle action (cf. \cite{Nakajima04}, \cite{Hitchin87}) on the space of representations
\begin{equation*}
e^{i \theta} \cdot \left( (x_a, x_{\bar{a}}) \right)_{a \in \mathcal{E}^{0,1}} = \left( (x_a, e^{i \theta} x_{\bar{a}}) \right)_{a \in \mathcal{E}^{0,1}} .
\end{equation*}
This action commutes with the action of $G_{\bf v}$ and therefore descends to the moduli space. If an equivalence class $[x]$ is fixed by the circle action, then the representation $\Vect(Q, {\bf v})$ decomposes into weight spaces, which are ordered $\Vect(Q, {\bf v_1}) \oplus \cdots \oplus \Vect(Q, {\bf v_n})$. In this way we see that the fixed point determines an expansion of the quiver.

To simplify the notation in line with the previous example, we write $B_2 \in \Rep(Q, {\bf v})^{0,1}$ for the sum of all homomorphisms corresponding to edges $a \in \mathcal{E}^{0,1}$ and $B_1 \in \Rep(Q, {\bf v})^{1,0}$ for the sum of all homomorphisms corresponding to edges in $a \in \overline{\mathcal{E}^{0,1}}$. Then the fixed point condition (cf. \cite{Hitchin87} and \cite{Nakajima04}) implies that $B_1$ decomposes into the sum of $B_1^j \in \Rep(Q, {\bf v_j})^{1,0}$ and $B_2$ decomposes into the sum of $B_2^j \in \Hom^1(Q, {\bf v_j}, {\bf v_{j+1}})^{0,1}$. This is represented by the diagram below, where for simplicity a single edge is used to denote each representation $B_1^j \in \Rep(Q, {\bf v_j})^{1,0}$ and $B_2^j \in \Hom^1(Q, {\bf v_j}, {\bf v_{j+1}})^{0,1}$.
\begin{equation*}
\xymatrixrowsep{0.5in}
\xymatrixcolsep{0.5in}
\xymatrix{
\bullet_{(Q, {\bf v_1})} \ar[r]^{B_1^1} \ar@`{(10,10),(-10,10)}_{B_2^1} & \bullet_{(Q, {\bf v_2})} \ar[r]^{B_1^2} \ar@`{(35,10),(15,10)}_{B_2^2} & \cdots \ar[r]^{B_1^{n-2}} & \bullet_{(Q, {\bf v_n})} \ar@`{(82,10),(62,10)}_{B_2^{n-1}} 
}
\end{equation*}
Therefore we see that the fixed point also determines a restriction of the above expansion. The relations are then a full restriction of those for the original hyperk\"ahler quiver for the same reason as the previous example.
\end{example}

We would also like to apply the above theory to the negative slice, which will be used in Corollary \ref{cor:thom-compatible} to describe the Thom class of the normal bundle to the Hecke correspondence inside a certain projective bundle. Example \ref{ex:neg-slice-restricted} shows that the negative slice can be written as a space of representations for a restricted quiver, however the restriction does not satisfy the conditions of Lemma \ref{lem:smoothness-restricts}. If the quiver has no loops then it is quite easy to recover the result of Lemma \ref{lem:smoothness-restricts}, and if the quiver does have loops then we need to add an extra condition which is as follows.

\begin{definition}\label{def:loop-neg-slice}
Given a set $\mathcal{R}$ of quadratic relations, for each loop $a \in \mathcal{E}$ such that $t(a) = h(a) = t(r)$ for some relation $r$ and for all paths $p = ba$ with $\lambda_p(r) \neq 0$ there is a unique loop $a' \in \mathcal{E}$ such that $t(a') = h(a') = h(r)$ and a path $p' = a'b$ such that $\lambda_{p'}(r) \neq 0$.
\end{definition}

\begin{remark}\label{rem:ADHM-loops}
This is a general condition which is satisfied in useful examples. For the ADHM quiver, since (with the notation of Example \ref{ex:handsaw-reduction}) the relation
\begin{equation*}
r = B_1 B_2 - B_2 B_1 + ab
\end{equation*}
contains two loops $B_1$, $B_2$. We have $t(B_1) = t(r)$ and there is a path $p = B_2 B_1$ with $\lambda_p(r) = -1 \neq 0$. There is a unique loop (which in this case is $B_1$ again) such that $h(B_1) = h(r)$ and a path $p' = B_1 B_2$ such that $\lambda_{p'}(r) = 1 \neq 0$. The same idea applies to the loop $B_2$.

Similarly, for the handsaw quiver, each loop $B_2^k$ is at the tail of a relation 
\begin{equation*}
r_k = B_1^k B_2^k - B_2^{k+1} B_1^k + a_{k+1} b_{k+1}
\end{equation*}
for $k = 1, \ldots, n-2$. Let $p = B_1^k B_2^k$ be the path with $B_2^k$ at the tail, and note that there is a unique loop $B_2^{k+1}$ such that $p' = -B_2^{k+1} B_1^k$ with $\lambda_{p'}(r_k) = -1 \neq 0$. 
\end{remark}

\begin{example}[The negative slice as a restricted quiver]\label{ex:neg-slice-restricted}
Given a framed quiver $Q$ with dimension vector ${\bf v}$, consider a reduction of structure group corresponding to ${\bf v} = {\bf v_1} + {\bf v_2}$ with the following restriction
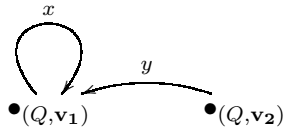
\begin{equation*}
\xymatrixrowsep{0.5in}
\xymatrixcolsep{0.5in}
\xymatrix{
\bullet_{(Q, {\bf v_1})} \ar@`{(-15,15),(15,15)}^{x}  & \bullet_{(Q, {\bf v_2})} \ar@/_0.8pc/[l]_{y} 
}
\end{equation*}
with $x \in \Rep(Q, {\bf v_1})$ and $y \in \Hom^1(Q, {\bf v_2}, {\bf v_1})$, where a single arrow is used to denote the representations $x$ and $y$ in order to simplify the above diagram. When $y \in \ker (\rho_x^\mathbb{C})^*$ is nonzero, then these representations appear in the negative slice \eqref{eqn:neg-slice-framed}. If $y \notin \ker (\rho_x^\mathbb{C})^*$, then \cite[Lem. 3.26]{Wilkin17} shows that there exists $g \in G_{\bf v}$ so that $g \cdot y$ so is in the negative slice.

The relation map $\nu : \Rep(Q, {\bf v}) \rightarrow \Rel(Q, {\bf v}, {\bf v}, \mathcal{R})$ decomposes into
\begin{equation*}
\nu_1 \oplus \nu_2 : \Rep(Q, {\bf v_1}) \oplus \Hom^1(Q, {\bf v_2}, {\bf v_1}) \rightarrow \Rel(Q, {\bf v_1}, {\bf v_1}, \mathcal{R}) \oplus \Rel(Q, {\bf v_2}, {\bf v_1}, \mathcal{R}) ,
\end{equation*}
where $\nu_1$ is the restriction of $\nu$ to $\Rep(Q, {\bf v_1})$. To define $\nu_2$, first note that since $y$ is nilpotent, then each path $p = a_n \cdots a_1$ is mapped to the homomorphism
\begin{align*}
\Rep(Q, {\bf v_1}) \oplus \Hom^1(Q, {\bf v_2}, {\bf v_1}) & \rightarrow \Rel(Q, {\bf v_2}, {\bf v_1}, \mathcal{R}) \\
 (x, y) & \mapsto (x,y)_p := x_{a_n} \cdots x_{a_2} y_{a_1} .
\end{align*}
Equivalently, let $p' = a_n \cdots a_2$ be the truncated path obtained by removing $a_1$ from $p$ and note that $(x,y)_p = x_{p'} y_{a_1}$.  Given $(\delta x, \delta y) \in \Hom^1(Q, {\bf v_1}, {\bf v_1}) \oplus \Hom^1(Q, {\bf v_2}, {\bf v_1})$, the derivative \eqref{eqn:path-derivative-def} of each path is
\begin{equation*}
dp_{(x,y)}(\delta x, \delta y) = dp_x'(\delta x) y_{a_1} + x_{p'} (\delta y)_{a_1} .
\end{equation*}

Now define
\begin{equation*}
\nu_2(x,y) = \bigoplus_{r \in \mathcal{R}} \sum_{p \in \mathcal{P}_{t(r), h(r)}} \lambda_p(r) \,  (x,y)_p 
\end{equation*}
and
\begin{align}\label{eqn:negative-relations}
\begin{split}
(d \nu_2)_{(x,y)}(\delta x, \delta y) & = \bigoplus_{r \in \mathcal{R}} \sum_{p \in \mathcal{P}_{t(r), h(r)}} \lambda_p(r) \, dp_{(x,y)}(\delta x, \delta y) \\
 & = \bigoplus_{r \in \mathcal{R}} \sum_{p \in \mathcal{P}_{t(r), h(r)}} \lambda_p(r) \,  dp_x'(\delta x) y_{a_1} + \bigoplus_{r \in \mathcal{R}} \sum_{p \in \mathcal{P}_{t(r), h(r)}} \lambda_p(r) \,  x_{p'} (\delta y)_{a_1} \\
 & =: (d \nu_2)_y(\delta x) + (d \nu_2)_x(\delta y) ,
\end{split}
\end{align}
where the above is used as the definition of $(d \nu_2)_y(\delta x)$ and $(d \nu_2)_x(\delta y)$. Now the derivative of the relation map is
\begin{align}\label{eqn:derivative-neg-slice-quiver}
\begin{split}
(d\nu_1 \oplus d\nu_2)_{(x,y)} : \Hom(Q, {\bf v_1}, {\bf v_1}) \oplus \Hom^1(Q, {\bf v_2}, {\bf v_1}) & \rightarrow \Rel(Q, {\bf v_1}, {\bf v_1}, \mathcal{R}) \oplus \Rel(Q, {\bf v_2}, {\bf v_1}, \mathcal{R}) \\
(\delta x, \delta y) & \mapsto \left( (d\nu_1)_x(\delta x), (d \nu_2)_y(\delta x) + (d \nu_2)_x(\delta y)  \right) .
\end{split}
\end{align}
\end{example}

This particular example and the notation for the derivative of $\nu$ will be used in Section \ref{sec:convolution-cup-product}. We would like to apply the deformation theory of Lemma \ref{lem:smoothness-restricts} to this negative slice quiver and hence the Hecke correspondence, however unfortunately the negative slice quiver is not a full restriction of the original quiver $Q$, since a deformation in $\Hom^1(\tilde{Q}', {\bf v_1 + v_2}, {\bf v_1 + v_2})$ of the form
\begin{equation}\label{eqn:full-deformations}
\xymatrixrowsep{0.5in}
\xymatrixcolsep{0.5in}
\xymatrix{
\bullet_{(Q, {\bf v_1})} \ar@`{(-15,15),(15,15)}^{\delta x_1} \ar@/_0.8pc/[r]_{\delta z}  & \bullet_{(Q, {\bf v_2})} \ar@/_0.8pc/[l]_{\delta y} \ar@`{(10,15),(40,15)}^{\delta x_2}  
}
\end{equation}
may contain relations with paths whose derivative has terms such as $y (\delta z)$ and $y (\delta x_2)$ which both contain exactly one edge in $\tilde{\mathcal{E}}' \setminus \mathcal{E}'$. The next lemma shows that the first type of path is not an obstacle towards proving that $d \nu_{(x+y)}$ is surjective, and the second type of path can be dealt with using the condition of Definition \ref{def:loop-neg-slice}.

\begin{lemma}\label{lem:stable-implies-adjoint-injective}
Let $Q = (\mathcal{E}, \mathcal{I})$ be a framed quiver with relations $\mathcal{R}$ and canonical stability parameter given by Definition \ref{def:aasp}. Suppose that for any dimension vector ${\bf v}$ we have $x$ stable implies $d \nu_x$ is surjective. 

Choose a dimension vector ${\bf v}$ and vertex $k \in \mathcal{I}$, and consider the restricted quiver associated to a reduction of structure group given by the negative slice quiver
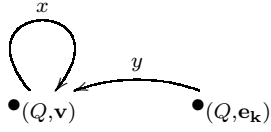
\begin{equation*}
\xymatrixrowsep{0.5in}
\xymatrixcolsep{0.5in}
\xymatrix{
\bullet_{(Q, {\bf v})} \ar@`{(-15,15),(15,15)}^{x}  & \bullet_{(Q, {\bf e_k})} \ar@/_0.8pc/[l]_{y} 
}
\end{equation*}
with $x \in \Rep(Q, {\bf v})^{\alpha-st}$ and suppose that $y \in S_x^-$ is nonzero. Then $d \nu'_{(x+y)}$ is surjective onto $\Rel_0(Q, {\bf v_1}, {\bf v_1}, \mathcal{R}) \oplus \Rel(Q, {\bf v_2}, {\bf v_1}, \mathcal{R})$.
\end{lemma}

\begin{proof}

Since $x$ is stable and $y \in S_x^-$ then the negative gradient flow \eqref{eqn:neg-grad-flow-def} takes $x+y$ down to a lower critical point for $\| \mu - \alpha \|^2$ on $\Rep(Q, {\bf v} + {\bf e_k}, \mathcal{R})$ \cite{Wilkin17}. Since $x$ is at the lowest non-minimal critical level for $\| \mu - \alpha \|^2$ in $\Rep(Q, {\bf v} + {\bf e_k}, \mathcal{R})$, then this lower critical point must be a minimiser for $\| \mu - \alpha \|^2$ and therefore $x+y$ is stable. In particular, this implies that 
\begin{equation}\label{eqn:full-relation-map}
\begin{tikzcd}
\Hom^1(\tilde{Q}', {\bf v + e_k}, {\bf v + e_k}) \arrow{r}{d\nu_{(x+y)}} & \Rel_0(\tilde{Q}', {\bf v + e_k}, {\bf v + e_k}, \tilde{\mathcal{R}}')
\end{tikzcd}
\end{equation}
is surjective for the unrestricted quiver $\tilde{Q}'$. To complete the proof we need to show that this implies that the induced homomorphism 
\begin{equation}\label{eqn:restricted-relation-map}
\begin{tikzcd}
\Hom^1(Q', {\bf v}, {\bf v}) \oplus \Hom^1(Q', {\bf e_k}, {\bf v}) \arrow{r}{d \nu_{(x+y)}'} & \Rel_0(Q', {\bf v}, {\bf v}, \mathcal{R}') \oplus \Rel(Q', {\bf e_k}, {\bf v}, \mathcal{R}')
\end{tikzcd}
\end{equation} 
from \eqref{eqn:relate-restricted-deformations} is surjective.


With reference to the notation for the deformations from \eqref{eqn:full-deformations}, first note that $d \nu_x(\delta x_1)$ defines a surjective map onto the first component $\Rel_0(Q', {\bf v}, {\bf v}, \mathcal{R}')$, since $x$ is stable. Therefore $d \nu'(\delta x_1)$ also defines a surjective map. Equivalently, this shows that the term $d\nu_y(\delta z)$ from the relations on $\tilde{Q}'$ does not affect the surjectivity of $d \nu'$ onto the first component.


The second component of $d \nu'_{(x+y)}$ mapping into $\Rep(Q', {\bf e_k}, {\bf v}, \mathcal{R}')$ is $(d \nu_y)(\delta x_1) + (d \nu_x)(\delta y)$ from \eqref{eqn:negative-relations}. If there are no loops at vertex $k \in \mathcal{I}$ then the deformation $\delta x_2$ in \eqref{eqn:full-deformations} is zero, and so the second component of $d \nu_{(x+y)}'$ is induced from the larger quiver and therefore surjective.


If there is a loop at vertex $k \in \mathcal{I}$ with scalar $\delta x_2 \in \mathbb{C}$, then Definition \ref{def:loop-neg-slice} implies that any relations of the form $y (\delta x_2) \in \Rel(\tilde{Q}', {\bf e_k}, {\bf v_1}, \tilde{\mathcal{R}}')$ can be cancelled by choosing $\delta x_1' = \frac{\lambda_p(r)}{\lambda_{p'}(r)} \delta x_2 \cdot \id$ in the corresponding relation $(\delta x_1') y$ appearing in $\Rel(Q', {\bf e_k}, {\bf v_1}, \mathcal{R}')$, which exists and is unique by the condition in Definition \ref{def:loop-neg-slice}. This can be done for each loop, and therefore surjectivity of $d \nu_{(x+y)}$ in \eqref{eqn:full-relation-map} implies surjectivity of $d \nu_{(x+y)}'$ in \eqref{eqn:restricted-relation-map}, which completes the proof.
\end{proof}

\begin{corollary}\label{cor:neg-slice-quiver-smooth}
With the same conditions as Lemma \ref{lem:stable-implies-adjoint-injective}, the moduli space associated to the negative slice quiver is smooth.
\end{corollary}

\begin{proof}
Lemma \ref{lem:stable-implies-adjoint-injective} shows that $d \nu'_{(x+y)}$ is surjective, after which smoothness follows from the standard deformation theory of the moduli space using the same method as Corollary \ref{cor:smoothness-restricts}.
\end{proof}



\section{The Morse spectral sequence}\label{sec:Morse-spectral}

In this section we recall some basic facts about the differentials and cup product on the spectral sequence associated to a Morse filtration of a topological space.

\subsection{Morse stratification}\label{subsec:Morse-stratification}

Given a topological space $M$ and a filtration $\emptyset = M_{-1} \subset M_0 \subset \cdots \subset M_n = M$, one can write down a spectral sequence where the terms on the first page are
\begin{equation*}
E_1^{p,q} = H^p(M_{q-p}, M_{q-p-1}) .
\end{equation*}
The goal is then to compute the terms $E_\infty^{p,q}$ which correspond to the graded cohomology of $M$. When $M$ is a manifold, then the Morse-theoretic point of view is to use a smooth function $f : M \rightarrow \mathbb{R}$ satisfying certain conditions (which we make precise below) to derive a \emph{Morse filtration} of $M$ and carry out the spectral sequence calculations by analysing the behaviour of $f$ near the critical points.

Suppose that for all $x \in M$ the negative gradient flow of $f$ with initial condition $x$, denoted $\varphi_t(x)$,  is defined for all $t \in [0, \infty)$ and that $\lim_{t \rightarrow \infty} \varphi_t(x)$ exists. For example, the Lojasiewicz inequality method \cite{Lojasiewicz84}, \cite{Simon83} shows that this is true when the function is analytic and the flow satisfies a compactness condition, which is true for the norm-square of a moment map on an affine variety (cf. \cite[Lem. 4.10]{Sjamaar98}). In this case one can construct a filtration $\emptyset = M_{-1} \subset M_0 \subset \cdots \subset M_n = M$ determined by connected components of the set of critical points of $f$.

If the function also has good properties near the critical sets (e.g. it is Morse, Morse-Bott \cite{Bott54} or minimally degenerate \cite{Kirwan84}) then for each critical set $C_j$ labelled by $j = 1,  \ldots, n$, there is a negative normal bundle $V_j \rightarrow C_j$ whose fibres correspond to the negative eigenspace of the Hessian at each critical point. Let $V_{j,0} = V_j \setminus C_j$ and let $\lambda_j := \rank_\mathbb{R} V_j$ be the index. The \emph{main theorem of Morse theory} shows that the terms on the first page of the spectral sequence can be written in terms of relative cohomology groups localised around the critical sets
\begin{equation*}
H^*(M_j, M_{j-1}) \cong H^*(V_j, V_{j,0}) ,
\end{equation*}
and the Thom isomorphism shows that these relative cohomology groups can be expressed in terms of Morse data (the critical set together with the Morse index)
\begin{equation*}
H^*(V_j, V_{j,0}) \cong H^{*-\lambda_j}(C_j) .
\end{equation*}

One then has to compute the differentials in the spectral sequence. If the gradient flow of the function has further good properties in the sense that the stable and unstable manifolds intersect transversely (e.g. the function is Morse-Smale or Morse-Bott-Smale) then one can express these differentials in terms of pullback/pushforward homomorphisms on spaces of flow lines, which is worked out in detail in \cite{AustinBraam95} using de Rham cohomology.

Therefore, if the function is sufficiently well-behaved, then Morse theory determines 
\begin{enumerate}

\item a canonical Morse filtration (which requires good compactness properties of the flow),

\item  the terms on the first page of the spectral sequence in terms of critical point data (which requires good behaviour of the function near the critical sets), and

\item the differentials of the spectral sequence in terms of flow line data (which requires a transversality condition).

\end{enumerate}

More generally, one can replace the manifold $M$ by a (possibly singular) subset $Z$ preserved by the flow $\varphi$ in the sense that if the initial condition $x_0 \in Z$, then $\varphi_t(x_0) \in Z$ for all $t$ such that $\varphi_t(x_0)$ is defined. Under some additional conditions, which are explained in detail in \cite{Wilkin19} and which we recall below, there is still a Morse filtration $\emptyset = Z_{-1} \subset Z_0 \subset \cdots \subset Z_n = Z$ and the terms on the first page of the spectral sequence can still be written in terms of relative cohomology groups localised around the critical sets.


Let $M$ be real analytic, let $f : M \rightarrow \mathbb{R}$ and let $\varphi_t(x)$ denote the time $t$ downwards gradient flow of $f$ with initial condition $x \in M$. Let $Z$ be a subset preserved by the gradient flow of $f$ in the sense that $x \in Z$ implies that $\varphi_t(x) \in Z$ for all $t$ such that $\varphi_t(x)$ exists. We can then define a critical point to be a fixed point of the flow. As above, label the critical sets by $C_j$ for $j = 0, \ldots, n$ and define
\begin{equation*}
W_j^- := \left\{ x \in Z \, \mid \, \lim_{t \rightarrow -\infty} \varphi_t(x) \in C_j \right\}, \quad W_{j,0}^- := W_{j}^- \setminus C_j .
\end{equation*}

The main result of \cite{Wilkin19} is that the following analog of the main theorem of Morse theory holds for $f : Z \rightarrow \mathbb{R}$ if certain conditions are satisfied.

\begin{theorem}[{\cite[Thm. 1.1 \& 1.3]{Wilkin19}}]\label{thm:main-thm-morse}
If the pair $(Z \subset M, f : M \rightarrow \mathbb{R})$ satisfies Conditions (1)--(5) of \cite{Wilkin19}, then there is a Morse filtration $\emptyset = Z_{-1} \subset Z_0 \subset \cdots \subset Z_n = Z$ and for any real numbers $a < b$ the following is true.
\begin{itemize}
\item If there are no critical values in $[a,b]$ then $Z_b \simeq Z_a$.

\item If $a$ and $b$ are not critical values and there is one critical value $c \in (a,b)$ with associated critical set $C_j$, then $Z_{j} \simeq Z_{j-1} \cup W_j^-$.   
\end{itemize}
Moreover, if $K$ is a compact Lie group acting on $Z$ and $f$ is $K$-invariant, then the above homotopy equivalences can be chosen to be $K$-equivariant.
\end{theorem}

Applying excision gives us the following isomorphism in $K$-equivariant cohomology
\begin{equation}\label{eqn:main-theorem-filtration}
H_K^*(Z_j, Z_{j-1}) \cong H_K^*(W_j^-, W_{j,0}^-) .
\end{equation}

Therefore the terms on the first page of the spectral sequence for the filtration $\emptyset = Z_{-1} \subset Z_0 \subset \cdots \subset Z_n = Z$ can be expressed in terms of relative cohomology groups defined in terms of local data around the critical sets.

The results of \cite{Wilkin19} also show that the above theorem is valid for the norm square of a moment map on an affine variety (see also \cite{wilkin-condition-4} for a more streamlined proof). In particular, we can apply this theorem in the general setting of representations of a quiver with a finite set of relations.

\subsection{A Thom homomorphism}\label{subsec:Thom-hom}



Now consider the case where the function $f:M \rightarrow \mathbb{R}$ on the ambient manifold $M$ is minimally degenerate in the sense of Kirwan \cite{Kirwan84} and Conditions (1)--(5) of \cite{Wilkin19} are satisfied for the restriction of $f$ to $Z \subset M$. The model example is where $M$ is the space of representations of a quiver, $f : M \rightarrow \mathbb{R}$ is the norm-square of a moment map as in Section \ref{subsec:moment-map-def} and $Z$ is the subvariety satisfying a finite set of relations on the quiver.

Since $f$ is minimally degenerate on $M$, then the unstable set $V_j$ of the critical set $C_j$ is a vector bundle. Let $\lambda_j : = \rank_\mathbb{R} V_j$. Then the Thom isomorphism shows that the terms on the first page of the spectral sequence for $M$ can be written as $H_K^p(M_j, M_{j-1}) \cong H_K^{p-\lambda_j}(C_j)$ and these fit into the following commutative diagram
\begin{equation}\label{eqn:smooth-spectral-sequence-Thom}
\begin{tikzcd}
H_K^p(M_j, M_{j-1}) \arrow{r} \arrow[leftrightarrow]{d}{\text{Thm.} \, \ref{thm:main-thm-morse}}[swap]{\cong} & H_K^p(M_j) \arrow{d}{restriction} \\
H_K^p(V_j, V_{j,0}) \arrow{r} \arrow[leftrightarrow]{d}{Thom}[swap]{\cong} & H_K(V_j) \arrow{d}{restriction} \\
H_K^{p-\lambda_j}(C_j) \arrow{r}{\smallsmile e} & H_K^p(C_j) 
\end{tikzcd}
\end{equation}
where the isomorphism $H_K^{p-\lambda_j}(C_j) \cong H_K^{p-\lambda_j}(V_j) \rightarrow H_K^p(V_j, V_{j,0})$ is given by pullback followed by cup product with the Thom class $\tau \in H_K^{\lambda_j}(V_j, V_{j,0})$.

Let $i : C_j \cap Z \hookrightarrow C_j$ denote the inclusion. In the following we will use $W_j$ to denote the unstable set of $C_j \cap Z$ of the restriction of $f$ to $Z$ and $V_j^Z := i^* V_j$ to denote the pullback of the negative normal bundle from the ambient manifold $M$. The Euler class of the bundle $V_j^Z \rightarrow C_j \cap Z$ is denoted $e_Z = i^* e$. Note that $W_j \subset V_j^Z$.

When we restrict to the subspace $Z$, the lower half of the above diagram becomes
\begin{equation*}
\begin{tikzcd}
H_K^p(W_j, W_{j,0}) \arrow{r} & H_K^p(W_j) \\
H_K^p(V_j^Z, V_{j,0}^Z) \arrow{u}{restriction} \arrow{r} & H_K^p(V_j^Z) \arrow[leftrightarrow]{u}{\cong} \\
H_K^{p-\lambda_j}(C_j \cap Z) \arrow[leftrightarrow]{u}{Thom}[swap]{\cong} \arrow{r}{\smallsmile e_Z} & H_K^p(C_j \cap Z) \arrow[leftrightarrow]{u}{\cong}
\end{tikzcd}
\end{equation*}
Therefore, we see that cup product with the Euler class $H_K^{p-\lambda_j}(C_j \cap Z) \stackrel{\smallsmile e_Z}{\xrightarrow{\hspace{0.8cm}}} H_K^p(C_j \cap Z)$ factors through the homomorphism $H_K^{p-\lambda_j}(C_j \cap Z) \rightarrow H_K^p(W_j, W_{j,0})$ in the sense that the following diagram commutes. 
\begin{equation*}
\begin{tikzcd}
H_K^p(W_j, W_{j,0}) \arrow{r} & H_K^p(W_j) \\
H_K^{p-\lambda_j}(C_j \cap Z) \arrow{r}{\smallsmile e_Z} \arrow{u} & H_K^p(C_j \cap Z) \arrow[leftrightarrow]{u}{\cong}
\end{tikzcd}
\end{equation*}
In particular, if the map $H_K^{p-\lambda_j}(C_j \cap Z) \overset{\smallsmile e_Z}{\xrightarrow{\hspace{0.8cm}}} H_K^p(C_j \cap Z)$ is injective, then so is the homomorphism $H_K^{p-\lambda_j}(C_j \cap Z) \rightarrow H_K^p(W_j, W_{j,0})$.

The above construction is summarised in the following result.

\begin{lemma}\label{lem:Thom-hom}
Let $\pi : V \rightarrow B$ be a vector bundle with $\lambda = \rank_\mathbb{R} V$, Euler class $e \in H_K^\lambda(B)$ and Thom class $\tau \in H_K^p(V, V \setminus B)$. Consider a subspace $i : W \hookrightarrow V$ such that $K \cdot i(W) = i(W)$ and the zero section $B \hookrightarrow V$ lies in $i(W)$. Define $V_0 := V \setminus B$ and $W_0 := W \setminus B$. Then there is a homomorphism
\begin{equation}\label{eqn:Thom-homomorphism}
H_K^{p - \lambda}(B) \xrightarrow{i^* \pi^*(\cdot) \smallsmile i^* \tau} H_K^p(W, W_0) \\
\end{equation}
Moreover, if $H_K^{p-\lambda}(B) \stackrel{\smallsmile e}{\longrightarrow} H_K^p(B)$ is injective, then so is the homomorphism \eqref{eqn:Thom-homomorphism}.
\end{lemma}

Therefore restricting the diagram \eqref{eqn:smooth-spectral-sequence-Thom} to the subspace $Z$ gives us the following diagram
\begin{equation*}
\begin{tikzcd}
H_K^p(Z_j, Z_{j-1}) \arrow{r} \arrow[leftrightarrow]{d}{\text{Thm.} \, \ref{thm:main-thm-morse}}[swap]{\cong} & H_K^p(Z_j) \arrow{d}{restriction} \\
H_K^p(W_j, W_{j,0}) \arrow{r}& H_K(W_j) \arrow[leftrightarrow]{d}{restriction}[swap]{\cong} \\
H_K^{p-\lambda_j}(C_j \cap Z) \arrow{r}{\smallsmile e_Z}  \arrow{u}[swap]{\text{Lem.} \, \ref{lem:Thom-hom}} & H_K^p(C_j \cap Z) 
\end{tikzcd}
\end{equation*}
from which we see that even though the terms on the first page of the spectral sequence for the filtration $\emptyset \subset Z_0 \subset \cdots \subset Z_n = Z$ may not be isomorphic to the terms $H_K^{p-\lambda_j}(C_j \cap Z)$, we still have a homomorphism $H_K^{p-\lambda_j}(C_j \cap Z) \rightarrow H_K^p(Z_j, Z_{j-1})$. If $H_K^{p-\lambda_j}(C_j \cap Z) \stackrel{\smallsmile e_Z}{\longrightarrow} H_K^p(C_j \cap Z)$ is injective, then so is $H_K^{p-\lambda_j}(C_j \cap Z) \rightarrow H_K^p(Z_j, Z_{j-1})$.

Finally, we have the following lemma which will be useful in Section \ref{sec:transversality}.

\begin{lemma}\label{lem:Euler-subbundle}
Let $E \cong E_1 \oplus E_2 \rightarrow B$ be a direct sum of vector bundles, and let 
\begin{equation*}
\nu_1 = \rank E_1, \quad \nu_2 = \rank E_2 , \quad \lambda = \nu_1 + \nu_2 = \rank E .
\end{equation*} 
Then the following diagram commutes, where the homomorphism in the top row is induced by inclusion, the left hand column is the Thom isomorphism for the bundle $E \rightarrow B$, the middle column is the Thom isomorphism for the bundle $E \rightarrow E_1$, the right hand column is the Thom isomorphism for $E_2$ and the bottom row is cup product with the Euler class of $E_1$.
\begin{equation*}
\begin{tikzcd}[column sep=1.5cm]
H^p(E, E_0) \arrow{r} \arrow[leftrightarrow]{dd}{\cong}[swap]{Thom} & H^p(E, E \setminus E_1) \arrow[leftrightarrow]{d}{\cong}[swap]{Thom} \arrow{r}{\cong}[swap]{homotopy} & H^p(E_2, E_{2,0})  \arrow[leftrightarrow]{dd}{\cong}[swap]{Thom} \\
 & H^{p-\nu_2}(E_1) \arrow[leftrightarrow]{d}{\cong}[swap]{homotopy} \\
H^{p-\lambda}(B) \arrow{r}{\smallsmile e_1} & H^{p-\nu_2}(B) \arrow{r}{=} & H^{p-\nu_2}(B)
\end{tikzcd}
\end{equation*}
\end{lemma}

\begin{proof}
The proof follows from the Whitney sum formula for the Euler class $e = e_1 \smallsmile e_2$ together with the homotopy equivalences in the above diagram.
\end{proof}

\subsection{A splitting lemma}

Consider the triple $(Z_{j+1}, Z_j, Z_{j-1})$, let $C_{j+1}$ and $C_j$ denote the associated critical sets in $Z_{j+1} \setminus Z_j$ and $Z_j \setminus Z_{j-1}$, and let $\mathcal{F}_{j,0}^{j+1,0}$ denote the space of flow lines between the lower critical set $C_j$ and the upper critical set $C_{j+1}$. The construction of the cup product in the sequel depends on a choice of splittings for the long exact sequences of different triples $(Z_{j+1}, Z_j, Z_{j-1})$, $(W_{j+1}, W_{j+1,0}, W_{j+1,0} \setminus \mathcal{F}_{j,0}^{j+1,0})$ and $(W_{j+1}, W_{j+1,0}, \emptyset)$. In this section we prove the following lemma which shows that these choices can be made to be compatible, which will be used to show that the different splittings used to define the cup product in Proposition \ref{prop:cup-product-main-theorem} can be chosen to be compatible so that the diagram \eqref{eqn:restriction-cup-product-relative} commutes. The lemma and its proof are elementary, however we were unable to find the precise statement we need in the literature, therefore we include all the details here in order to make the exposition self-contained. 

\begin{lemma}\label{lem:splitting-compatible}
Consider the following commutative diagram of short exact sequences of abelian groups
\begin{equation*}
\begin{tikzcd}
0 \arrow{r} & A \arrow{r}{j} \arrow{d}{\cong} & B \arrow{r}{k} \arrow{d}{r} & C \arrow{r} \arrow{d}{t} & 0 \\
0 \arrow{r} & A' \arrow{r}{j'} & B' \arrow{r}{k'} & C' \arrow{r} & 0 
\end{tikzcd}
\end{equation*}
and suppose that all groups are divisible abelian groups so that all short exact sequences split. Then given a choice of splitting homomorphism $i' : C' \rightarrow B'$ there exists a splitting homomorphism $i : C \rightarrow B$ such that $r \circ i = i' \circ t$.
\end{lemma}

\begin{proof}
Choose a splitting homomorphism $s : \im r \rightarrow B$ so that $r \circ s = \id$ and note that $r \circ j = j'$ implies that $\im j' \subset \im r$, therefore $s$ is defined on $\im j'$.

First we claim that $i' \circ t(c) \in \im r$ for all $c \in C$. To see this, note that
\begin{align*}
k' \circ i' \circ t(c) = t(c) = t \circ k(b) & = k' \circ r(b) \quad \text{for some $b \in B$ defined up to $\im j$} \\
 \Rightarrow \quad i' \circ t(c) - r(b) & \in \ker k' = \im j' = \im r \circ j \\
 \Rightarrow \quad i' \circ t(c) & = r(b+j(a)) \quad \text{for some $a \in A$} .
\end{align*}
Therefore $s \circ i' \circ t$ is well-defined and we can construct a splitting homomorphism $i := s \circ i' \circ t : C / \ker t \rightarrow B$ such that $r \circ i = i' \circ t$. It now only remains to define $i$ on $\ker t$.

Choose any splitting homomorphism $i'' : C \rightarrow B$ so that $k \circ i'' = \id$. If $c \in \ker t$, then we have
\begin{align*}
t \circ k \circ i''(c) & = 0 \\
\Rightarrow \quad k' \circ r \circ i''(c) & = 0 \\
\Rightarrow \quad r \circ i''(c) & \in \ker k' = \im j' = \im r \circ j .
\end{align*}
Therefore there exists $a \in A$ such that $r(i''(c) + j(a)) = 0$ and the choice of $a$ is unique since $r \circ j = j'$ is injective. Now define $i : \ker t \rightarrow \ker r$ by $i(c) := i''(c) + j(a) \in \ker r$, which trivially satisfies $r \circ i = i' \circ t$. It is easy to check that $i$ is a homomorphism. Together with $i : C / \ker t \rightarrow B$ defined above, since $C$ is a divisible abelian group then this determines $i : C \rightarrow B$ satisfying the conditions of the lemma.
\end{proof}

\subsection{Differentials for the spectral sequence}

Each differential on the first page of the spectral sequence is the composition of the following two homomorphisms
\begin{equation}\label{eqn:SS-differential}
\begin{tikzcd}
H^p(Z_j, Z_{j-1}) \arrow{r} \arrow[bend left=15]{rr}{d}& H^p(Z_j) \arrow{r} & H^{p+1}(Z_{j+1},Z_j) ,
\end{tikzcd}
\end{equation}
where the first homomorphism comes from the long exact sequence of the pair $(Z_j, Z_{j-1})$ and the second homomorphism is the connecting homomorphism from the long exact sequence of the pair $(Z_{j+1}, Z_j)$.

The construction of Section \ref{subsec:Thom-hom} shows that there are homomorphisms from the shifted cohomology of the critical sets that fit into the following diagram
\begin{equation*}
\begin{tikzcd}
H^p(Z_j, Z_{j-1}) \arrow{r} \arrow[bend left=15]{rr}{d} & H^p(Z_j) \arrow{r} & H^{p+1}(Z_{j+1},Z_j) \\
H^p(W_j, W_{j,0}) \arrow[leftrightarrow]{u}{\text{Thm.} \, \ref{thm:main-thm-morse}}[swap]{\cong} & & H^{p+1}(W_{j+1}, W_{j+1,0}) \arrow[leftrightarrow]{u}{\text{Thm.} \, \ref{thm:main-thm-morse}}[swap]{\cong} \\
H^{p-\lambda_j}(C_j) \arrow{u}{\text{Lem.} \, \ref{lem:Thom-hom}} & & H^{p-\lambda_{j+1} + 1}(C_{j+1}) \arrow{u}{\text{Lem.} \, \ref{lem:Thom-hom}}
\end{tikzcd}
\end{equation*}
It is then natural to ask the following questions.

\renewcommand\theenumi{\roman{enumi}}
\begin{enumerate}

\item Does the differential $d : H^p(Z_j, Z_{j-1}) \rightarrow H^{p+1}(Z_{j+1},Z_j)$ map the image of the left-hand column of the above diagram $H^{p-\lambda_j}(C_j) \rightarrow H^p(Z_j, Z_{j-1})$ to the image of the right-hand column $H^{p-\lambda_{j+1} + 1}(C_{j+1}) \rightarrow H^{p+1}(Z_{j+1},Z_j)$? 

\item If so, is it possible to construct an induced map $H^{p-\lambda_j}(C_j) \rightarrow H^{p-\lambda_{j+1} + 1}(C_{j+1})$?

\item If the answers to the above questions are positive, then does the induced map $H^{p-\lambda_j}(C_j) \rightarrow H^{p-\lambda_{j+1} + 1}(C_{j+1})$ have any topological meaning?

\end{enumerate}

For Morse-Bott-Smale functions on smooth spaces, the differential can be written in terms of pullback/pushforward homomorphisms involving spaces of flow lines (see for example \cite{AustinBraam95}). The main result of Section \ref{subsec:precise-differential} is that one can do this for spaces of representations of quivers and that the topological meaning of the induced map is that it can be expressed in terms of pullback/pushforward homomorphisms via spaces of flow lines, or equivalently (using the results of \cite{Wilkin17}) via the Hecke correspondence.

\subsection{The cup product on the spectral sequence for adjacent strata}\label{subsec:stratified-cup-product}

In this section we describe the cup product on the $E_1$ page of the spectral sequence induced from the cup product on the total space of $Z$ in the case where the differentials are all zero. This is purely algebraic and designed to develop the basis for subsequent sections in which we will use the analytic results of \cite{Wilkin17} and Section \ref{subsec:flow-lines} to show that 
\begin{enumerate}

\item Theorem \ref{thm:main-thm-morse} can be used to reduce the cup product homomorphism to cohomology groups localised around critical sets and spaces of flow lines (Section \ref{subsec:main-thm-cup-product}), and 

\item a modified transversality condition can be used to show that the cup product homomorphism can be constructed from pullback and pushforward homomorphisms between critical sets and spaces of flow lines connecting them (Section \ref{subsec:precise-cup-product}).

\end{enumerate}

\begin{remark}

For now we focus on the case of a cup product homomorphism between adjacent strata. It is easy to generalise the results of this section to non-adjacent strata, but to generalise the results of the later sections requires compactifying spaces of flow lines, which we defer to a subsequent paper.

\begin{itemize}

\item To simplify the exposition, \emph{for the remainder of this section we will assume that the filtration is perfect}, i.e. the long exact sequence for each pair $(Z_j, Z_{j-1})$ splits into short exact sequences, or equivalently the differentials in the spectral sequence are all zero. 

\item For notation, cohomology groups are written $H^*$ with coefficients in a divisible abelian group so that short exact sequences split, however all of the constructions work in the same way for equivariant cohomology $H_K^*$ with respect  to a group $K$ preserving each stratum $S_j = Z_j \setminus Z_{j-1}$.

\end{itemize}
\end{remark}

Given a triple $(Z, Z_j, Z_{j-1})$, perfection of the stratification shows that the following sequence is exact
\begin{equation*}
\begin{tikzcd}
0 \arrow{r} & H^*(Z, Z_j) \arrow{r} & H^*(Z, Z_{j-1}) \arrow{r} & H^*(Z_j, Z_{j-1}) \arrow{r} & 0 .
\end{tikzcd}
\end{equation*}
Given $\omega \in H^m(Z)$, there is an induced cup product homomorphism on the terms in the above exact sequence, which can be represented as follows.
\begin{equation}\label{eqn:cup-product-exact-strata}
\begin{tikzcd}
0 \arrow{r} & H^p(Z, Z_j) \arrow{r}{j} \arrow{d}{\smallsmile \omega} & H^p(Z, Z_{j-1}) \arrow{r}{k} \arrow{d}{\smallsmile \omega} & H^p(Z_j, Z_{j-1}) \arrow{r} \arrow{d}{\smallsmile \omega_j} & 0 \\
0 \arrow{r} & H^{p+m}(Z, Z_j) \arrow{r}{j} & H^{p+m}(Z, Z_{j-1}) \arrow{r}{k} & H^{p+m}(Z_j, Z_{j-1}) \arrow{r} & 0 
\end{tikzcd}
\end{equation}
where $\omega_j \in H^m(Z_j)$ denotes the restriction of $\omega$ to $Z_j$. 

Using the assumption on the coefficients, the above exact sequences split, and so the cup product induces a homomorphism
\begin{equation}\label{eqn:splitting-consequence}
\begin{tikzcd}
H^p(Z, Z_j) \oplus H^p(Z_j, Z_{j-1}) \arrow{r}{\smallsmile \omega} & H^{p+m}(Z, Z_j) \oplus H^{p+m}(Z_j, Z_{j-1}) .
\end{tikzcd}
\end{equation}
For the Morse-theoretic convolution of Section \ref{sec:morse-convolution}, we will focus attention on the off-diagonal component of the above homomorphism, denoted by
\begin{equation}\label{eqn:induced-upper-cup-product}
\begin{tikzcd}
H^p(Z_j, Z_{j-1}) \arrow{r}{c_\omega} & H^{p+m}(Z, Z_j)
\end{tikzcd}
\end{equation}
which is constructed as follows. First note that in writing \eqref{eqn:splitting-consequence} we have implicitly chosen a splitting of the exact sequences \eqref{eqn:cup-product-exact-strata}, in which there is an inclusion $i : H^p(Z_j, Z_{j-1}) \rightarrow H^p(Z, Z_{j-1})$ such that $k \circ i = \id$ and a projection $\pi : H^{p+m}(Z, Z_{j-1}) \rightarrow H^{p+m}(Z, Z_j)$ such that $\pi \circ j = \id$. 

Let $\eta \in H^p(Z_j, Z_{j-1})$ and define $\zeta = i(\eta)$. Then $\zeta \smallsmile \omega \in H^{p+m}(Z, Z_{j-1}) \cong H^{p+m}(Z, Z_j) \oplus H^{p+m}(Z_j, Z_{j-1})$ from which we then project to $\pi(\zeta \smallsmile \omega) \in H^{p+m}(Z, Z_j)$.

\begin{lemma}\label{lem:cup-product-strata}
With respect to the splitting homomorphisms $i, \pi$ defined above, the cup product homomorphism \eqref{eqn:induced-upper-cup-product} which is induced from cup product on the total space $\smallsmile \omega : H^p(Z) \rightarrow H^{p+m}(Z)$ is given by
\begin{equation*}
c_\omega(\eta) = \pi(i(\eta) \smallsmile \omega) .
\end{equation*}
\end{lemma}

\begin{remark}
If $\eta \smallsmile \omega_j = 0 \in H^{p+m}(Z_j, Z_{j-1})$ then this is the cup product
\begin{equation*}
H^p(Z_j, Z_{j-1}) \hookrightarrow H^p(Z, Z_{j-1}) \stackrel{\smallsmile \omega}{\longrightarrow} H^{p+m}(Z, Z_j) .
\end{equation*}
\end{remark}

\section{Using the main theorem to localise the differentials and cup product}\label{sec:main-theorem-localisation}

The construction of the differentials and cup product in the previous section is valid for a filtration $\emptyset = Z_{-1} \subset Z_0 \subset \cdots \subset Z_n = Z$ of a topological space $Z$. When the filtration is induced by a Morse filtration on an ambient manifold $M$ (cf. the assumptions of Theorem \ref{thm:main-thm-morse}), then one can interpret these homomorphisms in terms of local data around the critical sets and spaces of flow lines. The conditions of \cite{Wilkin19} are sufficient to define a Morse filtration and prove the isomorphism from Theorem \ref{thm:main-thm-morse}. The goal of this section is to show that these conditions are also sufficient to prove Propositions \ref{prop:differential-main-theorem} and \ref{prop:cup-product-main-theorem} describing the differentials and cup product in terms of local data around the critical sets.

\subsection{Localising the differential around the critical sets and spaces of flow lines}\label{subsec:main-thm-differential}

Let $c_j, c_{j+1}$ denote the critical values associated to the critical sets $C_j, C_{j+1}$. For each stratum $Z_{j+1}$, consider the subspace $N_j := Z_j \cup W_{j+1}$. Note that the flow induces homotopy equivalences $N_j \simeq f^{-1}(-\infty, c_{j+1}] \simeq Z_{j+1}$ and $N_j \setminus C_{j+1} \simeq f^{-1}(-\infty, c_j] \simeq Z_j$ (this requires Conditions (1)--(3) of \cite{Wilkin19}; see \cite[Prop. 2.4]{Wilkin19}). Applying these homotopy equivalences shows that the differential (corresponding to the top row in the following diagram) can be constructed as the composition of the homomorphisms in the second row of the following diagram. The vertical homomorphisms between the second and third rows are induced by restriction to $W_{j+1}$, and Theorem \ref{thm:main-thm-morse} then shows that this induces an isomorphism $H^{p+1}(N_j, N_j \setminus C_{j+1}) \cong H^{p+1}(W_{j+1}, W_{j+1,0})$ in the right-hand column. Note that $\mathcal{F}_j^{j+1}$ is used to denote the union of the space of points that flow up to $C_{j+1}$ and down to $C_j$ with the upper and lower critical sets, while $\mathcal{F}_{j,0}^{j+1,0}$ is the subset that excludes the upper and lower critical sets.

\begin{equation*}
\begin{tikzcd}
H^p(Z_j, Z_{j-1}) \arrow{r} \arrow{d}{\cong} & H^p(Z_j) \arrow{r} \arrow{d}{\cong} & H^{p+1}(Z_{j+1}, Z_j) \arrow{d}{\cong} \\
H^p(N_j \setminus C_{j+1}, N_j \setminus \mathcal{F}_j^{j+1}) \arrow{r} \arrow{d} & H^p(N_j \setminus C_{j+1}) \arrow{r} \arrow{d} & H^{p+1}(N_j, N_j \setminus C_{j+1}) \arrow{d}{\cong} \\
H^p(W_{j+1,0}, W_{j+1,0} \setminus \mathcal{F}_{j,0}^{j+1,0}) \arrow{r} & H^p(W_{j+1,0}) \arrow{r} & H^{p+1}(W_{j+1}, W_{j+1,0})
\end{tikzcd}
\end{equation*}

Therefore (up to the isomorphisms in the right-hand column of the above diagram), it follows from Theorem \ref{thm:main-thm-morse} that the differential in the spectral sequence factors through the homomorphisms in the bottom row, which can be expressed in terms of cohomology groups localised around the critical sets and spaces of flow lines.

\begin{proposition}\label{prop:differential-main-theorem}
Let $Z$ and $f : Z \rightarrow \mathbb{R}$ satisfy Conditions (1)--(5) of \cite{Wilkin19}, so that there is a Morse filtration $\emptyset = Z_{-1} \subset Z_0 \subset \cdots \subset Z_n = Z$ and Theorem \ref{thm:main-thm-morse} applies. Then the differential $H^p(Z_j, Z_{j-1}) \stackrel{d}{\longrightarrow} H^{p+1}(Z_{j+1}, Z_j)$ factors through the composition of the homomorphisms in the bottom row of the following diagram.
\begin{equation}\label{eqn:differential-main-theorem}
\begin{tikzcd}
H^p(Z_j, Z_{j-1}) \arrow{d}{restriction}  \arrow{rr}{d} & &  H^{p+1}(Z_{j+1}, Z_j) \arrow{d}{\text{Thm.} \, \ref{thm:main-thm-morse}}[swap]{\cong} \\
H^p(W_{j+1,0}, W_{j+1,0} \setminus \mathcal{F}_{j,0}^{j+1,0}) \arrow{r} & H^p(W_{j+1,0}) \arrow{r} & H^{p+1}(W_{j+1}, W_{j+1,0}) 
\end{tikzcd}
\end{equation}
\end{proposition}

In Section \ref{subsec:precise-differential} we will show that, in the presence of additional transversality assumptions, the above construction can be further refined using spaces of flow lines in analogy with the well-known theory for Morse-Bott-Smale functions on smooth manifolds (see for example \cite[Sec. 3]{AustinBraam95}).

\subsection{Localising the cup product around the critical sets and spaces of flow lines}\label{subsec:main-thm-cup-product}

Using a similar technique to the proof of Proposition \ref{prop:differential-main-theorem} for the differential, we can use the main theorem of Morse theory to express the cup product on the spectral sequence in terms of cohomology groups localised around the critical sets and spaces of flow lines. The terms in the diagram \eqref{eqn:cup-product-exact-strata} are isomorphic to the corresponding terms in the following diagram.
\begin{equation*}
\begin{tikzcd}[column sep=0.3cm]
0 \arrow{r} & H^{p}(N_j, N_j \setminus C_{j+1}) \arrow{r} \arrow{d}{\smallsmile \omega_{j+1}} & H^p(N_j, N_j \setminus \mathcal{F}_j^{j+1}) \arrow{r} \arrow{d}{\smallsmile \omega_{j+1}} & H^p(N_j \setminus C_{j+1}, N_j \setminus \mathcal{F}_j^{j+1}) \arrow{r} \arrow{d}{\smallsmile \omega_j} & 0 \\
0 \arrow{r} & H^{p+m}(N_j, N_j \setminus C_{j+1}) \arrow{r} & H^{p+m}(N_j \setminus C_{j+1}, N_j \setminus \mathcal{F}_j^{j+1}) \arrow{r} & H^{p+m}(N_j \setminus C_{j+1}, N_j \setminus \mathcal{F}_j^{j+1}) \arrow{r} & 0
\end{tikzcd}
\end{equation*}
In the above diagram and in the sequel we will abuse the notation and also use $\omega_j, \omega_{j+1}$ to denote the restriction of $\omega$ to subspaces of $Z_j$, $Z_{j+1}$ respectively. Restricting all of the above spaces to the intersection with $W_{j+1}$ determines homomorphisms from the terms in the above diagram to the corresponding terms in the following diagram, which are defined using spaces localised around the critical sets and spaces of flow lines.

{\small \begin{equation}\label{eqn:cup-product-exact-local}
\begin{tikzcd}[column sep=0.4cm]
0 \arrow{r} & H^{p}(W_{j+1}, W_{j+1,0}) \arrow{r}{j} \arrow{d}{\smallsmile \omega_{j+1}} & H^p(W_{j+1}, W_{j+1,0} \setminus \mathcal{F}_{j,0}^{j+1,0}) \arrow{r}{k} \arrow{d}{\smallsmile \omega_{j+1}} & H^p(W_{j+1,0}, W_{j+1,0} \setminus \mathcal{F}_{j,0}^{j+1,0}) \arrow{r} \arrow{d}{\smallsmile \omega_j} & 0 \\
0 \arrow{r} & H^{p+m}(W_{j+1}, W_{j+1,0}) \arrow{r}{j} & H^{p+m}(W_{j+1}, W_{j+1,0} \setminus \mathcal{F}_{j,0}^{j+1,0}) \arrow{r}{k} & H^{p+m}(W_{j+1,0}, W_{j+1,0} \setminus \mathcal{F}_{j,0}^{j+1,0}) \arrow{r} & 0
\end{tikzcd}
\end{equation}
}

The horizontal rows in the diagram above are the long exact sequences of the cohomology of the triple $(W_{j+1}, W_{j+1, 0}, W_{j+1,0} \setminus \mathcal{F}_{j,0}^{j+1,0})$. Using the assumption on the coefficients, there are splitting homomorphisms 
\begin{align*}
i : H^p(W_{j+1,0}, W_{j+1,0} \setminus \mathcal{F}_{j,0}^{j+1,0}) & \rightarrow H^p(W_{j+1}, W_{j+1,0} \setminus \mathcal{F}_{j,0}^{j+1,0}) \\
\text{and} \quad \pi : H^{p+m}(W_{j+1}, W_{j+1,0} \setminus \mathcal{F}_{j,0}^{j+1,0}) & \rightarrow H^{p+m}(W_{j+1}, W_{j+1,0})
\end{align*}
such that the cup product $\smallsmile \omega : H^p(W_{j+1}, W_{j+1,0} \setminus \mathcal{F}_{j,0}^{j+1,0}) \rightarrow H^{p+m}(W_{j+1}, W_{j+1,0} \setminus \mathcal{F}_{j,0}^{j+1,0})$ induces a homomorphism
\begin{multline*}
H^{p}(W_{j+1}, W_{j+1,0}) \oplus H^p(W_{j+1,0}, W_{j+1,0} \setminus \mathcal{F}_{j,0}^{j+1,0}) \\
 \stackrel{\smallsmile \omega_{j+1}}{\longrightarrow} H^{p+m}(W_{j+1}, W_{j+1,0}) \oplus H^{p+m}(W_{j+1,0}, W_{j+1,0} \setminus \mathcal{F}_{j,0}^{j+1,0}) .
\end{multline*}

Then the off-diagonal component
\begin{equation}\label{eqn:main-thm-upper-cup-product}
\begin{tikzcd}
H^p(W_{j+1,0}, W_{j+1,0} \setminus \mathcal{F}_{j,0}^{j+1,0}) \arrow{r}{c_\omega} & H^{p+m}(W_{j+1}, W_{j+1,0}) 
\end{tikzcd}
\end{equation}
of the cup product can be written in terms of the splitting homomorphisms as follows.

\begin{lemma}\label{lem:cup-product-main-thm}
Given $\omega \in H^m(Z)$, $\eta \in H^p(W_{j+1,0}, W_{j+1,0} \setminus \mathcal{F}_{j,0}^{j+1,0})$ and splitting homomorphisms as defined above, the cup product $\smallsmile \omega : H^*(Z) \rightarrow H^*(Z)$ induces a homomorphism \eqref{eqn:main-thm-upper-cup-product} given by
\begin{equation*}
c_\omega(\eta) = \pi( i(\eta) \smallsmile \omega_{j+1}) .
\end{equation*}
\end{lemma}

All of the homomorphisms in the diagrams \eqref{eqn:cup-product-exact-strata} and \eqref{eqn:cup-product-exact-local} are either induced by restriction to a subspace, or by cup product with a fixed class $\omega \in H^m(Z)$ or its restriction to a subspace. Therefore these homomorphisms commute with restriction, and so the following diagram, which is induced by the inclusion of triples $(W_{j+1}, W_{j+1,0}, W_{j+1,0} \setminus \mathcal{F}_{j,0}^{j+1,0}) \hookrightarrow (Z_{j+1}, Z_j, Z_{j-1})$, commutes. Equivalently, the cup product process of Lemma \ref{lem:cup-product-main-thm} is natural with respect to restriction of triples.

{\small \begin{equation}\label{eqn:restriction-cup-product-relative}
\begin{tikzcd}[column sep=0.5cm]
 &  & H^p(Z_{j+1}, Z_{j-1}) \arrow{r} \arrow{d}{\smallsmile \omega_{j+1}} & H^p(Z_j, Z_{j-1}) \arrow{r} \arrow{dd} & 0 \\
0 \arrow{r} & H^{p+m}(Z_{j+1}, Z_j) \arrow{r} \arrow[leftrightarrow]{dd}{\text{Thm. \ref{thm:main-thm-morse}}}[swap]{\cong} & H^{p+m}(Z_{j+1}, Z_{j-1})  &  \\
& & H^p(W_{j+1}, W_{j+1,0} \setminus \mathcal{F}_{j,0}^{j+1,0}) \arrow{r} \arrow{d}{\smallsmile \omega_{j+1}} & H^p(W_{j+1,0}, W_{j+1,0} \setminus \mathcal{F}_{j,0}^{j+1,0}) \arrow{r}  & 0 \\
0 \arrow{r} & H^{p+m}(W_{j+1}, W_{j+1,0}) \arrow{r} & H^{p+m}(W_{j+1}, W_{j+1,0} \setminus \mathcal{F}_{j,0}^{j+1,0}) 
\end{tikzcd}
\end{equation}
}

Lemma \ref{lem:splitting-compatible} shows that we can choose splittings of the exact sequence in the above diagram so that the diagram still commutes. Therefore Theorem \ref{thm:main-thm-morse} shows that the cup product construction of Lemma \ref{lem:cup-product-strata} is given (up to isomorphism) by first restricting $H^p(Z_j, Z_{j-1}) \rightarrow H^p(W_{j+1,0}, W_{j+1,0} \setminus \mathcal{F}_{j,0}^{j+1,0})$ to a pair localised around the upper critical set and the space of flow lines, and then applying the localised cup product construction of Lemma \ref{lem:cup-product-main-thm}. 

This is summarised in the following proposition.

\begin{proposition}\label{prop:cup-product-main-theorem}
Let $f : Z \rightarrow \mathbb{R}$ satisfy Conditions (1)--(5) of \cite{Wilkin19}, so there is a Morse filtration $\emptyset = Z_{-1} \subset Z_0 \subset \cdots \subset Z_n = Z$ and Theorem \ref{thm:main-thm-morse} applies. Use Lemma \ref{lem:splitting-compatible} to choose compatible splittings of each row in \eqref{eqn:restriction-cup-product-relative} and let $\omega \in H^m(Z)$ and let $\omega_j$ be the restriction to $H^m(Z_j)$. Then the cup product homomorphism $c(\omega) : H^p(Z_j, Z_{j-1}) \rightarrow H^{p+m}(Z_{j+1}, Z_j)$ of Lemma \ref{lem:cup-product-strata} factors through the bottom row in the following diagram.
\begin{equation}\label{eqn:cup-product-main-theorem}
\begin{tikzcd}[column sep=1cm]
H^p(Z_j, Z_{j-1}) \arrow[leftrightarrow]{d}{\cong}[swap]{\text{Thm.} \, \ref{thm:main-thm-morse}} \arrow{dr}{restriction} \arrow{rrr}{c_\omega \, \text{(Lem. \ref{lem:cup-product-strata})}} & & & H^{p+m}(Z_{j+1}, Z_j) \arrow[leftrightarrow]{d}{\cong}[swap]{\text{Thm.} \, \ref{thm:main-thm-morse}} \\
H^p(W_j, W_{j,0})  & H^p(W_{j+1,0}, W_{j+1,0} \setminus \mathcal{F}_{j,0}^{j+1,0}) \arrow{rr}{c_\omega \, \text{(Lem. \ref{lem:cup-product-main-thm})}} & & H^{p+m}(W_{j+1}, W_{j+1,0})  \\
\end{tikzcd}
\end{equation}
\end{proposition}

Therefore we see that the main theorem of Morse theory (Theorem \ref{thm:main-thm-morse}) allows us to interpret the cup product homomorphism via the maps in the bottom row of the above diagram. Equivalently, this expresses the cup product construction in terms of cohomology groups of spaces localised around the critical sets and spaces of flow lines.

In particular, up to the isomorphisms given by Theorem \ref{thm:main-thm-morse}, the cup product in $H^*(Z)$ corresponds to a homomorphism along the bottom row of the above diagram
\begin{equation}\label{eqn:local-cup-product-homomorphism}
H^p(W_j, W_{j,0}) \rightarrow H^p(W_{j+1,0}, W_{j+1,0} \setminus \mathcal{F}_{j,0}^{j+1,0}) \rightarrow H^{p+m}(W_{j+1}, W_{j+1,0}) .
\end{equation}

\begin{remark}\label{rem:local-cup-product-main-theorem}
It is sometimes useful to consider classes $\omega \in H^m(W_{j+1})$ that may not be in the image of the restriction $H^m(Z) \rightarrow H^m(W_{j+1})$. In this case we have the diagram
\begin{equation}\label{eqn:local-cup-product-main-theorem}
\begin{tikzcd}[column sep=1cm]
H^p(Z_j, Z_{j-1}) \arrow[leftrightarrow]{d}{\cong}[swap]{\text{Thm.} \, \ref{thm:main-thm-morse}} \arrow{dr}{restriction} \arrow[dashrightarrow]{rrr} & & & H^{p+m}(Z_{j+1}, Z_j) \arrow[leftrightarrow]{d}{\cong}[swap]{\text{Thm.} \, \ref{thm:main-thm-morse}} \\
H^p(W_j, W_{j,0})  & H^p(W_{j+1,0}, W_{j+1,0} \setminus \mathcal{F}_{j,0}^{j+1,0}) \arrow{rr}{c_\omega \, \text{(Lem. \ref{lem:cup-product-main-thm})}} & & H^{p+m}(W_{j+1}, W_{j+1,0})  \\
\end{tikzcd}
\end{equation}
for which there is an induced homomorphism $H^p(Z_j, Z_{j-1}) \rightarrow H^{p+m}(Z_{j+1}, Z_j)$, which may not come from a cup product in $H^*(Z)$. This is the case for the convolution homomorphism in Section \ref{subsec:cup-product-Thom-class}, where we can interpret this homomorphism as cup product on a product space $Z \times \Gr$, and the dashed arrow in the diagram above via pullback to this space.
\end{remark}

\section{A modified version of transversality}\label{sec:transversality}

The construction of the cup product and differentials from the previous two sections is valid for a (possibly singular) space $Z$ inside an ambient manifold $M$ satisfying the conditions of Theorem \ref{thm:main-thm-morse}. It then remains to interpret the topological meaning of the homomorphisms \eqref{eqn:differential-main-theorem} and \eqref{eqn:local-cup-product-homomorphism} in terms of critical sets and spaces of flow lines between them. In the Morse-Bott setting this requires the additional Morse-Bott-Smale transversality condition (cf. \cite{AustinBraam95}). 

In this section we return to the smooth space $M$ and show in Proposition \ref{prop:transversality-smooth-T2} that this transversality condition can be weakened so that the first homomorphism in \eqref{eqn:differential-main-theorem} and \eqref{eqn:local-cup-product-homomorphism} is isomorphic to pullback from the lower critical set $C_\ell$ to the space of flow lines $\mathcal{F}_{j,0}^{j+1,0}$ followed by cup product with a certain Euler class (cf. the diagrams \eqref{eqn:differential-T1-T2} and \eqref{eqn:cup-product-T1-T2}). The main example of interest is the space of representations of a quiver, for which this Euler class can be computed using Corollary \ref{cor:Euler-class}. 

It then remains to interpret the second homomorphism in \eqref{eqn:differential-main-theorem} and \eqref{eqn:local-cup-product-homomorphism}, which is the main result of Section \ref{sec:morse-convolution}, where we show that this can be expressed in terms of a pushforward homomorphism from the space of flow lines to the upper critical set.

\subsection{Modified transversality conditions}\label{subsec:modified-transversality}

First we recall the basic idea of the case where the function is Morse-Smale or Morse-Bott-Smale. These ideas have been developed by a number of authors such as Thom, Smale \cite{Smale60}, Witten \cite{Witten82} (see \cite{Bott88} for an overview) and the resulting Morse complex has been explained in detail by Austin and Braam \cite{AustinBraam95}.

Since the function is Morse or Morse-Bott, then the first page of the spectral sequence consists of terms
\begin{equation*}
H^p(M_j, M_{j-1}) \cong H^{p-\lambda_j}(C_j) .
\end{equation*}
The space of points on a flow line between the adjacent critical sets $C_j$ and $C_{j+1}$ is denoted $\mathcal{F}_{j,0}^{j+1,0}$ and the quotient by the $\mathbb{R}$-action of the flow is denoted $\tilde{\mathcal{F}}_{j,0}^{j+1,0}$, with associated upper and lower projection maps $\pi_\ell : \mathcal{F}_{j,0}^{j+1,0} \rightarrow C_j$, $\pi_u : \mathcal{F}_{j,0}^{j+1,0} \rightarrow C_{j+1}$ and $\tilde{\pi}_\ell : \tilde{\mathcal{F}}_{j,0}^{j+1,0} \rightarrow C_j$, $\tilde{\pi}_u : \tilde{\mathcal{F}}_{j,0}^{j+1,0} \rightarrow C_{j+1}$. The differential on the first page of the spectral sequence is then given by \cite[Thm. 3.1]{AustinBraam95}
\begin{align*}
d_1 : H^{p-\lambda_j}(C_j) & \rightarrow H^{p-\lambda_{j+1} + 1}(C_{j+1}) \\
 \eta & \mapsto (\tilde{\pi}_u)_* \tilde{\pi}_\ell^*(\eta) 
\end{align*}
and the cup product with $\omega \in H^m(M)$ is given by \cite[Thm. 3.11]{AustinBraam95}, which has a component mapping $H^*(C_j)$ to $H^*(C_{j+1})$ determined by
\begin{align*}
c(\omega) : H^{p-\lambda_j}(C_j) & \rightarrow H^{p-\lambda_{j+1}+m}(C_{j+1}) \\
 \eta & \mapsto (\pi_u)_* (\omega \smallsmile \pi_\ell^*(\eta))  .
\end{align*}

The proof given in \cite{AustinBraam95} uses de Rham cohomology, however our construction below will use singular cohomology since this behaves well when we pull back to a singular space in Theorem \ref{thm:cup-product-induces-convolution}.

The whole process described above works because the stable and unstable manifolds intersect transversely, however (as we have seen in Proposition \ref{prop:quivers-T2}) this is not satisfied for the norm-square of a moment map on the space of representations of a quiver, therefore we have to modify the transversality conditions to account for this.

Recall that if the stable and unstable manifolds intersect transversely, then the normal bundle of the stratum $S_j$ inside the ambient manifold $M$ restricts to the normal bundle of $\mathcal{F}_{j,0}^{j+1,0} = S_j \cap W_{j+1,0}$ inside $W_{j+1,0}$. In this section we will show that a modified form of the above construction of the differentials and cup product is still true if a weaker form of transversality holds, where the normal bundle of $\mathcal{F}_{j,0}^{j+1,0}$ inside $W_{j+1,0}$ is only a subbundle of the restriction of the normal bundle of $S_j$ to the space of flow lines.

\begin{definition}\label{def:transversality-conditions}

Let $M$ be a Riemannian manifold and $f : M \rightarrow \mathbb{R}$ a smooth function satisfying Conditions (1)--(5) of \cite{Wilkin19}. The spaces of flow lines satisfy \emph{weak transversality} if the following conditions hold. 

\begin{enumerate}

\item[{\bf (T1)}] The space of flow lines $\mathcal{F}_{j,0}^{j+1,0}$ has a tubular neighbourhood in $W_{j+1,0}$, denoted $D_j \rightarrow \mathcal{F}_{j,0}^{j+1,0}$. 

\item[{\bf (T2)}] The stratum $W_j^+$ has a tubular neighbourhood in $M$ denoted $\tilde{V}_j \rightarrow W_j^+$, which restricts to a disk bundle $V_j \rightarrow \mathcal{F}_{j,0}^{j+1,0}$ such that $D_j$ is a subbundle of $V_j$.

\end{enumerate}

\end{definition}

Of course, this definition is only useful if these conditions are satisfied for a class of interesting examples. Lemma \ref{lem:flow-lines-tubular-nbhd} and Proposition \ref{prop:quivers-T2} show that this is indeed the case for the Morse filtration of $\| \mu - \alpha \|^2$ on the space of representations of a quiver.

\subsection{Consequences of the transversality conditions}\label{subsec:cup-product-T1-T2}

The transversality conditions (T1) and (T2) allow us to carry out the following construction, which will be used in the next two sections to localise the differentials and cup products from the spectral sequence to spaces of flow lines.

First consider the pair $(M_j, M_{j-1})$. The convergence of the upwards flow determines two sub-pairs $(W_j, W_{j,0})$ and $(W_{j+1,0}, W_{j+1,0} \setminus \mathcal{F}_{j,0}^{j+1,0})$. \emph{A priori} these two pairs are not related, in fact they do not even intersect, however the main theorem of Morse theory does determine a homomorphism $H^p(W_j, W_{j,0}) \rightarrow H^p(W_{j+1,0}, W_{j+1,0} \setminus \mathcal{F}_{j,0}^{j+1,0})$ defined by composing the isomorphism $H^p(W_j, W_{j,0}) \cong H^p(M_j, M_{j-1})$ with the restriction homomorphism $H^p(M_j, M_{j-1}) \rightarrow H^p(W_{j+1,0}, W_{j+1,0} \setminus \mathcal{F}_{j,0}^{j+1,0})$.
\begin{equation}\label{eqn:main-thm-diagram}
\begin{tikzcd}
& H^p(M_j, M_{j-1}) \arrow[leftrightarrow]{dl}{\cong}[swap]{\text{Thm.} \, \ref{thm:main-thm-morse}} \arrow{dr} \\
H^p(W_j, W_{j,0})  & & H^p(W_{j+1,0}, W_{j+1,0} \setminus \mathcal{F}_{j,0}^{j+1,0}) 
\end{tikzcd}
\end{equation}

If we also impose Kirwan's minimal degeneracy condition \cite{Kirwan84}, then the Thom isomorphism shows that $H^p(W_j, W_{j,0}) \cong H^{p-\lambda_j}(C_j)$. The first transversality condition (T1) shows that there is also a Thom isomorphism $H^p(W_{j+1,0}, W_{j+1,0} \setminus \mathcal{F}_{j,0}^{j+1,0}) \cong H^p(D_j, D_j \setminus \mathcal{F}_{j,0}^{j+1,0}) \cong H^{p-\nu_j}(\mathcal{F}_{j,0}^{j+1,0})$, where $\nu_j$ is the codimension of $\mathcal{F}_{j,0}^{j+1,0} \subset V_{j+1,0}$. These maps fit together in the following diagram, where the dashed arrow in the bottom row is induced by the maps in the rest of the diagram.
\begin{equation}\label{eqn:T1-diagram}
\begin{tikzcd}
& H^p(M_j, M_{j-1}) \arrow[leftrightarrow]{dl}{\cong}[swap]{\text{Thm.} \, \ref{thm:main-thm-morse}}  \arrow{dr}{restriction} \\
H^p(W_j, W_{j,0}) \arrow[leftrightarrow]{d}{\cong}[swap]{Thom}  & & H^p(W_{j+1,0}, W_{j+1,0} \setminus \mathcal{F}_{j,0}^{j+1,0}) \arrow[leftrightarrow]{d}{\cong}[swap]{Thom} \\
H^{p-\lambda_j}(C_j)  \arrow[dashrightarrow]{rr} & & H^{p-\nu_j}(\mathcal{F}_{j,0}^{j+1,0})
\end{tikzcd}
\end{equation}

It now remains to understand the topological meaning of the homomorphism corresponding to the dashed arrow, which is the content of Proposition \ref{prop:transversality-smooth-T2} below.

The minimal degeneracy condition implies that the strata have normal bundles that restrict to the negative normal bundle on the critical set, and so we can use the Thom isomorphism $H^p(W_j, W_{j,0}) \cong H^{p-\lambda_j}(C_j)$. Moreover, on a smooth affine or projective variety this negative normal bundle extends to a bundle $\tilde{V}_j \rightarrow S_j$ over the entire stratum (cf. \cite[Sec. 4]{Kirwan84}) such that $H^p(W_j, W_{j,0}) \cong H^p(\tilde{V}_j, \tilde{V}_{j,0})$. Since the space of flow lines is contained in the stratum $\mathcal{F}_{j,0}^{j+1,0} \subset S_j$, then $\tilde{V}_j \rightarrow S_j$ restricts to a bundle over the space of flow lines, which we denote by $V_j \rightarrow \mathcal{F}_{j,0}^{j+1,0}$. 
\begin{equation*}
\begin{tikzcd}[column sep=15mm]
H^p(W_j, W_{j,0}) \arrow[leftrightarrow]{r}{\cong} \arrow[leftrightarrow]{d}{\cong}[swap]{Thom} & H^p(\tilde{V}_j, \tilde{V}_{j,0}) \arrow{r} \arrow[leftrightarrow]{d}{\cong}[swap]{Thom} & H^p(V_j, V_{j,0}) \arrow[leftrightarrow]{d}{\cong}[swap]{Thom} \\
H^{p-\lambda_j}(C_j) \arrow[leftrightarrow]{r}{\cong} \arrow[bend right = 15]{rr}[swap]{pullback} & H^{p-\lambda_j}(S_j) \arrow{r}{restriction} & H^{p-\lambda_j}(\mathcal{F}_{j,0}^{j+1,0}) 
\end{tikzcd}
\end{equation*}

Now if transversality condition (T2) is satisfied, then the tubular neighbourhood of $\mathcal{F}_{j,0}^{j+1,0}$ inside $W_{j+1,0}$ is a disk bundle $D \rightarrow \mathcal{F}_{j,0}^{j+1,0}$, which is a subbundle of $V_j$. Then Lemma \ref{lem:Euler-subbundle} shows that we can augment diagram \eqref{eqn:T1-diagram} with a homomorphism corresponding to cup product with the Euler class of the normal bundle of $D \subset V_j$.
\begin{equation}\label{eqn:T2-diagram}
\begin{tikzcd}
& H^p(M_j, M_{j-1})  \arrow[leftrightarrow]{dl}{\cong}[swap]{\text{Thm.} \, \ref{thm:main-thm-morse}}  \arrow{dr}{restriction} \arrow[bend left=10]{drr}{restriction} \\
H^p(W_j, W_{j,0}) \arrow[leftrightarrow]{d}{\cong}[swap]{Thom} \arrow[leftrightarrow]{r}{\cong} & H^p(\tilde{V}_j, \tilde{V}_{j,0}) \arrow[leftrightarrow]{d}{\cong}[swap]{Thom} \arrow{r} & H^p(V_j, V_{j,0}) \arrow{r} \arrow[leftrightarrow]{d}{\cong}[swap]{Thom} & H^p(W_{j+1,0}, W_{j+1,0} \setminus \mathcal{F}_{j,0}^{j+1,0}) \arrow[leftrightarrow]{d}{\cong}[swap]{Thom} \\
H^{p-\lambda_j}(C_j) \arrow[leftrightarrow]{r}{\cong} & H^{p-\lambda_j}(S_j) \arrow{r}{restriction}  & H^{p-\lambda_j}(\mathcal{F}_{j,0}^{j+1,0}) \arrow{r}{\smallsmile e} & H^{p-\nu_j}(\mathcal{F}_{j,0}^{j+1,0})
\end{tikzcd}
\end{equation}

Therefore we see that, up to the isomorphisms $H^{p-\nu_j}(\mathcal{F}_{j,0}^{j+1,0}) \cong H^p(W_{j+1,0}, W_{j+1,0} \setminus \mathcal{F}_{j,0}^{j+1,0})$ and $H^{p-\lambda_j}(C_j) \cong H^p(W_j, W_{j,0})$ (the left and right hand columns of the above diagram), the composition of homomorphisms
\begin{equation*}
H^p(W_j, W_{j,0}) \stackrel{\cong}{\rightarrow} H^p(M_j, M_{j-1}) \rightarrow H^p(W_{j+1,0}, W_{j+1,0} \setminus \mathcal{F}_{j,0}^{j+1,0}),
\end{equation*}
which appears in the construction of the differential from Proposition \ref{prop:differential-main-theorem} and the cup product from Proposition \ref{prop:cup-product-main-theorem}, factors through the composition of homomorphisms 
\begin{equation*}
\begin{tikzcd}
H^{p-\lambda_j}(C_j) \arrow{r}{pullback} & H^{p-\lambda_j}(\mathcal{F}_{j,0}^{j+1,0}) \arrow{r}{\smallsmile e} & H^{p-\nu_j}(\mathcal{F}_{j,0}^{j+1,0})
\end{tikzcd}
\end{equation*}
from the bottom row of the diagram \eqref{eqn:T2-diagram}. The above result is summarised in the following proposition.

\begin{proposition}\label{prop:transversality-smooth-T2}
Let $M$ be a manifold and $f : M \rightarrow \mathbb{R}$ a minimally degenerate smooth function satisfying Conditions (1)--(5) of \cite{Wilkin19} as well as the transversality conditions (T1) and (T2), and let $\emptyset = M_{-1} \subset M_0 \subset \cdots \subset M_n = M$ be the associated Morse filtration. Then the following diagram commutes.
\begin{equation}\label{eqn:singular-T2-diagram}
\begin{tikzcd}
& H^p(M_j, M_{j-1}) \arrow[leftrightarrow]{dl}{\cong}[swap]{\text{Thm.} \, \ref{thm:main-thm-morse}} \arrow{dr}{restriction} \\
H^{p}(W_j, W_{j,0}) & & H^p(W_{j+1,0}, W_{j+1} \setminus \mathcal{F}_{j,0}^{j+1,0})  \\
H^{p-\lambda_j}(C_j)  \arrow[leftrightarrow]{u}{\text{Thom}}[swap]{\cong} \arrow{r}{pullback} & H^{p-\lambda_j}(\mathcal{F}_{j,0}^{j+1,0}) \arrow{r}{\smallsmile e} & H^{p-\nu_j}(\mathcal{F}_{j,0}^{j+1,0}) \arrow[leftrightarrow]{u}{\text{Thom}}[swap]{\cong} 
\end{tikzcd}
\end{equation}
\end{proposition}

Therefore the diagram \eqref{eqn:differential-main-theorem} for the differential can be augmented with another row corresponding to the critical sets and spaces of flow lines

\begin{equation}\label{eqn:differential-T1-T2}
\begin{tikzcd}[column sep=6mm]
H^p(M_j, M_{j-1}) \arrow{dr}{restriction} \arrow[leftrightarrow]{d}{\cong}[swap]{\text{Thm.} \, \ref{thm:main-thm-morse}} \arrow{rrr}{d} & & &  H^{p+1}(M_{j+1}, M_j) \arrow[leftrightarrow]{d}{\cong}[swap]{\text{Thm.} \,  \ref{thm:main-thm-morse}} \\
H^p(W_j, W_{j,0})  & H^p(W_{j+1,0}, W_{j+1,0} \setminus \mathcal{F}_{j,0}^{j+1,0}) \arrow{r} & H^p(W_{j+1,0}) \arrow{r} & H^{p+1}(W_{j+1}, W_{j+1,0}) \\
 & H^{p-\nu_j}(\mathcal{F}_{j,0}^{j+1,0}) \arrow[leftrightarrow]{u}{\text{Thom}}[swap]{\cong} \arrow[dashrightarrow]{rr} & & H^{p+1-\lambda_{j+1}}(C_{j+1}) \arrow[leftrightarrow]{u}{\text{Thom}}[swap]{\cong} \\
H^{p-\lambda_j}(C_j) \arrow[leftrightarrow]{uu}{\text{Thom}}[swap]{\cong} \arrow{r}{pullback} & H^{p-\lambda_j}(\mathcal{F}_{j,0}^{j+1,0}) \arrow{u}{\smallsmile e}  & & 
\end{tikzcd}
\end{equation}
and the diagram \eqref{eqn:cup-product-main-theorem} for the cup product has an analogous augmentation
\begin{equation}\label{eqn:cup-product-T1-T2}
\begin{tikzcd}[column sep=1cm]
H^p(M_j, M_{j-1}) \arrow[leftrightarrow]{d}{\cong}[swap]{\text{Thm.} \, \ref{thm:main-thm-morse}} \arrow{dr}{restriction} \arrow{rrr}{c(\omega) \, \text{(Lem. \ref{lem:cup-product-strata})}} & & & H^{p+m}(M_{j+1}, M_j) \arrow[leftrightarrow]{d}{\cong}[swap]{\text{Thm.} \, \ref{thm:main-thm-morse}} \\
H^p(W_j, W_{j,0})  & H^p(W_{j+1,0}, W_{j+1,0} \setminus \mathcal{F}_{j,0}^{j+1,0}) \arrow{rr}{c(\omega) \, \text{(Lem. \ref{lem:cup-product-main-thm})}} & & H^{p+m}(W_{j+1}, W_{j+1,0})  \\
 & H^{p-\nu_j}(\mathcal{F}_{j,0}^{j+1,0}) \arrow[leftrightarrow]{u}{\text{Thom}}[swap]{\cong} \arrow[dashrightarrow]{rr} & & H^{p+m-\lambda_{j+1}}(C_{j+1}) \arrow[leftrightarrow]{u}{\text{Thom}}[swap]{\cong} \\
H^{p-\lambda_j}(C_j) \arrow[leftrightarrow]{uu}{\text{Thom}}[swap]{\cong} \arrow{r}{pullback} & H^{p-\lambda_j}(\mathcal{F}_{j,0}^{j+1,0}) \arrow{u}{\smallsmile e}  & & 
\end{tikzcd}
\end{equation}

\begin{remark}
Lemma \ref{lem:equiv-Euler-general} and Corollary \ref{cor:Euler-class} show that the Euler class in the above diagrams can be explicitly computed in terms of the critical point data.
\end{remark}

\section{Relationship with convolution in Borel-Moore homology}\label{sec:BM-convolution}

It now remains to interpret the topological meaning of the dashed arrow in diagrams \eqref{eqn:differential-T1-T2} and \eqref{eqn:cup-product-T1-T2}. The aim of this section is to explain how how this is related to pushforward in Borel-Moore homology. There are two cases

\begin{itemize}

\item pushforward on a sphere bundle (Lemma \ref{lem:pushforward-sphere-bundle}), which will correspond to the dashed arrow in the diagram \eqref{eqn:differential-T1-T2} for the differential (Proposition \ref{prop:differential-pushforward}), and 

\item pushforward on a projective bundle (Lemma \ref{lem:pushforward-projective-bundle}) which will correspond to the dashed arrow in the diagram \eqref{eqn:cup-product-T1-T2} (Proposition \ref{prop:cup-product-pushforward}).

\end{itemize}

\subsection{Pullback and pushforward in Borel-Moore homology}\label{subsec:BM-pullback-pushforward}

Let $M_1, M_2$ be manifolds of dimension $d_1, d_2$ respectively, and let $N \subset M_1 \times M_2$ be an embedded submanifold of codimension $d$ such that the projection $N \stackrel{i}{\hookrightarrow} M_1 \times M_2 \stackrel{p_2}{\longrightarrow} M_2$ is proper. Then, following \cite[Sec. 2.7]{ChrissGinzburg97}, there is a pullback map in Borel-Moore homology
\begin{align}\label{eqn:BM-pullback}
\begin{split}
H_p^{BM}(M_2) & \rightarrow H_{p+d_1-d}^{BM}(N) \\
c & \mapsto [N] \cap ([M_1] \boxtimes c)
\end{split}
\end{align}
where $\cap : H_*^{BM}(M_1 \times M_2) \rightarrow H_{*-d}^{BM}(N)$ denotes restriction with supports (cf. \cite[Sec. 2.6.21]{ChrissGinzburg97}) and $\boxtimes$ denotes the K\"unneth isomorphism $H_*^{BM}(M_1) \otimes H_*^{BM}(M_2) \rightarrow H_*^{BM}(M_1 \times M_2)$. 

First note that the homomorphism $H_p^{BM}(M_2) \rightarrow H_{p+d_1}^{BM}(M_1 \times M_2)$ given by $c \mapsto [M_1] \boxtimes c$ is Poincar\'e dual to the pullback in cohomology $H^{d_2-p}(M_2) \rightarrow H^{d_2-p}(M_1 \times M_2)$. By definition (cf. \cite[Sec. 2.6.21]{ChrissGinzburg97}), restriction with supports to a submanifold $H_k^{BM}(M_1 \times M_2) \rightarrow H_{k-d}^{BM}(N)$ is Poincar\'e dual to restriction in cohomology, and therefore the composition \eqref{eqn:BM-pullback} of these two homomorphisms in Borel-Moore homology is Poincar\'e dual to the composition of these two pullback homomorphisms. This is summarised in the following lemma.

\begin{lemma}\label{lem:PD-pullback}
The following diagram commutes, where the top row is the pullback homomorphism \eqref{eqn:BM-pullback}, the bottom row is pullback in ordinary cohomology, and the vertical arrows are given by Poincar\'e duality.
\begin{equation}
\begin{tikzcd}
H_p^{BM}(M_2) \arrow{r} \arrow{d}{\cong}[swap]{P.D.} & H_{p+d_1-d}^{BM}(N) \arrow{d}{\cong}[swap]{P.D.} \\
H^{d_2-p}(M_2) \arrow{r} & H^{d_2-p}(N)
\end{tikzcd}
\end{equation}
\end{lemma}

In order to relate the Poincar\'e dual of pushforward in Borel-Moore homology to the Morse-theoretic constructions in the subsequent sections, we need the following lemmas. The proofs use well-known results from the theory, but the statements are needed for the connection with Morse theory in the next section, so we state and prove everything here for completeness. The first is for the pushforward of a submanifold by the inclusion map.
\begin{lemma}\label{lem:pushforward-submanifold}
Let $M$ be a manifold of dimension $n$ and $i : N \hookrightarrow M$ a submanifold of dimension $s$ with a tubular neighbourhood $U$. Then the pushforward $i_* : H_p^{BM}(N) \rightarrow H_p^{BM}(M)$ is Poincar\'e dual to the composition of the following homomorphisms
\begin{equation*}
H^{s}(N) \stackrel{\cong}{\longrightarrow} H^{n}(U, U \setminus N) \stackrel{\cong}{\longrightarrow} H^{n}(M, M \setminus N) \rightarrow H^{n}(M) ,
\end{equation*}
where the first isomorphism is the Thom isomorphism, the second isomorphism is excision and the final homomorphism comes from the long exact sequence of the pair $(M, M \setminus N)$.  
\end{lemma}

\begin{proof}
This is a special case of the pushforward from the long exact sequence of Borel-Moore homology for the triple $N \hookrightarrow M \hookleftarrow M \setminus N$, which is defined using Poincar\'e duality (cf. \cite[Sec. 2.6.9]{ChrissGinzburg97}).
\end{proof}

The next lemma is for the case of a sphere bundle associated to a vector bundle. 
\begin{lemma}\label{lem:pushforward-sphere-bundle}
Let $B$ be a manifold of dimension $n$, let $V \rightarrow B$ be a real vector bundle of rank $r$, define $V_0 = V \setminus B$ and let $\pi : S \rightarrow B$ denote the associated sphere bundle. Then the pushforward $\pi_* : H_p^{BM}(S) \rightarrow H_p^{BM}(B)$ is Poincar\'e dual to the composition of the following homomorphisms
\begin{equation*}
H^{n+r-1-p}(S) \stackrel{\cong}{\longrightarrow} H^{n+r-1-p}(V_0) \rightarrow H^{n+r-p}(V, V_0) \stackrel{\cong}{\longrightarrow} H^{n-p}(B)
\end{equation*}
where the first homomorphism is homotopy equivalence, the second is the connecting homomorphism in the long exact sequence of the pair $(V, V_0)$ and the third is the pushforward from the Thom isomorphism.
\end{lemma}

\begin{proof}
This follows from the fact that, on a smooth manifold, the pushforward in the Thom-Gysin sequence is given by integration over the fibres (cf. \cite[Prop. 14.33]{BottTu82}), which is Poincar\'e dual to the pushforward in Borel-Moore homology.
\end{proof}

Let $B$ be a manifold of dimension $\dim_\mathbb{R} B = n$ and let $\pi : E \rightarrow B$ be a smooth fibre bundle with smooth compact fibres diffeomorphic to $F$ with $\dim_\mathbb{R} F = r$. Suppose that the conditions of the Leray-Hirsch theorem are satisfied, so that there exist classes $\zeta_1, \ldots, \zeta_m \in H^*(E)$ such that for any fibre $i : F \hookrightarrow E$ the classes $i^*(\zeta_1), \ldots, i^*(\zeta_m)$ generate $H^*(F)$ as a group. In the following we choose the ordering so that $i^*(\zeta_m)$ generates the top dimensional cohomology $H^r(F)$. For example, these conditions are satisfied when $E$ is the projectivisation of a complex vector bundle, where the classes $\zeta_1, \ldots, \zeta_m$ are powers of the Chern classes of the tautological line bundle over $E$.

The following lemma describes the pushforward map in this setting.

\begin{lemma}\label{lem:pushforward-projective-bundle}
The pushforward in Borel-Moore homology $\pi_* : H_p^{BM}(E) \rightarrow H_p^{BM}(B)$ is Poincar\'e dual to the pushforward by fibre integration
\begin{align}\label{eqn:cohomology-pushforward}
\begin{split}
\pi_* : H^{n+r-p}(E) & \rightarrow H^{n-p}(B) \\
\pi_* \left( \pi^*(\beta) \smallsmile \zeta_k \right) & = \begin{cases} \beta & \text{if $k=m$} \\ 0 & \text{if $0 \leq k < m$} \end{cases}
\end{split}
\end{align}
\end{lemma}

\begin{proof}
Let $b \in H_p^{BM}(B)$ be Poincar\'e dual to $\beta \in H^{n-p}(B)$ and let $c_k \in H_{n+r-\ell}^{BM}(E)$ be Poincar\'e dual to $\zeta_k \in H^\ell(E)$ for $k = 1, \ldots, m$. 

Since $E$ and $B$ are manifolds and the fibres of $\pi : E \rightarrow B$ are compact manifolds, then pullback to a fibre bundle in Borel-Moore homology is Poincar\'e dual to pullback in ordinary cohomology (cf. \cite[p102]{ChrissGinzburg97}), and so 
\begin{equation*}
\pi^*(b) \in H_{p+r}^{BM}(E) \quad \text{is Poincar\'e dual to} \quad \pi^*(\beta) \in H^{n-p}(E) .
\end{equation*} 
The intersection pairing in Borel-Moore homology is Poincar\'e dual to cup product in cohomology (cf. \cite[Sec. 2.6.15]{ChrissGinzburg97}), therefore
\begin{equation*}
\pi^*(b) \smallfrown c_k \quad \text{is Poincar\'e dual to} \quad \pi^*(\beta) \smallsmile \zeta_k \quad \text{for all $k = 1, \ldots, m$.}
\end{equation*} 
The following special case of the projection formula in Borel-Moore homology
\begin{equation*}
\pi_*(\pi^*(b) \smallfrown c_k) = b \smallfrown \pi_*(c_k) \quad \text{where $b \in H_*^{BM}(B)$ and $c_k \in H_*^{BM}(E)$} ,
\end{equation*}
is valid when $E$ and $B$ are manifolds (cf. \cite[(2.6.29]{ChrissGinzburg97}).

If $\zeta_k \in H^\ell(E)$ for $k < m$, then $\ell < r$ (since the fibre $F$ is smooth and compact). Then the Poincar\'e dual $c_k \in H_{n+r-\ell}^{BM}(E)$ satisfies $\pi_*(c_k) = 0$ for dimensional reasons, therefore $\pi_*(\pi^*(b) \smallfrown c_k) = b \smallfrown \pi_*(c_k) = 0$.

If $k=m$, then $\pi^*(\beta) \smallsmile \zeta_m \in H^{n+r-p}(E)$ is Poincar\'e dual to $\pi^*(b) \smallfrown c_m \in H_p^{BM}(E)$ and so the projection formula shows that the pushforward is $\pi_*(\pi^*(b) \smallfrown c_m) = b \smallfrown \pi_*(c_m) \in H_p^{BM}(B)$. 

Therefore it remains to show that $\pi_*(c_m) = [B]$, which we will do by showing that $\pi_*(c_m)$ restricts to the fundamental class on each trivialisation $\pi^{-1}(U) \rightarrow U$. First recall that restriction to an open subset in Borel-Moore homology is induced from pushforward in ordinary relative homology (cf. \cite[Sec. 2.6.9]{ChrissGinzburg97}) which appears as the horizontal homomorphisms in the following diagram
 \begin{equation}\label{eqn:pushforward-restriction-commutes}
\begin{tikzcd}
H_*^{BM}(E) \cong H_*(\overline{E}, \overline{E} \setminus E) \arrow{r}{j_*} \arrow{d}{\pi_*} & H_*(\overline{E}, \overline{E} \setminus \pi^{-1}(U)) \cong H_*^{BM}(\pi^{-1}(U)) \arrow{d}{\pi_*} \\
H_*^{BM}(B) \cong H_*(\overline{B}, \overline{B} \setminus B) \arrow{r}{j_*} & H_*(\overline{B}, \overline{B} \setminus U) \cong H_*^{BM}(U)
\end{tikzcd}
\end{equation}
where $j$ is used to denote both the inclusion $(\bar{E}, \bar{E} \setminus E) \hookrightarrow (\bar{E}, \bar{E} \setminus \pi^{-1}(U))$ and the inclusion $(\bar{B}, \bar{B} \setminus B) \hookrightarrow (\bar{B}, \bar{B} \setminus U)$. Each homomorphism is induced from a continuous map of pairs, therefore the diagram commutes.

Restriction in Borel-Moore homology can also be defined using Poincar\'e duality as part of the construction of the long exact sequence for Borel-Moore homology (cf. \cite[Sec. 2.6.9]{ChrissGinzburg97}). Since all the spaces $E$, $B$, $\pi^{-1}(U)$ and $U$ are manifolds, then the top and bottom rows in the above diagram are induced via pullback in cohomology
\begin{equation*}
\begin{tikzcd}
H_*^{BM}(E) \cong H^{n+r-*}(E) \arrow{r}{j^*} & H^{n+r-*}(\pi^{-1}(U)) \cong H_*^{BM}(\pi^{-1}(U)) \\
H_*^{BM}(B) \cong H^{n-*}(B) \arrow{r}{j^*} & H^{n-*}(U) \cong H_*^{BM}(U)
\end{tikzcd}
\end{equation*}
In particular, the classes $\zeta_k \in H^{n+r-*}(E)$ restrict to classes $\xi_k \in H^{n+r-*}(\pi^{-1}(U))$ for which the Leray-Hirsch property implies that they restrict to generators of the cohomology of each fibre. Since $\pi^{-1}(U) \cong U \times F \rightarrow U$ is trivial, then the pushforward in Borel-Moore homology has an explicit description, which maps the Poincar\'e dual of $\xi_m$ to the fundamental class $[U]$. Therefore in diagram \eqref{eqn:pushforward-restriction-commutes} we have $\pi_* \circ j_*(c_m) = [U]$. Commutativity of \eqref{eqn:pushforward-restriction-commutes} then implies $[U] = j_* \circ \pi_*(c_m)$, and therefore $\pi_*(c_m)$ restricts to the fundamental class $[U]$ of each trivialisation. This is the defining property of the fundamental class $[B]$, which must therefore be equal to $\pi_*(c_m)$. Therefore
\begin{equation*}
\pi_*(\pi^*(b) \smallfrown c_m) = b \smallfrown [B] = b
\end{equation*}
is Poincar\'e dual to $\beta = \pi_*(\pi^*(\beta) \smallsmile \zeta_m)$ defined using \eqref{eqn:cohomology-pushforward}, which completes the proof. 
\end{proof}

Now we consider a special case of the above construction when the fibre bundle is the projectivisation of a complex vector bundle. In this setting we will show that the pushforward is given by a construction analogous to that for the cup product from in Sections \ref{subsec:stratified-cup-product} and \ref{subsec:main-thm-cup-product}.

Let $B$ be a manifold, let $V \rightarrow B$ be a complex vector bundle of complex rank $r$, and let $V_0$ denote the vector bundle minus the zero section. Suppose that there is an action of $S^1$ on $V$ fixing the zero section such that the projection $\pi : V \rightarrow B$ is $S^1$-equivariant and the action on each fibre is linear with weight one. The Atiyah-Bott lemma then implies that the horizontal rows of the following diagram are exact
\begin{equation}\label{eqn:AB-diagram}
\begin{tikzcd}
0 \arrow{r} & H_{S^1}^p(V, V_0) \arrow{r} \arrow{d}{\cong} & H_{S^1}^p(V) \arrow{r} \arrow{d}{\cong} & H_{S^1}^p(V_0) \arrow{r} \arrow{d}{\cong} & 0 \\
0 \arrow{r} & H_{S^1}^{p-2r}(B)\arrow{r}{\smallsmile e} & H_{S^1}^p(B) \arrow{r}{\pi^*} & H^p(\mathbb{P} V) \arrow{r} & 0
\end{tikzcd}
\end{equation}
where $e$ denotes the $S^1$-equivariant Euler class of $V$. Since $e \in H_{S^1}^{2r}(B)$ then it follows from the exactness of the bottom row of the above diagram that $H_{S^1}^p(B) \cong H^p(\mathbb{P} V)$ for all $p < 2r-1$.

Let $\xi \in H^2(\mathbb{P}V) \cong H_{S^1}^2(V_0)$ denote the first Chern class of the projective bundle $\pi : \mathbb{P}V \rightarrow B$. For each $k = 0, \ldots, r-1$, the restriction of $\xi^k$ to any fibre $\mathbb{P}^{r-1} \cong \mathbb{P}V_x$ is a generator of $H^{2k}(\mathbb{P}^{r-1})$. The Leray-Hirsch theorem then determines an isomorphism of groups 
\begin{align}\label{eqn:Leray-Hirsch}
\begin{split}
\sum_{\ell+2k = m} H^\ell(B) \otimes H^{2k}(\mathbb{P}^{r-1}) & \rightarrow H^m(\mathbb{P}V) \\
\sum \beta \otimes i^*(\xi^k) & \mapsto \pi^*(\beta) \smallsmile \xi^k  .
\end{split}
\end{align}

For $k = 1, \ldots, r$, define $\gamma_k := \pi^* c_k(V) \in H^{2k}(\mathbb{P}V)$. Recall that, since the action of $S^1$ fixes the zero section and acts freely with weight one on the fibres, then $H_{S^1}^*(B) \cong H^*(B)[\xi]$ and $H_{S^1}^*(V_0) \cong H^*(\mathbb{P}V) \cong H^*(B)[\xi] / \sim$ with the relation 
\begin{equation}\label{eqn:projective-relation}
\xi^r + c_1(V) \smallsmile \xi^{r-1} + \cdots + c_r(V) \sim 0 ,
\end{equation}
where the isomorphism is given by pullback $H_{S^1}^*(B) \rightarrow H_{S^1}^*(V_0) \cong H^*(\mathbb{P}V)$ (cf. \cite[Ch. VII]{ACGH85}). 

Therefore, in the exact sequence \eqref{eqn:AB-diagram}, we see that the equivariant Euler class is
\begin{equation*}
e(V) = \pm \left( \xi^r + c_1(V) \xi^{r-1} + \cdots + c_r(V) \right) \in H_{S^1}^*(B) \cong H^*(B)[\xi] .
\end{equation*}

This Euler class and the Leray-Hirsch theorem can be used to construct an explicit splitting of the bottom row of \eqref{eqn:AB-diagram} as follows. There is a splitting that respects the group structure in the graded cohomology rings
\begin{equation}\label{eqn:relative-splitting-inclusion}
\begin{tikzcd}
i : H^*(\mathbb{P}V) \cong H^*(B)[\xi]/\sim \arrow{r} & H_{S^1}^*(B) \cong H^*(B)[\xi]
\end{tikzcd}
\end{equation}
given by first using the Leray-Hirsch isomorphism, and then defining a homomorphism of groups
\begin{align*}
H^*(\mathbb{P}V) \cong H^*(B) \otimes H^*(\mathbb{P}^{r-1}) & \rightarrow H_{S^1}^*(B) \cong H^*(B) \otimes H^*(\mathsf{BU}(1)) \\
i(\pi^*(\beta) \smallsmile \xi^k) & = \beta \otimes \xi^k .
\end{align*}

Now we define the splitting homomorphism for the map $\smallsmile e$ in \eqref{eqn:AB-diagram}
\begin{equation}\label{eqn:relative-splitting-projection}
pr : H_{S^1}^p(B) \rightarrow H_{S^1}^{p-2r}(B)
\end{equation}
by contracting with the Euler class as follows. Given $\eta = \beta \otimes \xi^k \in H_{S^1}^*(B) \cong H^*(B) \otimes H^*(\mathsf{BU}(1))$, if $k \geq r$ then we can rewrite
\begin{equation*}
\eta = (\beta \otimes \xi^{k-r}) \otimes \xi^r = (\beta \otimes \xi^{k-r}) \smallsmile e + \text{lower order terms}
\end{equation*}
where ``lower order'' means that the power of $\xi$ is less than $k$. Repeating this process, we obtain $\eta = \eta' \smallsmile e + \eta''$, where $\eta''$ is a sum of terms of the form $\beta'' \otimes \xi^k$ with $k \leq r-1$. Now define
\begin{equation*}
pr(\eta) = \eta'
\end{equation*}
and note that $pr \circ i = 0$, therefore this is a well-defined splitting of the sequence \eqref{eqn:AB-diagram}.

From Lemma \ref{lem:pushforward-projective-bundle}, the pushforward $H^*(\mathbb{P}V) \rightarrow H^*(B)$ which is Poincar\'e dual to the pushforward in Borel-Moore homology is then given by
\begin{equation}\label{eqn:pushforward}
\pi_*(\pi^*(\beta) \smallsmile \xi^k) \mapsto \beta \smallsmile \pi_*(\xi^k) = \begin{cases} \beta & k = r-1 \\ 0 & \text{otherwise}. \end{cases}
\end{equation}

The following lemma shows that this can be interpreted in terms of the cup product construction given above.

\begin{lemma}\label{lem:pushforward-cup-product}
With the same notation and conditions as above, the pushforward $H^p(\mathbb{P}V) \rightarrow H^{p-2r}(B)$ is given by
\begin{equation}\label{eqn:pushforward-cup-product}
\pi_*(\eta) = pr(i(\eta) \otimes \xi) ,
\end{equation}
where $pr$ and $i$ are the splitting homomorphisms defined above for the exact sequence \eqref{eqn:AB-diagram}.
\end{lemma}

\begin{proof}
From \eqref{eqn:pushforward}, the pushforward is nonzero if and only if $\eta = \pi^*(\beta) \smallsmile \xi^{r-1}$, in which case $\pi_*(\eta) = \beta$ and $i(\eta) \otimes \xi = \beta \otimes \xi^r \in H_{S^1}^*(B)$. Therefore the above construction of the splitting $pr$ shows that $pr(i(\eta) \otimes \xi) = pr(\beta \otimes \xi^r) = \beta = \pi_*(\eta)$. 

In the case that $\eta = \pi^*(\beta) \smallsmile \xi^k$ with $k < r-1$, then $\pi_*(\eta) = 0$ and $pr(\beta \otimes \xi^{k+1}) = 0$.
\end{proof}

\begin{remark}\label{rem:Chern-class}
Therefore we see that the pushforward is given by a process analogous to that of the cup product construction of Lemma \ref{lem:cup-product-main-thm}. This is made precise in Proposition \ref{prop:cup-product-pushforward} below, which shows that the cup product on the Morse complex factors through this pushforward. In the proof of Proposition \ref{prop:cup-product-pushforward} we will use the above general construction in the following context of representations of quivers.

Consider two adjacent critical sets $C_\ell$ (lower) and $C_u$ (upper) with associated dimension vectors ${\bf v_\ell}$ and ${\bf v_u}$ from \cite[Prop. 3.13]{Wilkin17} and let $W_u^-$ be the unstable bundle of $C_u$ (equivalently, use the homeomorphism of Theorem \ref{thm:slice-homeo} to replace $W_u^-$ with the negative slice bundle $S_u^-$ in the following). Suppose also that ${\bf v} = {\bf v_\ell}$ so that $C_\ell$ minimises $\| \mu - \alpha \|^2 : \Rep(Q, {\bf v}) \rightarrow \mathbb{R}$. Then $C_\ell / K_{\bf v_\ell} \cong \mathcal{M}(Q, {\bf v_\ell})$ and $C_u / K_{\bf v_u}$ is homotopy equivalent to $\mathcal{M}(Q, {\bf v_u})$ by Corollary \ref{cor:critical-level-isomorphism}. In equivariant cohomology we have
\begin{equation*}
H_{K_{\bf v_\ell}}^*(C_\ell) \cong H^*(\mathcal{M}(Q, {\bf v_\ell})) \otimes H^*(\BU(1)) \quad \text{and} \quad H_{K_{\bf v_\ell}}^*(C_u) \cong H^*(\mathcal{M}(Q, {\bf v_u})) \otimes H^*(\BU(1))^{\otimes 2}
\end{equation*}
since the scalar multiples of the identity in $K_{\bf v_\ell}$ act trivially on $\Rep(Q, {\bf v_\ell})$ and each upper critical point determines a reduction of structure group from $K_{\bf v_\ell}$ to $K_{\bf v_u} \times \U(1)$.

This second factor of $\U(1)$ acts freely on $S_{u,0}^-$ and the projection $S_{u,0}^- \rightarrow C_u$ is $K_{\bf v_\ell}$-equivariant, so this bundle descends to the quotient by $K_{\bf v_\ell}$ 
\begin{equation*}
S_{u,0}^-/K_{\bf v_\ell} \rightarrow \mathcal{M}(Q, {\bf v_u}) ,
\end{equation*}
with an induced map in cohomology
\begin{equation*}
H^*(\mathcal{M}(Q, {\bf v_u})) \otimes H^*(\BU(1)) \rightarrow H^*(S_{u,0}^- / K_{\bf v_\ell}) .
\end{equation*}

Lemma \ref{lem:cup-product-pushforward} and Proposition \ref{prop:cup-product-pushforward} use the first Chern class of the projective bundle $\xi \in H^2(S_{u,0}^- / K_{\bf v_\ell})$.
\end{remark}

\begin{remark}\label{rem:image-pushforward}
\emph{A priori} the pushforward lies in $H_{K_{\bf v_\ell}}^*(C_u) \cong H^*(\mathcal{M}(Q, {\bf v_u})) \otimes H^*(\BU(1)$, however we see from the above construction that the image of \eqref{eqn:pushforward} lies in $H^*(\mathcal{M}(Q, {\bf v_u}))$. This will be important in the final step of the convolution construction of \eqref{eqn:full-convolution-cup-product-diagram}, where the image of the convolution lies in $H^*(\mathcal{M}(Q, {\bf v_u}))$.
\end{remark}

\section{The differential and cup product as convolution operators}\label{sec:morse-convolution}

The aim of this section is to relate convolution in Borel-Moore homology from the previous section to the differential and cup product from Section \ref{sec:transversality}. In particular, we will show that the dashed arrows in diagrams \eqref{eqn:differential-T1-T2} and \eqref{eqn:cup-product-T1-T2} are Poincar\'e dual to pushforward in Borel-Moore homology. Coupled with the pullback homomorphisms from \eqref{eqn:differential-T1-T2} and \eqref{eqn:cup-product-T1-T2}, this will show that the differential and cup product are Poincar\'e dual to convolution in Borel-Moore homology.

\subsection{Expressing the differentials in terms of pullback and pushforward homomorphisms}\label{subsec:precise-differential}

Now we focus on the following part of diagram \eqref{eqn:differential-T1-T2} and interpret the dashed arrow in terms of pushforward in Borel-Moore homology.
\begin{equation}\label{eqn:zoom-in-differential-T1-T2}
\begin{tikzcd}
H^p(W_{j+1,0}, W_{j+1,0} \setminus \mathcal{F}_{j,0}^{j+1,0}) \arrow{r} & H^p(W_{j+1,0}) \arrow{r} & H^{p+1}(W_{j+1}, W_{j+1,0}) \\
H^{p-\nu_j}(\mathcal{F}_{j,0}^{j+1,0}) \arrow[leftrightarrow]{u}{Thom}[swap]{\cong} \arrow[dashrightarrow]{rr} & & H^{p+1-\lambda_{j+1}}(C_{j+1}) \arrow[leftrightarrow]{u}{Thom}[swap]{\cong}
\end{tikzcd}
\end{equation}

Following the notation of \cite{AustinBraam95}, let $\tilde{\mathcal{F}}_{j,0}^{j+1,0} = \mathcal{F}_{j,0}^{j+1,0} / \mathbb{R}$ denote the quotient by the $\mathbb{R}$-action of the flow and note that this is a homotopy equivalence, so $H^*(\mathcal{F}_{j,0}^{j+1,0}) \cong H^*(\tilde{\mathcal{F}}_{j,0}^{j+1,0})$.

\begin{proposition}\label{prop:differential-pushforward}
Consider a minimally degenerate Morse function satisfying the conditions of Theorem \ref{thm:main-thm-morse} and the weak transversality conditions of Definition \ref{def:transversality-conditions}. Then the homomorphism corresponding to the dashed arrow in \eqref{eqn:zoom-in-differential-T1-T2} is Poincar\'e dual to the pushforward in Borel-Moore homology associated to the proper projection map $\tilde{\mathcal{F}}_{j,0}^{j+1,0} \rightarrow C_{j+1}$ taking a point on a flow line to the critical point which is the limit of the upwards flow.
\end{proposition}

\begin{proof}
After using the flow to define a homotopy $W_{j+1,0} \simeq W_{j+1,0} / \mathbb{R}$ of $W_{j+1,0}$ with a sphere bundle, Lemma \ref{lem:pushforward-submanifold} implies that the composition 
\begin{equation*}
H^{p-\nu_j}(\tilde{\mathcal{F}}_{j,0}^{j+1,0}) \rightarrow H^p(W_{j+1,0}, W_{j+1,0} \setminus \tilde{\mathcal{F}}_{j,0}^{j+1,0}) \rightarrow H^p(W_{j+1,0})
\end{equation*}
is Poincar\'e dual to pushforward in Borel-Moore homology for the inclusion map $\tilde{\mathcal{F}}_{j,0}^{j+1,0} \hookrightarrow W_{j+1,0}$. Lemma \ref{lem:pushforward-sphere-bundle} then shows that the composition
\begin{equation*}
H^p(W_{j+1,0}) \rightarrow H^{p+1}(W_{j+1}, W_{j+1,0}) \rightarrow H^{p+1-\lambda_{j+1}}(C_{j+1})
\end{equation*}
is Poincar\'e dual to pushforward for the projection of the sphere bundle $W_{j+1,0} \rightarrow C_{j+1}$. Therefore diagram \eqref{eqn:zoom-in-differential-T1-T2} becomes
\begin{equation*}
\begin{tikzcd}[column sep=15mm]
H^p(W_{j+1,0}, W_{j+1,0} \setminus \mathcal{F}_{j,0}^{j+1,0}) \arrow{r} & H^p(W_{j+1,0}) \arrow{r} \arrow{dr}[swap]{pushforward} & H^{p+1}(W_{j+1}, W_{j+1,0}) \\
H^{p-\nu_j}(\mathcal{F}_{j,0}^{j+1,0}) \arrow[leftrightarrow]{u}{\text{Thom}}[swap]{\cong} \arrow[dashrightarrow]{rr} \arrow{ur}[swap]{pushforward} & & H^{p+1-\lambda_{j+1}}(C_{j+1}) \arrow[leftrightarrow]{u}{\text{Thom}}[swap]{\cong}
\end{tikzcd}
\end{equation*}
and so the dashed arrow is a composition of homomorphisms which are Poincar\'e dual to pushforward in Borel-Moore homology. This composition is then the pushforward associated to the projection $\tilde{\mathcal{F}}_{j,0}^{j+1,0} \rightarrow C_{j+1}$ mapping a point on a flow line to the associated upper critical point.
\end{proof}

Therefore the spectral sequence differential between adjacent critical sets is (up to isomorphism given by the main theorem of Morse theory and the Thom isomorphism) given by a homomorphism $H^{p-\lambda_j}(C_j) \rightarrow H^{p+1-\lambda_{j+1}}(C_{j+1})$ which is the composition of pullback to the space of flow lines, cup product with the equivariant Euler class from Corollary \ref{cor:Euler-class} and pushforward to the upper critical set.

\begin{equation}\label{eqn:differential-final}
\begin{tikzcd}[column sep=6mm]
H^p(Z_j, Z_{j-1}) \arrow{dr}{restriction} \arrow[leftrightarrow]{d}{\cong}[swap]{\text{Thm.} \, \ref{thm:main-thm-morse}} \arrow{rrr}{d} & & &  H^{p+1}(Z_{j+1}, Z_j) \arrow[leftrightarrow]{d}{\cong}[swap]{\text{Thm.} \,  \ref{thm:main-thm-morse}} \\
H^p(W_j, W_{j,0})  & H^p(W_{j+1,0}, W_{j+1,0} \setminus \mathcal{F}_{j,0}^{j+1,0}) \arrow{r} & H^p(W_{j+1,0}) \arrow{r} & H^{p+1}(W_{j+1}, W_{j+1,0}) \\
 & H^{p-\nu_j}(\mathcal{F}_{j,0}^{j+1,0}) \arrow[leftrightarrow]{u}{Thom}[swap]{\cong} \arrow{rr}{pushforward} & & H^{p+1-\lambda_{j+1}}(C_{j+1}) \arrow[leftrightarrow]{u}{Thom}[swap]{\cong} \\
H^{p-\lambda_j}(C_j) \arrow[leftrightarrow]{uu}{Thom}[swap]{\cong} \arrow{r}{pullback} & H^{p-\lambda_j}(\mathcal{F}_{j,0}^{j+1,0}) \arrow{u}{\smallsmile i^* e}  & & 
\end{tikzcd}
\end{equation}

\subsection{Expressing the cup product in terms of pullback and pushforward homomorphisms}\label{subsec:precise-cup-product}

In this section we use equivariant cohomology with respect to a circle action on $W_{j+1}$ satisfying the conditions of the Atiyah-Bott Lemma \cite[Prop. 13.4]{AtiyahBott83}, namely that the action is free on $W_{j+1,0}$ and fixes the zero section $C_{j+1}$. We also assume that this action comes from a circle subgroup of a compact group acting on the total space $M$ and preserving the filtration $M_0 \subset \cdots \subset M_n = M$. These assumptions are satisfied in applications where $M$ is a symplectic manifold with a Hamiltonian group action and moment map $\mu$, and the associated filtration is determined by the gradient flow of $\| \mu \|^2$ (cf. \cite{Kirwan84}).

First recall the diagram \eqref{eqn:cup-product-T1-T2}. For the remainder of this section we focus on the following subdiagram and interpret the topological meaning of the dashed arrow as the composition of cup product followed by pushforward in Borel-Moore homology
\begin{equation}\label{eqn:zoom-in-cup-product-T1-T2}
\begin{tikzcd}[column sep=1cm]
H_K^p(W_{j+1,0}, W_{j+1,0} \setminus \mathcal{F}_{j,0}^{j+1,0}) \arrow{rr}{c_\omega \, \text{(Lem. \ref{lem:cup-product-main-thm})}} & & H_K^{p+m}(W_{j+1}, W_{j+1,0})  \\
 H_K^{p-\nu_j}(\mathcal{F}_{j,0}^{j+1,0}) \arrow[leftrightarrow]{u}{\text{Thom}}[swap]{\cong} \arrow[dashrightarrow]{rr} & & H_K^{p+m-\lambda_{j+1}}(C_{j+1}) \arrow[leftrightarrow]{u}{\text{Thom}}[swap]{\cong} 
\end{tikzcd}
\end{equation}

The first result says that, in the case where we take cup product $c_\xi$ in the top row of \eqref{eqn:zoom-in-cup-product-T1-T2} with respect to the class $\xi$ from Remark \ref{rem:Chern-class}, then the dashed arrow in \eqref{eqn:zoom-in-cup-product-T1-T2} is Poincar\'e dual to pushforward in Borel-Moore homology in the bottom row.
\begin{lemma}\label{lem:cup-product-pushforward}
With the same assumptions as Proposition \ref{prop:differential-pushforward}, let $\xi \in H^2(\mathbb{P} W_{j+1})$ denote the first Chern class of the projective bundle $\mathbb{P}W_{j+1} \rightarrow C_{j+1}$. Then the pushforward in Borel-Moore homology associated to the proper projection map $\tilde{\mathcal{F}}_{j,0}^{j+1,0} / S^1 \rightarrow C_{j+1}$ which takes an $S^1$ orbit to the critical point given by the upwards flow is Poincar\'e dual to the homomorphism corresponding to the dashed arrow in \eqref{eqn:zoom-in-cup-product-T1-T2} for the case $\omega = \xi$.
\end{lemma}

\begin{proof}
Lemma \ref{lem:pushforward-submanifold} (see also Remark \ref{rem:Chern-class}) shows that the homomorphism
\begin{equation*}
H_K^{p-\nu_j}(\mathcal{F}_{j,0}^{j+1,0}) \rightarrow H_K^p(W_{j+1,0}, W_{j+1,0} \setminus \mathcal{F}_{j,0}^{j+1,0}) \rightarrow H_K^p(W_{j+1,0})
\end{equation*}
is Poincar\'e dual to pushforward in Borel-Moore homology corresponding to inclusion of the submanifold $\mathcal{F}_{j,0}^{j+1,0} \hookrightarrow W_{j+1,0}$. Therefore it remains to push forward to the critical set $C_{j+1}$.

Lemma \ref{lem:pushforward-cup-product} shows that there are canonical splitting homomorphisms \eqref{eqn:relative-splitting-inclusion} and \eqref{eqn:relative-splitting-projection} for the short exact sequences associated to the pair $(W_{j+1}, W_{j+1,0})$ such that pushforward is given by the homomorphism \eqref{eqn:pushforward-cup-product}. 

Given the Chern class $\xi \in H^2(\mathbb{P} W_{j+1})$, the cup product map $c_\xi$ from Lemma \ref{lem:cup-product-main-thm} comes from diagram \eqref{eqn:restriction-cup-product-relative}, which can be augmented with another restriction homomorphism to the short exact sequence of the pair $(W_{j+1}, W_{j+1,0})$ in the bottom row of the following diagram

{\small \begin{equation*}
\begin{tikzcd}[column sep=0.5cm]
& & H_{K}^p(W_{j+1}, W_{j+1,0} \setminus \mathcal{F}_{j,0}^{j+1,0}) \arrow{r} \arrow{d}{\smallsmile \xi} & H_{K}^p(W_{j+1,0}, W_{j+1,0} \setminus \mathcal{F}_{j,0}^{j+1,0}) \arrow{r} \arrow{dd} & 0 \\
0 \arrow{r} & H_{K}^{p+2}(W_{j+1}, W_{j+1,0}) \arrow{r} \arrow[leftrightarrow]{dd}{=} & H_{K}^{p+2}(W_{j+1}, W_{j+1,0} \setminus \mathcal{F}_{j,0}^{j+1,0}) \\
& & H_{K}^p(W_{j+1}) \arrow{r} \arrow{d}{\smallsmile \xi} & H_{K}^p(W_{j+1,0}) \arrow{r}  & 0 \\
0 \arrow{r} & H_{K}^{p+2}(W_{j+1}, W_{j+1,0}) \arrow{r} & H_{K}^{p+2}(W_{j+1}) 
\end{tikzcd}
\end{equation*}
}

Recall from Lemma \ref{lem:splitting-compatible} that splitting homomorphisms for each row of the above diagram can be chosen to be compatible with these canonical splitting homomorphisms for the short exact sequences associated to the pair $(W_{j+1}, W_{j+1,0})$. In particular, the cup product homomorphism from Lemma \ref{lem:cup-product-main-thm} restricts to the pushforward homomorphism \eqref{eqn:pushforward-cup-product} from Lemma \ref{lem:pushforward-cup-product}. Therefore the dashed arrow in \eqref{eqn:zoom-in-cup-product-T1-T2} corresponds to the composition of pushforward homomorphisms.

\begin{equation*}
\begin{tikzcd}[column sep=1cm]
H_{K}^p(W_{j+1,0}, W_{j+1,0} \setminus \mathcal{F}_{j,0}^{j+1,0}) \arrow{r} & H_{K}^p(W_{j+1,0}) \arrow{rr}{c_\xi \, \text{(Lem. \ref{lem:pushforward-cup-product})}} \arrow{drr}[swap]{pushforward}  & & H_{K}^{p+2}(W_{j+1}, W_{j+1,0})  \\
 H_K^{p-\nu_j}(\mathcal{F}_{j,0}^{j+1,0}) \arrow[leftrightarrow]{u}{\text{Thom}}[swap]{\cong} \arrow[dashrightarrow]{rrr} \arrow{ur}[swap]{pushforward} & & & H_{K}^{p+2-\lambda_{j+1}}(C_{j+1}) \arrow[leftrightarrow]{u}{\text{Thom}}[swap]{\cong} 
\end{tikzcd}
\end{equation*}

Therefore,  in the case $\omega = \xi$, the cup product homomorphism for the triple $(W_{j+1}, W_{j+1,0}, W_{j+1,0} \setminus \mathcal{F}_{j,0}^{j+1,0})$ from Proposition \ref{prop:cup-product-main-theorem} restricts to the pushforward given by \eqref{eqn:pushforward-cup-product}.
\end{proof}

More generally, one can take cup product with $\omega = \eta \smallsmile \xi$ for arbitrary $\eta \in H_K^{m-2}(Z)$. 

\begin{proposition}\label{prop:cup-product-pushforward}
With the same assumptions as Lemma \ref{lem:cup-product-pushforward}, the cup product on the Morse complex with a class $\omega = \eta \smallsmile \xi \in H_K^*(W_{j+1})$ is induced by pullback from the lower critical set to the space of flow lines, then cup product with the Euler class from Corollary \ref{cor:Euler-class}, then cup product with the class $\eta$ followed by pushforward to the upper critical set.
\end{proposition}

Therefore the dashed arrow in diagram \eqref{eqn:cup-product-T1-T2} can be filled in as follows
\begin{equation}\label{eqn:cup-product-final}
\begin{tikzcd}[column sep=0.6cm]
H_K^p(M_j, M_{j-1}) \arrow[leftrightarrow]{d}{\cong}[swap]{\text{Thm.} \, \ref{thm:main-thm-morse}} \arrow{dr}{restriction} & & & H_K^{p+m}(M_{j+1}, M_j) \arrow[leftrightarrow]{d}{\cong}[swap]{\text{Thm.} \, \ref{thm:main-thm-morse}} \\
H_K^p(W_j, W_{j,0})  & H_K^p(W_{j+1,0}, W_{j+1,0} \setminus \mathcal{F}_{j,0}^{j+1,0}) \arrow{rr}{c(\eta \smallsmile \xi) \, \text{(Lem. \ref{lem:cup-product-main-thm})}} & & H_K^{p+m}(W_{j+1}, W_{j+1,0})  \\
 & H_K^{p-\nu_j}(\mathcal{F}_{j,0}^{j+1,0}) \arrow[leftrightarrow]{u}{\text{Thom}}[swap]{\cong} \arrow{r}{\smallsmile \eta} & H_K^{p+m-\nu_j}(\mathcal{F}_{j,0}^{j+1,0}) \arrow{r}{push}[swap]{forward} & H_K^{p+m-\lambda_{j+1}}(C_{j+1}) \arrow[leftrightarrow]{u}{\text{Thom}}[swap]{\cong} \\
H_K^{p-\lambda_j}(C_j) \arrow[leftrightarrow]{uu}{\text{Thom}}[swap]{\cong} \arrow{r}{pullback} & H_K^{p-\lambda_j}(\mathcal{F}_{j,0}^{j+1,0}) \arrow{u}{\smallsmile e}  & & 
\end{tikzcd}
\end{equation}

From the above, we see that the class $\xi \in H_K^2(W_{j+1})$ can be used to turn the Morse-theoretic cup product from Proposition \ref{prop:cup-product-main-theorem} into pushforward (Lemma \ref{lem:cup-product-pushforward}) or, more generally, cup product followed by pushforward (Proposition \ref{prop:cup-product-pushforward}). Two questions now remain. Firstly, what happens when we restrict to the subset of representations of a quiver which satisfy a given set of relations? Secondly, what is the topological meaning of the induced homomorphism $H_K^p(M_j, M_{j-1}) \rightarrow H_K^{p+m}(M_{j+1}, M_j)$ in the diagram \eqref{eqn:cup-product-final}? 

These questions are answered in the next section, where we show that cup product with the right choice of class $\eta$ in \eqref{eqn:cup-product-final} pulls back to a homomorphism Poincar\'e dual to convolution in Borel-Moore homology (Theorem \ref{thm:cup-product-induces-convolution}) and that we can interpret this in the Morse spectral sequence by pulling back from $M$ to a product $M \times \Gr$ and applying the cup product in the Morse spectral sequence of this product space (diagram \eqref{eqn:full-convolution-cup-product-diagram}).

\section{Convolution in Borel-Moore homology is determined by the cup product}\label{sec:convolution-cup-product}

Now we return to the case of representations of a quiver with complete quadratic relations (Definition \ref{def:complete-relations}) that are fully restricted (Definition \ref{def:fully-restricted}) from the relations for a Nakajima quiver. If there are loops then we also impose the condition of Definition \ref{def:loop-neg-slice}. As a consequence of these conditions, the subset of stable representations is smooth (Corollary \ref{cor:smoothness-restricts}) and the space of flow lines between adjacent critical sets is smooth (Corollary \ref{cor:neg-slice-quiver-smooth}). These conditions are not too restrictive, as this setup includes Nakajima quivers without loops (Example \ref{ex:nakajima-complete}), the ADHM quiver (Remark \ref{rem:ADHM-loops}) and handsaw quivers (Examples \ref{ex:handsaw-complete} and \ref{ex:handsaw-reduction}), but is not limited to these examples. 

The goal of this section is to prove Theorem \ref{thm:cup-product-induces-convolution}, which uses the cup product on the Morse complex to show that the convolution homomorphism of \cite[Sec. 2(i)]{Nakajima97} and \cite[(2.7.14)]{ChrissGinzburg97}  in Borel-Moore homology can be constructed by first applying cup product on the Morse complex of the product of the ambient smooth space $M = \Rep(Q, {\bf v})$ with a Grassmannian and then restricting to the singular subset $\nu^{-1}(0)$. The key is to take cup product with the correct class, which in this case is the Thom class of a certain submanifold of the space of flow lines in $\Rep(Q, {\bf v})$.

For the remainder of this section, fix two adjacent critical sets $C_u$ (upper) and $C_\ell$ (lower) corresponding to the respective strata $M_{j+1} \setminus M_j$ and $M_j \setminus M_{j-1}$, with associated dimension vectors ${\bf v_u}$ and ${\bf v_\ell} = {\bf v_u} + {\bf e_k}$ and use the homotopy equivalence of Corollary \ref{cor:reduce-to-zero} so that $C_u / K_{\bf v_u} \cong \mathcal{M}(Q, {\bf v_u}) \times \{0\}$ and $C_\ell / K_{\bf v_\ell} \cong \mathcal{M}(Q, {\bf v_\ell}) \times \{0\}$.

In constructing the cup product on the Morse complex, the results of Section \ref{subsec:cup-product-T1-T2} show that it is sufficient to consider the case where $C_\ell$ is the minimum of $\| \mu - \alpha \|^2$ on $\Rep(Q, {\bf v})$, so that ${\bf v} = {\bf v_\ell}$. If ${\bf v} > {\bf v_\ell}$ so that $C_\ell$ is non-minimal then the construction follows diagram \eqref{eqn:cup-product-T1-T2}, which includes the additional step in which we take cup product with the Euler class of Corollary \ref{cor:Euler-class}. Therefore, \emph{for the remainder of this section we fix ${\bf v} = {\bf v_\ell} = {\bf v_u} + {\bf e_k}$}.


First note that the lower critical set $C_\ell$ determines a reduction of structure group from $K_{\bf v}$ to $K_{\bf v_\ell} \times K_{{\bf v} - {\bf v}_\ell}$ with respect to which (modulo $K_{\bf v}$) a critical representation can be written as $x_\ell = x_\ell^{(1)} + x_\ell^{(2)}$. Corollary \ref{cor:reduce-to-zero} shows that there is a $K_{\bf v_\ell}$-equivariant homotopy equivalence to the subset of representations for which $x_\ell^{(2)} = 0$.

Now recall that the upper critical set determines a reduction of structure group (Corollary \ref{cor:critical-level-isomorphism}) so that 
\begin{equation*}
H_{K_{\bf v_\ell}}^*(C_u) \cong H_{K_{\bf v_u}}^*(C_u) \otimes H^*(BK_{\bf e_k})
\end{equation*}
with $C_u / K_{\bf v_u} \cong \mathcal{M}(Q, {\bf v_u}) \times \{0\}$ and $K_{\bf e_k} \cong \U(1)$. The circle subgroup $\{ \id \} \times \mathsf{U}(1) \subset K_{\bf v_u} \times \mathsf{U}(1)$ then acts freely on $\mathcal{F}_{\ell,0}^{u,0}$. 

Now that we have fixed the lower and upper critical sets, then we drop this notation from the space of flow lines and use $\mathcal{F} = \mathcal{F}_{\ell,0}^{u,0}$.

\subsection{Construction of the Thom class}

Given this choice of upper critical set $C_u$, recall the definition of the negative slice bundle $\pi_u : S_u^- \rightarrow C_u$ \eqref{eqn:neg-slice-bundle}, and the unstable bundle \eqref{eqn:unstable-bundle-def}, which we denote by $W_u^-$. Theorem \ref{thm:slice-homeo} determines a homeomorphism of pairs
\begin{equation*}
(S_u^-, S_{u,0}^-) \stackrel{\cong}{\leftrightarrow} (W_u^-, W_{u,0}^-) .
\end{equation*}
The space of flow lines $\mathcal{F} \subset W_{u,0}^-$ is then mapped into the negative slice bundle, and we abuse the notation by also using $\mathcal{F} \subset S_{u,0}^-$ to denote the image. The above homeomorphism then extends to a homeomorphism of triples
\begin{equation*}
(S_u^-, S_{u,0}^-, S_{u,0}^- \setminus \mathcal{F}) \stackrel{\cong}{\leftrightarrow} (W_u^-, W_{u,0}^-, W_{u,0}^- \setminus \mathcal{F}) .
\end{equation*}

The negative slice bundle is a linearisation of the unstable set, which simplifies the deformation theory associated to the subsets $\mathcal{B}$, $N$ and $T$ defined below. We can construct normal bundles inside the negative slice and then use the above homeomorphism to map these into the unstable set $W_u^-$. The homeomorphism guarantees that these images will also have well-defined normal bundles.

Define 
\begin{equation*}
\mathcal{B} = \mathcal{F} \cap \nu^{-1}(0) \subset S_u^-.
\end{equation*}
Since $\nu^{-1}(0)$ is closed and preserved by the flow, then any flow line in $\mathcal{B}$ has both upper and lower limits in $\nu^{-1}(0)$. Conversely, if the lower limit is in $\nu^{-1}(0)$ then Lemma \ref{lem:isomorphic-to-lower-limit} shows that all points in the flow line are isomorphic to the lower limit and hence in $\nu^{-1}(0)$. The result of \cite[Thm. 4.35]{Wilkin17} shows that $\mathcal{B} / K_{{\bf v}_\ell}$ is the Hecke correspondence for the quiver varieties $\mathcal{M}(Q, {\bf v}_\ell, \mathcal{R})$ and $\mathcal{M}(Q, {\bf v}_u, \mathcal{R})$.

Now define 
\begin{equation*}
T = \pi_u^{-1}(\nu^{-1}(0) \cap C_u) \cap \mathcal{F} \subset S_u^- 
\end{equation*}
and
\begin{equation*}
N = \{ x \in \mathcal{F} \, \mid \, \nu(x) = \nu(\pi_u(x)) \} \subset S_u^- .
\end{equation*}
Note that
\begin{equation*}
T \cap N = \{ x \in \mathcal{F} \, \mid \, \nu(x) = \nu(\pi_u(x)) = 0 \} = \mathcal{F} \cap \nu^{-1}(0) = \mathcal{B} .
\end{equation*}

\begin{lemma}
$\mathcal{F}$ and $T$ are manifolds, and $N \subset \mathcal{F}$ and $\mathcal{B} \subset T$ are submanifolds with the same codimension $d$.
\end{lemma}

\begin{proof}
Recall that the fibres of $\mathcal{F}$ are determined by the cokernel of \eqref{eqn:flow-lines-cokernel}. The upper critical set $C_{\bf v_u}$ is preserved by the action of $K_{\bf v_u}$ and (within the stratum $\Rep(Q, {\bf v_u})^{st}$) the normal bundle to the critical set has fibres given by the infinitesimal action of $i \mathfrak{k}_{\bf v_u} \subset \mathfrak{g}_{\bf v_u}$. The deformation complex describing the tangent space of $\mathcal{F}$ is then given by
\begin{equation*}
\begin{tikzcd}
i \mathfrak{k}_{\bf v_u} \arrow{r}{\rho_{(x_u+y)}^\mathbb{C}} & \Hom^1(Q, {\bf e_k}, {\bf v_u}) \oplus \Hom^1(Q, {\bf v_u}, {\bf v_u}) .
\end{tikzcd}
\end{equation*}
Stability of $x_u$ implies that (modulo the scalar multiples of the identity in $\mathfrak{k}_{\bf v_u}$) this homomorphism is injective. A point of $N$ is given by a pair $x_u+ y \in S_{u}^-$, where $x_u \in C_u$ and $\nu(x_u+y) = \nu(x_u)$. Therefore Lemma \ref{lem:linearise-slice} shows that $x_u + y \in N$ if and only if $d \nu_{x_u}(y) = 0$ and so the derivative of this condition defines the tangent bundle of $N$, which is given by the middle cohomology of the following complex
\begin{equation*}
\begin{tikzcd}[column sep = 2cm]
i \mathfrak{k}_{\bf v_u}  \arrow{r}{\rho_{(x_u+y)}^\mathbb{C}} & \Hom^1(Q, {\bf v_u}, {\bf v_u}) \oplus \Hom^1(Q, {\bf e_k}, {\bf v_u})   \arrow{r}{d \nu_{x_u} + d\nu_y} & \Rel_0(Q, {\bf e_k}, {\bf v_u}, \mathcal{R}) ,
\end{tikzcd}
\end{equation*}
where (using the notation of \eqref{eqn:derivative-neg-slice-quiver} in Example \ref{ex:neg-slice-restricted})
\begin{equation}\label{eqn:sum-relations-N}
(d \nu_{x_u} + d\nu_y)(\delta x, \delta y) = d \nu_{x_u}(\delta y) + d\nu_y(\delta x) .
\end{equation}
Now $x_u + y \in \Rep(Q, {\bf v_\ell})$ is stable by Lemma \ref{lem:isomorphic-to-lower-limit} and therefore Lemma \ref{lem:stable-implies-adjoint-injective} shows that the adjoint of $d\nu_{(x_u+y)}$ is injective. Now note that the adjoint of the homomorphism above is the restriction of $d\nu_{(x_u+y)}^*$ to $\Rel(Q, {\bf e_k}, {\bf v_u}, \mathcal{R})$ and therefore must also be injective. Therefore $N \subset \mathcal{F}$ is smooth and the normal bundle is determined by the image of the adjoint of \eqref{eqn:sum-relations-N} in the above deformation complex.

A similar construction works for $\mathcal{B}$. For hyperk\"ahler quivers and handsaw quivers, Nakajima has already shown that $\mathcal{B}$ is smooth (cf. \cite[Sec. 5]{Nakajima98} and \cite[Sec. 5]{Nakajima12}) using deformation theory inside $\mathcal{M}(Q, {\bf v_u}) \times \mathcal{M}(Q, {\bf v_u} + {\bf e_k})$, however here the goal is to relate the normal bundles of $\mathcal{B} \subset N$ and $T \subset \mathcal{F}$ and so we instead use the deformation theory of $\mathcal{B} \subset N \subset \mathcal{F}$. The associated deformation complex is
\begin{equation*}
\begin{tikzcd}[column sep=1.2cm]
i \mathfrak{k}_{\bf v_u}  \arrow{r}{\rho_{(x+y)}^\mathbb{C}} & \Hom^1(Q, {\bf v_u}, {\bf v_u}) \oplus \Hom^1(Q, {\bf e_k}, {\bf v_u}) \arrow{r}{d \nu_{(x_u+y)}} & \Rel_0(Q, {\bf v_u}, {\bf v_u}, \mathcal{R}) \oplus \Rel(Q, {\bf e_k}, {\bf v_u}, \mathcal{R}) 
\end{tikzcd}
\end{equation*}
where the extra condition $x_u \in \nu^{-1}(0)$ determines an extra term in the derivative compared to the above \eqref{eqn:sum-relations-N}
\begin{equation*}
d \nu_{(x_u+y)}(\delta x, \delta y) = \left( d \nu_{x_u}(\delta x), d \nu_{x_u}(\delta y) + d\nu_y(\delta x) \right) .
\end{equation*}
Again, since $x_u+y$ is stable, then the adjoint of $d \nu_{(x_u+y)}$ is injective (modulo scalar multiples of the identity) by Lemma \ref{lem:stable-implies-adjoint-injective} and so the normal bundle of $\mathcal{B} \subset \mathcal{F}$ is determined by the image of $d \nu_{(x_u+y)}^*$. 

In order to determine the normal bundle of $\mathcal{B} \subset T$, the definition $T = \pi_u^{-1}(\nu^{-1}(0) \cap C_u) \cap \mathcal{F}$ determines an extra \emph{a priori} condition $d \nu_{x_u}(\delta x_u) = 0$ on the tangent space which leads to the following deformation complex
\begin{equation*}
\begin{tikzcd}[column sep=2cm]
i \mathfrak{k}_{\bf v_u}  \arrow{r}{\rho_{(x_u+y)}^\mathbb{C}} & \Hom^1(Q, {\bf v_u}, {\bf v_u}) \oplus \Hom^1(Q, {\bf e_k}, {\bf v_u}) \arrow{r}{d \nu_{x_u} + d\nu_y} & \Rel(Q, {\bf e_k}, {\bf v_u}, \mathcal{R}) .
\end{tikzcd}
\end{equation*}

For the same reason as before, this homomorphism is surjective and the normal bundle of $\mathcal{B} \subset T$ has fibres given by the image of the adjoint $(d \nu_{x_u} + d\nu_y)^*$, which is injective and therefore has the same dimension as the image of the adjoint of \eqref{eqn:sum-relations-N}, which corresponds to the fibres of the normal bundle of $N \subset \mathcal{F}$
\end{proof}

Let $\tau_\mathcal{B} \in H^d(T, T \setminus \mathcal{B})$ and $\tau_N \in H^d(\mathcal{F}, \mathcal{F} \setminus N)$ be the Thom classes associated to the respective inclusions $\mathcal{B} \hookrightarrow T$ and $N \hookrightarrow \mathcal{F}$. Denote the inclusion of pairs by $i_1 : (T, T \setminus \mathcal{B}) \hookrightarrow (\mathcal{F}, \mathcal{F} \setminus N)$.

\begin{corollary}\label{cor:thom-compatible}
The Thom classes satisfy $\tau_\mathcal{B} = i_1^* \tau_N$.
\end{corollary}

\begin{proof}
This follows from the previous proof by noting that the normal bundle of $\mathcal{B} \subset T$ is the image of
\begin{equation*}
\begin{tikzcd}[column sep = 2cm]
\Rel(Q, {\bf e_k}, {\bf v_u}, \mathcal{R}) \arrow{r}{(d \nu_{x_u} + d\nu_y)^*} & \Hom^1(Q, {\bf v_u}, {\bf v_u}) \oplus \Hom^1(Q, {\bf e_k}, {\bf v_u}) \hookrightarrow \Hom^1(Q, {\bf v}, {\bf v}), 
\end{tikzcd}
\end{equation*}
where $\Hom^1(Q, {\bf v}, {\bf v})$ is a globally defined trivial bundle over $\mathcal{B} \subset T \subset \mathcal{F}$. Therefore the normal bundle of $\mathcal{B} \subset T$ is the pullback of the normal bundle of $N \subset \mathcal{F}$ via the inclusion $\mathcal{B} \hookrightarrow N$, which implies the relation $\tau_\mathcal{B} = i^* \tau_N$ on the Thom classes.
\end{proof}


Now let $\Gr({\bf v_u}, {\bf v_\ell})$ denote the Grassmannian of ${\bf v_u}$ planes in ${\bf v_\ell}$, or equivalently injective homomorphisms $\Vect(Q, {\bf v_u}) \hookrightarrow \Vect(Q, {\bf v_\ell})$ modulo $G_{\bf v_u}$. Consider the product $\Rep(Q, {\bf v}_\ell)^{st} \times \Gr({\bf v}_u, {\bf v}_\ell)$ with the induced group action of $K_{\bf v_\ell}$ on $\Gr({\bf v}_u, {\bf v}_\ell)$. Using the Hermitian structure, a subspace in $\Gr({\bf v}_u, {\bf v}_\ell)$ corresponds to an orthogonal projection $\Pi : \Vect(Q, {\bf v}_\ell) \rightarrow \Vect(Q, {\bf v}_u)$. Therefore there is a continuous projection
\begin{align*}
p_u : \Rep(Q, {\bf v}_\ell)^{st} \times  \Gr({\bf v}_u, {\bf v}_\ell) & \rightarrow \Rep(Q, {\bf v}_u) \\
(x, \Pi) & \mapsto \Pi x \Pi .
\end{align*}
Define the open subset $U \subset \Rep(Q, {\bf v}_\ell)^{st} \times  \Gr({\bf v}_u, {\bf v}_\ell)$ as the preimage of $\Rep(Q, {\bf v}_u)^{st}$ by this projection. There is a continuous map $U \rightarrow \Rep(Q, {\bf v}_\ell)^{st}$ given by projection onto the first factor. Note that each $x \in S_{u,0}^-$ comes equipped with a canonical element of $\Gr({\bf v}_u, {\bf v}_\ell)$ determined by reduction of structure group at the upper critical point, and therefore the inclusion $\mathcal{F} \hookrightarrow S_{u,0} \hookrightarrow \Rep(Q, {\bf v}_\ell)^{st}$ lifts as follows
\begin{equation}\label{eqn:canonical-lift}
\begin{tikzcd}
 & & U \arrow{d} \arrow[hookrightarrow]{r} & \Rep(Q, {\bf v_\ell}) \times \Gr({\bf v_u}, {\bf v_\ell}) \\
 \mathcal{F} \arrow[hookrightarrow]{r} & S_{u,0}^- \arrow[hookrightarrow]{r} \arrow[dashrightarrow]{ur} & \Rep(Q, {\bf v}_\ell)^{st} 
\end{tikzcd}
\end{equation}

There is also a projection onto $\Hom^1(Q, {\bf e}_k, {\bf v}_u)$ given by $\Pi x \Pi^\perp$. For each $(x, \Pi) \in U$, denote the representations associated to these projections by $x_u := \Pi x \Pi \in \Rep(Q, {\bf v}_u)^{st}$ and $y := \Pi x \Pi^\perp \in \Hom^1(Q, {\bf e}_k, {\bf v}_u)$. 

The projection $\Pi$ also determines a choice of subgroup $K_{\bf v_u} \subset K_{\bf v_\ell}$ and, with respect to the action of $G_{\bf v_u}$ on $x$, there is an orthogonal decomposition
\begin{equation*}
\Hom^1(Q, {\bf e}_k, {\bf v_u}) \cong \im \rho_{x_u}^\mathbb{C} \oplus \ker (\rho_{x_u}^\mathbb{C})^* .
\end{equation*}
Given $y = \Pi x \Pi^\perp$ as above, let $y_h$ denote the component in $\ker (\rho_{x_u}^\mathbb{C})^*$ with respect to the above decomposition. Stability of $x + y$ then implies that $y_h \neq 0$, since otherwise one can construct an isomorphism between $x$ and a representation with zero component in $\Hom^1(Q, {\bf e}_k, {\bf v_u})$, which is then clearly unstable. Now define 
\begin{equation}\label{eqn:tilde-N-def}
\tilde{N} := \left\{ (x, \Pi) \in U \, \mid \, \nu(x_u) = \nu(x_u + y_h) \right\} .
\end{equation}

\begin{lemma}\label{lem:pullback-Thom-class}
$\tilde{N}$ is a submanifold of $U$. With respect to the canonical lift $\mathcal{F} \hookrightarrow U$ from \eqref{eqn:canonical-lift}, the intersection $\tilde{N} \cap \mathcal{F}$ is transverse in $U$ and the inclusions of pairs
\begin{equation*}
(T, T \setminus \mathcal{B}) \stackrel{i_1}{\hookrightarrow} (\mathcal{F}, \mathcal{F} \setminus N) \stackrel{i_2}{\hookrightarrow} (U, U \setminus \tilde{N})
\end{equation*}
satisfy
\begin{equation}\label{eqn:Thom-class-pullback}
\tau_\mathcal{B} = (i_2 \circ i_1)^* \tau_{\tilde{N}} \quad \text{and} \quad \tau_N = i_2^* \tau_{\tilde{N}} ,
\end{equation}
where $\tau_{\tilde{N}}$ is the Thom class of $\tilde{N} \hookrightarrow U$.
\end{lemma}

\begin{proof}

The derivative of the condition $\nu(x_u) = \nu(x_u + y_h)$ is
\begin{equation*}
\begin{tikzcd}[column sep = 2cm]
T_{(x, \Pi)} \left( \Rep(Q, {\bf v_\ell})^{st} \times \Gr({\bf v_u}, {\bf v_\ell}) \right) \arrow{r}{d \nu_{x_u} + d\nu_y} & \Rel(Q, {\bf e_k}, {\bf v_u}, \mathcal{R})
\end{tikzcd}
\end{equation*}
Again, stability of $x_u$ implies that this is surjective, and therefore the adjoint is injective. In particular we see that the normal bundle of $\tilde{N} \subset U$ has fibres given by the image of the adjoint and that the image of these fibres in $T_{(x, \Pi)} \left( \Rep(Q, {\bf v_\ell})^{st} \times \Gr({\bf v_u}, {\bf v_\ell}) \right)$ restricts to the fibres of the normal bundle of $N \subset \mathcal{F}$ defined in the previous proof. 


Therefore the same argument as above shows that the normal bundle of $N \subset \mathcal{F}$ is the pullback of the normal bundle of $\tilde{N} \subset U$ via the inclusion $N \hookrightarrow \tilde{N}$, which implies the relation \eqref{eqn:Thom-class-pullback} on the Thom classes.
\end{proof}

\subsection{Global generators for the Thom class}

Now consider the image of the Thom class $\tau_N \in H_K^d(\mathcal{F}, \mathcal{F} \setminus N)$ under the homomorphism $H_K^d(\mathcal{F}, \mathcal{F} \setminus N) \rightarrow H_K^d(\mathcal{F})$. In the sequel we will denote this image by $\tau_N'$ and the goal of this section is to prove Corollary \ref{cor:global-generation-Thom}, which shows that $\tau_N'$ is generated by a global class in $H_{K_{\bf v_\ell}}^*(\Rep(Q, {\bf v}_\ell) \times  \Gr({\bf v}_u, {\bf v}_\ell))$.

Now we show that the inclusion $U \hookrightarrow \Rep(Q, {\bf v}_\ell)^{st} \times  \Gr({\bf v}_u, {\bf v}_\ell)$ induces a surjection in cohomology. Coupled with the well-known surjection
\begin{equation*}
H^*(BK_{\bf v_\ell} \times BK_{\bf v_\ell}) \twoheadrightarrow H_{K_{\bf v_\ell}}^*(\Rep(Q, {\bf v}_\ell)^{st} \times  \Gr({\bf v}_u, {\bf v}_\ell))
\end{equation*}
(using the method of Kirwan \cite{Kirwan84}), this implies that the image of the Thom class $\tau_{\tilde{N}}$ under the inclusion $H^*(U, U \setminus \tilde{N}) \rightarrow H^*(U)$ is generated by a class in $H^*(BK_{\bf v_\ell} \times BK_{\bf v_\ell})$.

\begin{lemma}
The inclusion $U \hookrightarrow \Rep(Q, {\bf v}_\ell) \times \Gr({\bf v}_u, {\bf v}_\ell)$ induces a surjection in $K_{\bf v_\ell}$-equivariant cohomology.
\end{lemma}

\begin{proof}

Choose a fixed $\Pi \in \Gr({\bf v_u}, {\bf v_\ell})$. We have
\begin{equation*}
\Gr({\bf v_u}, {\bf v_\ell}) \cong K_{\bf v_\ell} / \left( K_{\bf v_u} \times \U(1) \right) 
\end{equation*}
and so
\begin{equation*}
H_{K_{\bf v_\ell}}^* \left( \Rep(Q, {\bf v}_\ell) \times \Gr({\bf v}_u, {\bf v}_\ell) \right) \cong H_{K_{\bf v_u} \times \U(1)}^*(\Rep(Q, {\bf v_\ell})) ,
\end{equation*}
where $K_{\bf v_u} \times \U(1)$ acts via the choice of $\Pi \in \Gr({\bf v_u}, {\bf v_\ell})$. Equivalently, there is a canonical decomposition 
\begin{equation*}
\Rep(Q, {\bf v_\ell}) \cong \Rep(Q, {\bf v_u}) \oplus \Hom^1(Q, {\bf e_k}, {\bf v_u}) \oplus \Hom^1(Q, {\bf v_u}, {\bf e_k}) \oplus \Hom^1(Q, {\bf e_k}, {\bf e_k}) .
\end{equation*}
and $K_{\bf v_u} \times \U(1)$ acts via the obvious induced action. Kirwan's method \cite{Kirwan84} on $\Rep(Q, {\bf v_u})$ then shows that there is a surjection
\begin{multline*}
H_{K_{\bf v_u} \times \U(1)}^*(\Rep(Q, {\bf v_\ell})) \\
 \rightarrow H_{K_{\bf v_u} \times \U(1)}^*\left( \Rep(Q, {\bf v_u})^{\alpha-st} \oplus \Hom^1(Q, {\bf e_k}, {\bf v_u}) \oplus \Hom^1(Q, {\bf v_u}, {\bf e_k}) \oplus \Hom^1(Q, {\bf e_k}, {\bf e_k}) \right) .
\end{multline*}

It remains to restrict to the subset where $x_\ell \in \Rep(Q. {\bf v_\ell})$ is stable. We denote this by $U_\Pi$ with respect to the above choice of $\Pi \in \Gr({\bf v_u}, {\bf v_\ell})$, and note that $U \cong K_{\bf v_\ell} \times_{K_{\bf v_u} \times \U(1)} U_\Pi$. When the component $x_u \in \Rep(Q, {\bf v_u})$ is stable, then $x_\ell$ is stable if and only if the component in $\ker (\rho_{x_u}^\mathbb{C})^* \subset \Hom^1(Q, {\bf e_k}, {\bf v_u})$ is nonzero. Therefore the above projection is homotopy equivalent to a vector bundle over $\Rep(Q, {\bf v_u})^{\alpha-st}$ with fibre $\ker (\rho_{x_u}^\mathbb{C})^*$ via a homotopy equivalence that preserves the subset where $x_\ell$ is stable, and the subgroup $\{ \id\} \times \U(1) \subset K_{\bf v_u} \times \U(1)$ fixes the base and acts freely on the nonzero fibres. Therefore the Atiyah-Bott lemma implies that 
\begin{multline*}
H_{K_{\bf v_u} \times \U(1)}^*\left( \Rep(Q, {\bf v_u})^{\alpha-st} \oplus \Hom^1(Q, {\bf e_k}, {\bf v_u}) \oplus \Hom^1(Q, {\bf v_u}, {\bf e_k}) \oplus \Hom^1(Q, {\bf e_k}, {\bf e_k}) \right) \\
\rightarrow H_{K_{\bf v_u} \times \U(1)}^*(U_\Pi)
\end{multline*}
is surjective. The result then follows from the isomorphism
\begin{equation*}
H_{K_{\bf v_u} \times \U(1)}^*(U_\Pi) \cong H_{K_{\bf v_\ell}}^*(U) . \qedhere
\end{equation*}
\end{proof}


\begin{corollary}\label{cor:global-generation-Thom}
There is a class $\tilde{\tau} \in H_K^*(\Rep(Q, {\bf v_\ell}) \times \Gr({\bf v_u}, {\bf v_\ell}))$ such that pullback by the inclusion $i : \mathcal{F} \hookrightarrow \Rep(Q, {\bf v_\ell}) \times \Gr({\bf v_u}, {\bf v_\ell})$ from \eqref{eqn:canonical-lift} satisfies $i^* \tilde{\tau} = \tau_N'$.
\end{corollary}

\begin{proof}
The previous lemma shows that the image of the Thom class $\tau_{\tilde{N}} \in H_{K_{\bf v_\ell}}^*(U, U \setminus \tilde{N}) \rightarrow H_{K_{\bf v_\ell}}^*(U)$ is in the image of the pullback by inclusion $U \hookrightarrow \Rep(Q, {\bf v_\ell}) \times \Gr({\bf v_u}, {\bf v_\ell})$. Together with Lemma \ref{lem:pullback-Thom-class}, the commutative diagram

\begin{equation*}
\begin{tikzcd} 
& H_{K_{\bf v_\ell}}^*(\Rep(Q, {\bf v_\ell}) \times \Gr({\bf v_u}, {\bf v_\ell})) \arrow{d} \\
H_{K_{\bf v_\ell}}^*(U, U \setminus \tilde{N}) \arrow{r} \arrow{d} & H_{K_{\bf v_\ell}}^*(U) \arrow{d} \\
H_{K_{\bf v_\ell}}^*(\mathcal{F}, \mathcal{F} \setminus N) \arrow{r} & H_{K_{\bf v_\ell}}^*(\mathcal{F})
\end{tikzcd}
\end{equation*}
shows that $\tau_N'$ is in the image of the pullback homomorphism $H_{K_{\bf v_\ell}}^*(\Rep(Q, {\bf v_\ell}) \times \Gr({\bf v_u}, {\bf v_\ell})) \rightarrow H_{K_{\bf v_\ell}}^*(\mathcal{F})$.
\end{proof}

The following result is used in the proof of Lemma \ref{lem:global-generation-pushforward}.

\begin{corollary}\label{cor:U-surjectivity}
Inclusion $U \hookrightarrow \Rep(Q, {\bf v}_\ell)^{\alpha-st} \times \Gr({\bf v}_u, {\bf v}_\ell)$ induces a surjection in $K_{\bf v_\ell}$-equivariant cohomology.
\end{corollary}

\begin{proof}
There is a sequence of inclusions
\begin{equation*}
U \hookrightarrow \Rep(Q, {\bf v}_\ell)^{\alpha-st} \times \Gr({\bf v}_u, {\bf v}_\ell) \hookrightarrow \Rep(Q, {\bf v}_\ell) \times \Gr({\bf v}_u, {\bf v}_\ell)
\end{equation*}
together with the previous lemma shows that the induced homomorphism
\begin{equation*}
H_{K_{\bf v_\ell}}^* \left( \Rep(Q, {\bf v}_\ell) \times \Gr({\bf v}_u, {\bf v}_\ell) \right) \rightarrow H_{K_{\bf v_\ell}}^*(U)
\end{equation*}
is surjective. Therefore 
\begin{equation*}
H_{K_{\bf v_\ell}}^* \left( \Rep(Q, {\bf v}_\ell)^{\alpha-st} \times \Gr({\bf v}_u, {\bf v}_\ell) \right) \rightarrow H_{K_{\bf v_\ell}}^*(U)
\end{equation*}
is also surjective.
\end{proof}

\subsection{Global generation of the first Chern class of the negative slice bundle}

Now we show that the Chern class $\xi \in H_{\U(1)}^2(W_{u,0}^-) \cong H_{\U(1)}^2(S_{u,0}^-) \cong H^2(\mathbb{P} S_u^-)$ used in Lemma \ref{lem:cup-product-pushforward} is generated by a class in $H_K^2(\Rep(Q, {\bf v}) \times \Gr({\bf v}_u, {\bf v}_\ell)$.

\begin{lemma}\label{lem:global-generation-pushforward}
There exists $\tilde{\xi} \in H_K^2(\Rep(Q, {\bf v}) \times \Gr({\bf v}_u, {\bf v}_\ell))$ such that $i^* \tilde{\xi} = \xi$, where $i : S_{u,0}^- \hookrightarrow \Rep(Q, {\bf v}) \times \Gr({\bf v}_u, {\bf v}_\ell))$ is the inclusion from \eqref{eqn:canonical-lift} and $\xi$ is the Chern class from Remark \ref{rem:Chern-class} and Lemma \ref{lem:cup-product-pushforward}. 
\end{lemma}

\begin{proof}
Recall from the construction of the submanifold $\tilde{N} \subset U$ from \eqref{eqn:tilde-N-def} that given $(x, \Pi) \in U$ the projection $\Pi$ induces representations
\begin{equation*}
x_u := \Pi x \Pi \in \Rep(Q, {\bf v_u})^{\alpha-st} \quad \text{and} \quad y = \Pi x \Pi^\perp \in \Hom^1(Q, {\bf e_k}, {\bf v_u}) .
\end{equation*} 
The decomposition $\Hom^1(Q, {\bf e_k}, {\bf v_u}) \cong \ker (\rho_{x_u}^\mathbb{C})^* \oplus \im \rho_x^\mathbb{C}$ then determines a well defined vector bundle $\tilde{S}$ over $U$ with fibres $\ker (\rho_{x_u}^\mathbb{C})^* \subset \Hom^1(Q, {\bf e_k}, {\bf v_u}) \subset \Hom^1(Q, {\bf v}, {\bf v})$. Moreover, $y \in \Hom^1(Q, {\bf e_k}, {\bf v_u})$ determines a canonical nonzero $y_h \in \ker (\rho_{x_u}^\mathbb{C})^*$, which then determines a line in $\ker (\rho_{x_u}^\mathbb{C})^*$ and hence a line subbundle $\tilde{L} \subset \tilde{S}$.

On the image of the negative slice $S_{u,0}^- \hookrightarrow U$, this line bundle $\tilde{L}$ pulls back to the line bundle $L$ which is the tautological bundle of $\mathbb{P}S_u^-$ and hence the first Chern class of $\tilde{L}$ restricts to the first Chern class of $L$, which is the class $\xi$ used in Lemma \ref{lem:cup-product-pushforward}.

Coupled with the surjectivity from Corollary \ref{cor:U-surjectivity}, this proves the existence of $\tilde{\xi} \in H_K^2(\Rep(Q, {\bf v}) \times \Gr({\bf v}_u, {\bf v}_\ell))$ such that $i^* \tilde{\xi} = \xi$.
\end{proof}

\subsection{Cup product with the Thom class on the Morse complex}\label{subsec:cup-product-Thom-class}

Now we are in a position to describe the cup product with the Thom class on the Morse complex and then use the above results to give a Morse-theoretic construction of convolution in Borel-Moore homology. To simplify the notation in the following, we use $\Gr$ in place of $\Gr({\bf v}_u, {\bf v}_\ell)$.

Lemma \ref{lem:pullback-Thom-class} shows that the Thom class $\tau_N$ is in the image of $\tau_{\tilde{N}} \in H_K^d(U, U \setminus \tilde{N})$. The image $\tau_{\tilde{N}}' \in H_K^d(U)$ is in the image of the pullback homomorphism from $H_{K_{\bf v_\ell}}^*(\Rep(Q, {\bf v_\ell}) \times \Gr({\bf v_u}, {\bf v_\ell})) \rightarrow H_{K_{\bf v_\ell}}^*(U)$ by Corollary \ref{cor:U-surjectivity}, and therefore $\tau_N'$ is in the image of the pullback
\begin{equation*}
H_{K_{\bf v_\ell}}^*(\Rep(Q, {\bf v_\ell}) \times \Gr({\bf v_u}, {\bf v_\ell})) \rightarrow H_{K_{\bf v_\ell}}^*(\mathcal{F}) .
\end{equation*}
Therefore the cup product 
\begin{align*}
H_{K_{\bf v_\ell}}^k(W_u, W_u \setminus \mathcal{F}_{\ell,0}^{u,0}) \times H_{K_{\bf v_\ell}}^d(W_u) & \rightarrow H_{K_{\bf v_\ell}}^{k+d}(W_u, W_u \setminus \mathcal{F}_{\ell,0}^{u,0}) \\
(\omega, \tau_N') & \mapsto \omega \smallsmile \tau_N'
\end{align*}
is the restriction of cup product
\begin{align*}
H_{K_{\bf v_\ell}}^k(M_{j+1} \times \Gr, M_{j-1} \times \Gr) \times H_{K_{\bf v_\ell}}^d(M_{j+1} \times \Gr) & \rightarrow H_{K_{\bf v_\ell}}^{k+d}(M_{j+1} \times \Gr, M_{j-1} \times \Gr) \\
(\tilde{\omega} , \tilde{\tau}) & \mapsto \tilde{\omega} \smallsmile \tilde{\tau} .
\end{align*}

Consider the Chern class $\xi$ from Remark \ref{rem:Chern-class}. Proposition \ref{prop:cup-product-pushforward} then shows that cup product with $\tau_N' \smallsmile \xi$ on the space of flow lines induces a homomorphism $H_K^p(M_j, M_{j-1}) \rightarrow H_K^{p+m}(M_{j+1}, M_j)$ between the terms on the first page of the Morse spectral sequence. 

\begin{equation}\label{eqn:cup-product-Thom}
\begin{tikzcd}[column sep=0.6cm]
H_K^p(M_j, M_{j-1}) \arrow[leftrightarrow]{d}{\cong}[swap]{\text{Thm.} \, \ref{thm:main-thm-morse}} \arrow{dr}{restriction} \arrow[dashrightarrow]{rrr} & & & H_K^{p+m}(M_{j+1}, M_j) \arrow[leftrightarrow]{d}{\cong}[swap]{\text{Thm.} \, \ref{thm:main-thm-morse}} \\
H_K^p(W_j, W_{j,0})  & H_K^p(W_{j+1,0}, W_{j+1,0} \setminus \mathcal{F}_{j,0}^{j+1,0}) \arrow{rr}{c(\tau_N' \smallsmile \xi) \, \text{(Lem. \ref{lem:cup-product-main-thm})}} & & H_K^{p+m}(W_{j+1}, W_{j+1,0})  \\
 & H_K^{p-\nu_j}(\mathcal{F}_{j,0}^{j+1,0}) \arrow[leftrightarrow]{u}{\text{Thom}}[swap]{\cong} \arrow{r}{\smallsmile \tau_N'} & H_K^{p+m-\nu_j}(\mathcal{F}_{j,0}^{j+1,0}) \arrow{r}{push}[swap]{forward} & H_K^{p+m-\lambda_{j+1}}(C_{j+1}) \arrow[leftrightarrow]{u}{\text{Thom}}[swap]{\cong} \\
H_K^{p-\lambda_j}(C_j) \arrow[leftrightarrow]{uu}{\text{Thom}}[swap]{\cong} \arrow{r}{pullback} & H_K^{p-\lambda_j}(\mathcal{F}_{j,0}^{j+1,0}) \arrow{u}{\smallsmile e}  & & 
\end{tikzcd}
\end{equation}

Now we can explain the topological meaning of this induced homomorphism. The above results show that the cup product homomorphism $c(\tau_N' \smallsmile \xi)$ is the restriction of the analogous homomorphism on $M \times \Gr$, since the classes $\tau_N'$ and $\xi$ are the pullback of the corresponding classes $\tilde{\tau}$ and $\tilde{\xi}$ via the inclusion $S_u \hookrightarrow M \times \Gr$ (Corollary \ref{cor:global-generation-Thom} and Lemma \ref{lem:global-generation-pushforward}). Therefore naturality of the cup product process of Lemma \ref{lem:cup-product-main-thm} implies that the following diagram commutes.
\begin{equation*}
\begin{tikzcd}[column sep = 2.5cm]
H_K^p(M_j \times \Gr, M_{j-1} \times \Gr) \arrow{d} \arrow{r}{c(\tilde{\tau} \smallsmile \tilde{\xi}) \, \text{(Lem. \ref{lem:cup-product-main-thm})}} & H_K^{p+m}(M_{j+1} \times \Gr, M_j \times \Gr)  \arrow{d} \\
H_K^p(W_{j+1,0}, W_{j+1,0} \setminus \mathcal{F}_{j,0}^{j+1,0}) \arrow{r}{c(\tau_N' \smallsmile \xi) \, \text{(Lem. \ref{lem:cup-product-main-thm})}} & H_K^{p+m}(W_{j+1}, W_{j+1,0}) 
\end{tikzcd}
\end{equation*}
Naturality of cup product then implies that the restriction $H_K^p(M_j, M_{j-1}) \rightarrow H_K^p(W_{j+1,0}, W_{j+1,0} \setminus \mathcal{F}_{j,0}^{j+1,0})$ factors through the pullback to the product $H_K^p(M_j, M_{j-1}) \rightarrow H_K^p(M_j \times \Gr, M_{j-1} \times \Gr)$ by the projection $M \times \Gr \rightarrow M$. Therefore the diagram \eqref{eqn:cup-product-Thom} can be augmented as follows

{\small 
\begin{equation}\label{eqn:cup-product-Thom-augmented}
\begin{tikzcd}[column sep=0.6cm]
& H_K^p(M_j \times \Gr, M_{j-1} \times \Gr) \arrow{dd} \arrow{rr}{c(\tilde{\tau} \smallsmile \tilde{\xi})} & & H_K^{p+m}(M_{j+1} \times \Gr, M_j \times \Gr) \arrow[bend right = 75]{dd} \\
H_K^p(M_j, M_{j-1}) \arrow{ur} \arrow[leftrightarrow]{d}{\cong}[swap]{\text{Thm.} \, \ref{thm:main-thm-morse}} \arrow{dr}{restriction} \arrow[crossing over, dashrightarrow]{rrr} & & & H_K^{p+m}(M_{j+1}, M_j) \arrow{u} \arrow[leftrightarrow]{d}{\cong}[swap]{\text{Thm.} \, \ref{thm:main-thm-morse}} \\
H_K^p(W_j, W_{j,0})  & H_K^p(W_{j+1,0}, W_{j+1,0} \setminus \mathcal{F}_{j,0}^{j+1,0}) \arrow{rr}{c(\tau_N' \smallsmile \xi) \, \text{(Lem. \ref{lem:cup-product-main-thm})}} & & H_K^{p+m}(W_{j+1}, W_{j+1,0})  \\
 & H_K^{p-\nu_j}(\mathcal{F}_{j,0}^{j+1,0}) \arrow[leftrightarrow]{u}{\text{Thom}}[swap]{\cong} \arrow{r}{\smallsmile \tau_N'} & H_K^{p+m-\nu_j}(\mathcal{F}_{j,0}^{j+1,0}) \arrow{r}{push}[swap]{forward} & H_K^{p+m-\lambda_{j+1}}(C_{j+1}) \arrow[leftrightarrow]{u}{\text{Thom}}[swap]{\cong} \\
H_K^{p-\lambda_j}(C_j) \arrow[leftrightarrow]{uu}{\text{Thom}}[swap]{\cong} \arrow{r}{pullback} & H_K^{p-\lambda_j}(\mathcal{F}_{j,0}^{j+1,0}) \arrow{u}{\smallsmile e}  & & 
\end{tikzcd}
\end{equation}
}
Therefore we see that the induced homomorphism $H_K^p(M_j \times \Gr, M_{j-1} \times \Gr) \rightarrow H_K^{p+m}(M_{j+1} \times \Gr, M_j \times \Gr)$ pulls back to the cup product homomorphism on the first page of the spectral sequence for $M \times \Gr$.

Now recall the convolution construction of \cite[Sec. 2(i)]{Nakajima97} and \cite[(2.7.14)]{ChrissGinzburg97} in Borel-Moore homology. Lemmas \ref{lem:PD-pullback} and \ref{lem:pushforward-projective-bundle} show that this corresponds via Poincar\'e duality to the composition of the following homomorphisms in cohomology
\begin{equation}\label{eqn:convolution-singular}
\begin{tikzcd}[column sep=1cm]
H_{K_{\bf v_\ell}}^p(C_\ell \cap \nu^{-1}(0)) \arrow{r}{\pi_\ell^*} & H_{K_{\bf v_\ell}}^p(\mathcal{B}) \arrow{r}{(\pi_u)_*} & H_{K_{\bf v_\ell}}^{p+d-\lambda}(C_{u} \cap \nu^{-1}(0)) .
\end{tikzcd}
\end{equation}

(Recall from \cite{Wilkin17} that the Hecke correspondence is $\mathcal{B} / K_{\bf v_\ell}$.)

Since the pushforward $H_{K_{\bf v_\ell}}^p(\mathcal{B}) \rightarrow  H_{K_{\bf v_\ell}}^{p+d-\lambda_u}(C_{u} \cap \nu^{-1}(0))$ factors through the inclusion $\mathcal{B} \hookrightarrow T$, then the convolution is given by the composition of the following homomorphisms
\begin{equation}\label{eqn:convolution}
\begin{tikzcd}[column sep=1cm]
H_{K_{\bf v_\ell}}^p(C_\ell \cap \nu^{-1}(0)) \arrow{r}{\pi_\ell^*} & H_{K_{\bf v_\ell}}^p(\mathcal{B}) \arrow{r}{i_*} & H_{K_{\bf v_\ell}}^{p+d}(T) \arrow{r}{(\pi_u)_*} & H_{K_{\bf v_u}}^{p+d-\lambda}(C_{u} \cap \nu^{-1}(0)) .
\end{tikzcd}
\end{equation} 
Note that the image of the pushforward 
\begin{equation*}
H_{K_{\bf v_\ell}}^{*+\lambda}(T) \stackrel{(\pi_u)_*}{\longrightarrow} H_{K_{\bf v_\ell}}^{*}(C_{u} \cap \nu^{-1}(0))  \cong H^{*}(\mathcal{M}(Q, {\bf v_u}) \otimes H^*(\BU(1)) .
\end{equation*}
lies in $H^{*}(\mathcal{M}(Q, {\bf v_u}) \hookrightarrow H^{*}(\mathcal{M}(Q, {\bf v_u}) \otimes H^*(\BU(1))$ (see Remark \ref{rem:image-pushforward}). Therefore we use $K_{\bf v_u}$-equivariant cohomology for the final term of the above diagram.

Corollary \ref{cor:thom-compatible} shows that cup product with the Thom class commutes with restriction
\begin{equation*}
\begin{tikzcd}
H_{K_{\bf v_\ell}}^p(\mathcal{F}) \arrow{r}{i^*} \arrow{d}{\smallsmile \tau_N} & H_{K_{\bf v_\ell}}^p(T) \arrow{d}{\smallsmile \tau_{\mathcal{B}}} \\
H_{K_{\bf v_\ell}}^{p+d}(\mathcal{F}, \mathcal{F} \setminus N) \arrow{r}{i^*} \arrow{d} & H_{K_{\bf v_\ell}}^{p+d}(T, T \setminus \mathcal{B}) \arrow{d} \\
H_{K_{\bf v_\ell}}^{p+d}(\mathcal{F}) \arrow{r}{i^*} & H_{K_{\bf v_\ell}}^{p+d}(T)  
\end{tikzcd}
\end{equation*}

The convolution \eqref{eqn:convolution} is the middle column of the following diagram. Taking the quotient by $K_{\bf v_\ell}$ gives the homomorphisms in the right hand column, which is Poincar\'e dual to the convolution used by Nakajima \cite[Sec. 2(i)]{Nakajima97} and Chriss \& Ginzburg \cite[(2.7.14)]{ChrissGinzburg97}. The left-hand column is the bottom row of the diagram \eqref{eqn:cup-product-Thom-augmented}, which therefore corresponds to the cup product on the Morse complex in the case where $C_\ell$ is a global minimum for $\| \mu -\alpha \|^2$.

{\small
\begin{equation}\label{eqn:full-convolution-cup-product-diagram}
\begin{tikzcd}[column sep=0.5cm]
& H_{K_{\bf v_\ell}}^p(C_\ell) \arrow{d}[swap]{\pi_\ell^*} \arrow{rrr}{i^*} & & & H_{K_{\bf v_\ell}}^p(C_\ell \cap \nu^{-1}(0)) \arrow{d}{\pi_\ell^*} \arrow{dl}[swap]{\pi_\ell^*} \arrow[leftrightarrow]{r}{\cong} & H^p(\mathcal{M}(Q, {\bf v_\ell}, \mathcal{R})) \arrow{d}{pullback} \\
& H_{K_{\bf v_\ell}}^p(\mathcal{F}) \arrow{d}[swap]{\smallsmile \tau_N} \arrow{rr}{i^*} &  & H_{K_{\bf v_\ell}}^p(T) \arrow{r}{i^*} \arrow{dr}[swap]{\smallsmile \tau_{\mathcal{B}}}& H_{K_{\bf v_\ell}}^p(\mathcal{B}) \arrow{d}{\cong}[swap]{Thom}   \arrow[leftrightarrow]{r}{\cong} & H^p(\mathcal{B}/K_{\bf v_\ell}) \arrow{ddd}{pushforward} \\
 & H_{K_{\bf v_\ell}}^{p+d}(\mathcal{F}, \mathcal{F} \setminus N) \arrow{d} \arrow{rrr}{i^*} & & & H_{K_{\bf v_\ell}}^{p+d}(T, T \setminus \mathcal{B}) \arrow{d} \\
 & H_{K_{\bf v_\ell}}^{p+d}(\mathcal{F}) \arrow{d}{(\pi_u)_*} \arrow{rrr}{i^*} & & & H_{K_{\bf v_\ell}}^{p+d}(T) \arrow{d}{(\pi_u)_*} \\
 & H_{K_{\bf v_u}}^{p+d-\lambda_{u}}(C_{u}) \arrow{rrr}{i^*} & & & H_{K_{\bf v_u}}^{p+d-\lambda_{u}}(C_{u} \cap \nu^{-1}(0)) \arrow[leftrightarrow]{r}{\cong} & H^{p+d-\lambda_u}(\mathcal{M}(Q, {\bf v_u}, \mathcal{R}))
\end{tikzcd}
\end{equation}
}

When ${\bf v_\ell} < {\bf v}$ then there is an additional step of cup product with the Euler class of Corollary \ref{cor:Euler-class} that relates the above construction to the cup product of \eqref{eqn:cup-product-Thom-augmented}.

Combined with the results of the previous section, we have proved the following.

\begin{theorem}\label{thm:cup-product-induces-convolution}
The cup product $H_K^p(M_j \times \Gr, M_{j-1} \times \Gr) \stackrel{\smallsmile (\tilde{\tau} \smallsmile \tilde{\xi})}{\longrightarrow} H_K^{p+d}(M_{\ell+1} \times \Gr, M_\ell \times \Gr)$ on the Morse complex for $M \times \Gr$ restricts to the local cup product with $\tau_N' \smallsmile \xi$ in the Morse complex for $M$ (the left-hand column of \eqref{eqn:full-convolution-cup-product-diagram}. Restricting this to $\nu^{-1}(0)$ determines the middle column of \eqref{eqn:full-convolution-cup-product-diagram}, which is then equivalent to convolution \eqref{eqn:convolution}.
\end{theorem}


\begin{thebibliography}{10}

\bibitem{ACGH85}
E.~Arbarello, M.~Cornalba, P.~A. Griffiths, and J.~Harris.
\newblock {\em Geometry of algebraic curves. {V}ol. {I}}, volume 267 of {\em
  Grundlehren der Mathematischen Wissenschaften [Fundamental Principles of
  Mathematical Sciences]}.
\newblock Springer-Verlag, New York, 1985.

\bibitem{AtiyahBott83}
M.~F. Atiyah and R.~Bott.
\newblock The {Y}ang-{M}ills equations over {R}iemann surfaces.
\newblock {\em Philos. Trans. Roy. Soc. London Ser. A}, 308(1505):523--615,
  1983.

\bibitem{AtiyahBott84}
M.~F. Atiyah and R.~Bott.
\newblock The moment map and equivariant cohomology.
\newblock {\em Topology}, 23(1):1--28, 1984.

\bibitem{AustinBraam95}
D.~M. Austin and P.~J. Braam.
\newblock Morse-{B}ott theory and equivariant cohomology.
\newblock In {\em The {F}loer memorial volume}, volume 133 of {\em Progr.
  Math.}, pages 123--183. Birkh\"auser, Basel, 1995.

\bibitem{Bott54}
Raoul Bott.
\newblock Nondegenerate critical manifolds.
\newblock {\em Ann. of Math. (2)}, 60:248--261, 1954.

\bibitem{Bott88}
Raoul Bott.
\newblock Morse theory indomitable.
\newblock {\em Inst. Hautes \'{E}tudes Sci. Publ. Math.}, (68):99--114 (1989),
  1988.

\bibitem{BottTu82}
Raoul Bott and Loring~W. Tu.
\newblock {\em Differential forms in algebraic topology}, volume~82 of {\em
  Graduate Texts in Mathematics}.
\newblock Springer-Verlag, New York-Berlin, 1982.

\bibitem{ChrissGinzburg97}
Neil Chriss and Victor Ginzburg.
\newblock {\em Representation theory and complex geometry}.
\newblock Birkh\"auser Boston Inc., Boston, MA, 1997.

\bibitem{Crawley01}
William Crawley-Boevey.
\newblock Geometry of the moment map for representations of quivers.
\newblock {\em Compositio Math.}, 126(3):257--293, 2001.

\bibitem{Hitchin87}
N.~J. Hitchin.
\newblock The self-duality equations on a {R}iemann surface.
\newblock {\em Proc. London Math. Soc. (3)}, 55(1):59--126, 1987.

\bibitem{King94}
A.~D. King.
\newblock Moduli of representations of finite-dimensional algebras.
\newblock {\em Quart. J. Math. Oxford Ser. (2)}, 45(180):515--530, 1994.

\bibitem{Kirwan84}
Frances~Clare Kirwan.
\newblock {\em Cohomology of quotients in symplectic and algebraic geometry},
  volume~31 of {\em Mathematical Notes}.
\newblock Princeton University Press, Princeton, NJ, 1984.

\bibitem{Kronheimer89}
P.~B. Kronheimer.
\newblock The construction of {ALE} spaces as hyper-{K}\"ahler quotients.
\newblock {\em J. Differential Geom.}, 29(3):665--683, 1989.

\bibitem{KronheimerNakajima90}
Peter~B. Kronheimer and Hiraku Nakajima.
\newblock Yang-{M}ills instantons on {ALE} gravitational instantons.
\newblock {\em Math. Ann.}, 288(2):263--307, 1990.

\bibitem{Lojasiewicz84}
S.~{\L}ojasiewicz.
\newblock Sur les trajectoires du gradient d'une fonction analytique.
\newblock In {\em Geometry seminars, 1982--1983 (Bologna, 1982/1983)}, pages
  115--117. Univ. Stud. Bologna, Bologna, 1984.

\bibitem{McGertyNevins18}
Kevin McGerty and Thomas Nevins.
\newblock Kirwan surjectivity for quiver varieties.
\newblock {\em Invent. Math.}, 212(1):161--187, 2018.

\bibitem{Nakajima94}
Hiraku Nakajima.
\newblock Instantons on {ALE} spaces, quiver varieties, and {K}ac-{M}oody
  algebras.
\newblock {\em Duke Math. J.}, 76(2):365--416, 1994.

\bibitem{Nakajima97}
Hiraku Nakajima.
\newblock Heisenberg algebra and {H}ilbert schemes of points on projective
  surfaces.
\newblock {\em Ann. of Math. (2)}, 145(2):379--388, 1997.

\bibitem{Nakajima98}
Hiraku Nakajima.
\newblock Quiver varieties and {K}ac-{M}oody algebras.
\newblock {\em Duke Math. J.}, 91(3):515--560, 1998.

\bibitem{Nakajima01}
Hiraku Nakajima.
\newblock Quiver varieties and finite-dimensional representations of quantum
  affine algebras.
\newblock {\em J. Amer. Math. Soc.}, 14(1):145--238 (electronic), 2001.

\bibitem{Nakajima04}
Hiraku Nakajima.
\newblock Quiver varieties and {$t$}-analogs of {$q$}-characters of quantum
  affine algebras.
\newblock {\em Ann. of Math. (2)}, 160(3):1057--1097, 2004.

\bibitem{Nakajima12}
Hiraku Nakajima.
\newblock Handsaw quiver varieties and finite {$W$}-algebras.
\newblock {\em Mosc. Math. J.}, 12(3):633--666, 669--670, 2012.

\bibitem{PflaumWilkin19}
Markus~J. Pflaum and Graeme Wilkin.
\newblock Equivariant control data and neighborhood deformation retractions.
\newblock {\em Methods Appl. Anal.}, 26(1):13--36, 2019.

\bibitem{Reineke03}
Markus Reineke.
\newblock The {H}arder-{N}arasimhan system in quantum groups and cohomology of
  quiver moduli.
\newblock {\em Invent. Math.}, 152(2):349--368, 2003.

\bibitem{Simon83}
Leon Simon.
\newblock Asymptotics for a class of nonlinear evolution equations, with
  applications to geometric problems.
\newblock {\em Ann. of Math. (2)}, 118(3):525--571, 1983.

\bibitem{Sjamaar98}
Reyer Sjamaar.
\newblock Convexity properties of the moment mapping re-examined.
\newblock {\em Adv. Math.}, 138(1):46--91, 1998.

\bibitem{Smale60}
Stephen Smale.
\newblock Morse inequalities for a dynamical system.
\newblock {\em Bull. Amer. Math. Soc.}, 66:43--49, 1960.

\bibitem{wilkin-condition-4}
Graeme Wilkin.
\newblock Local behaviour of the gradient flow of an analytic function near the
  unstable set of a critical point.
\newblock Available at \texttt{https://arxiv.org/abs/1904.08045}.

\bibitem{Wilkin17}
Graeme Wilkin.
\newblock Moment map flows and the {H}ecke correspondence for quivers.
\newblock {\em Adv. Math.}, 320:730--794, 2017.

\bibitem{Wilkin19}
Graeme Wilkin.
\newblock Equivariant {M}orse theory for the norm-square of a moment map on a
  variety.
\newblock {\em Int. Math. Res. Not. IMRN}, (15):4730--4763, 2019.

\bibitem{Witten82}
Edward Witten.
\newblock Supersymmetry and {M}orse theory.
\newblock {\em J. Differential Geometry}, 17(4):661--692 (1983), 1982.

\end{thebibliography}

\end{document}